\documentclass[11pt,a4paper,english,reqno
]{amsart}
\usepackage{verbatim}
\usepackage{amsthm}
\usepackage{amscd}
\usepackage{amstext}
\usepackage{amssymb}
\usepackage[all]{xy}

\usepackage[dvips]{graphicx}
\usepackage[colorinlistoftodos]{todonotes}
\usepackage[colorlinks=true,linkcolor=red,citecolor=blue]{hyperref}
\usepackage{tabu}
\usepackage{mathdots}

\usepackage{scalerel}

\usepackage[colorinlistoftodos]{todonotes}
\usepackage{calrsfs}
\usepackage{tikz}
\usetikzlibrary{positioning,intersections,cd,matrix, arrows, backgrounds}
\usetikzlibrary{shapes,shapes.geometric,shapes.misc}
\usetikzlibrary{decorations}
\usetikzlibrary{decorations.pathreplacing}
\newcommand\Ar[3]{\ar[from={#1}, to={#2}, #3]}

\makeatletter
\@namedef{subjclassname@2010}{%
	\textup{2010} Mathematics Subject Classification}
\makeatother

\numberwithin{equation}{section}
\numberwithin{figure}{section}
\theoremstyle{plain}

\newtheorem{thm-nonum}{Theorem}

\newtheorem*{thm-nn}{Theorem}
\newtheorem*{prp-nn}{Proposition}

\usepackage{amscd, amssymb,}
\xyoption{2cell}
\UseTwocells

\newcommand{\xyR}[1]{%
	\xydef@\xymatrixrowsep@{#1}}

\newcommand{\xyC}[1]{%
	\xydef@\xymatrixcolsep@{#1}}

\newcommand{\rt}{\rightarrow}

\newcommand{\st}{\stackrel}

\newcommand{\al}{\alpha}
\newcommand{\la}{\lambda}
\newcommand{\La}{\Lambda}

\newcommand{\bbN}{\mathbb{N}}

\newcommand{\CA}{\mathcal{A} }
\newcommand{\CC}{\mathcal{C} }
\newcommand{\DD}{\mathcal{D} }

\newcommand{\CF}{\mathcal{F} }
\newcommand{\CG}{\mathcal{G} }

\newcommand{\CK}{\mathcal{K} }

\newcommand{\CM}{\mathcal{M} }
\newcommand{\CN}{\mathcal{N} }

\newcommand{\CX}{\mathcal{X} }

\newcommand{\CU}{\mathcal{U}}

\newcommand{\CB}{\mathcal{B} }

\newcommand{\Mod}{{\mathrm{Mod\mbox{-}}}}

\newcommand{\mmod}{{\mathrm{mod}}\mbox{-}}

\newcommand{\prj}{{\mathrm{prj}\mbox{-}}}

\newcommand{\Loc}{{\mathrm{Loc}}}

\newcommand{\op}{{\mathrm{op}}}

\newcommand{\add}{{\mathrm{add}\mbox{-}}}

\newcommand{\Add}{{\mathrm{Add}\mbox{-}}}

\newcommand{\id}{{\mathrm{id}}}

\newcommand{\Coker}{{\mathrm{Coker}}}
\newcommand{\Ker}{{\mathrm{Ker}}}

\newcommand{\Hom}{{\mathrm{Hom}}}
\newcommand{\Ext}{{\mathrm{Ext}}}

\newcommand\Cdot{\centerdot}
\newcommand\down{_{\Cdot}}
\newcommand\up{^{\Cdot}}
\newcommand\ya[1]{\xrightarrow{#1}}
\newcommand\ang[1]{{\langle #1 \rangle}}

\newcommand{\vertmap}[3]{%
{\footnotesize
\begin{tikzcd}
#1\\#3
\Ar{1-1}{2-1}{"#2"}
\end{tikzcd}}}
\newcommand{\longto}[2]{\xrightarrow[\ \ #2\ \ ]{\ \ #1\ \ }}
\newcommand\blank{\mbox{-}}
\newcommand\umod{\underline{\operatorname{mod}}\blank}
\newcommand\uMod{\underline{\operatorname{Mod}}\blank}
\newcommand\sfC{\mathsf{C}}
\newcommand\sfD{\mathsf{D}}

\newcommand\calA{\mathcal{A}}
\newcommand\calB{\mathcal{B}}
\newcommand\calC{\mathcal{C}}
\newcommand\calD{\mathcal{D}}
\newcommand\calF{\mathcal{F}}

\newcommand\calP{\mathcal{P}}
\newcommand\calS{\mathcal{S}}
\newcommand\calT{\mathcal{T}}

\newcommand\calK{\mathcal{K}}
\newcommand\frkU{\mathfrak{U}}
\newcommand\rmH{\mathrm{H}}

\renewcommand\al{\alpha}
\newcommand\be{\beta}

\newcommand\de{\delta}
\newcommand\ep{\varepsilon}

\newcommand\et{\eta}
\renewcommand\th{\theta}

\renewcommand\la{\lambda}
\newcommand\ro{\rho}
\newcommand\si{\sigma}

\newcommand\ph{\phi}

\renewcommand\La{\Lambda}
\newcommand\om{\omega}

\newcommand\Th{\Theta}

\newcommand\Aut{\operatorname{Aut}}
\newcommand\supp{\operatorname{supp}}
\newcommand\dom{\operatorname{dom}}
\newcommand\cod{\operatorname{cod}}
\newcommand\cpt{\operatorname{cpt}}
\newcommand\obj{\operatorname{obj}}

\newcommand\ind{\mathrm{ind\mbox{-}}}
\newcommand\udl{\underline}
\newcommand\ovl{\overline}
\newcommand\Ds{\bigoplus}
\newcommand\inv{^{-1}}
\newcommand\To{\Rightarrow}
\newcommand\idty{1\kern-.25em{\text{{\rm l}}}} 
\newcommand\ds{\oplus}

\newcommand\kCat{k\text{-}\mathbf{Cat}}
\newcommand\GCat{G\text{-}\mathbf{Cat}}
\newcommand\kCAT{k\text{-}\mathbf{CAT}}
\newcommand\GCAT{G\text{-}\mathbf{CAT}}

\newcommand\adj{\dashv}

\renewcommand\CF{\operatorname{fp}}

\theoremstyle{plain}
\newtheorem{theorem}{Theorem}[section]

\newtheorem{lemma}[theorem]{Lemma}

\newtheorem{proposition}[theorem]{Proposition}

\newtheorem{notation}[theorem]{Notation}

\newtheorem{clm}{Claim}

\newtheorem*{clm-nn}{Claim}

\theoremstyle{definition}
\newtheorem{definition}[theorem]{Definition}
\newtheorem{example}[theorem]{Example}

\newtheorem{construction}[theorem]{Construction}

\newtheorem{remark}[theorem]{Remark}

\def\repr[#1;#2;#3;#4;#5]{
	\left(
	\begin{matrix}#1\\#2\end{matrix}
	#3
	\begin{matrix}#4\\#5\end{matrix}
	\right)}
\def\bmat#1{\begin{bmatrix} #1 \end{bmatrix}}

\numberwithin{equation}{section}
\numberwithin{figure}{section}

\newcommand\reallywidehat[1]{\arraycolsep=0pt\relax%
	\begin{array}{c}
		\stretchto{
			\scaleto{
				\scalerel*[\widthof{\ensuremath{#1}}]{\kern-.5pt\bigwedge\kern-.5pt}
				{\rule[-\textheight/2]{1ex}{\textheight}} 
			}{\textheight} %
		}{0.7ex}\\           
		#1\\                 
		\rule{-1ex}{0ex}
	\end{array}}

\newcommand{\pr}{\operatorname{pr}}
\newcommand{\DS}{\bigoplus}

\newcommand\Nname[1]{|[alias=#1]|}

\newcommand{\sk}{\operatorname{sk}}

\newcommand{\Kb}{\operatorname{\calK^{\mathrm{b}}}}

\begin{document}

\title[Unifying constructions of precoverings]%
{
2-categorical approach to unifying constructions of precoverings and  its applications}


\author[R.~Hafezi, H.~Asashiba and  M.H.~Keshavarz ]%
{Rasool Hafezi, Hideto Asashiba and   Mohammad Hossein Keshavarz}

\thanks{HA is partially supported by Grant-in-Aid for Scientific Research
 (C) 18K03207 from JSPS, and by Osaka Central Advanced Mathematical Institute
(MEXT Joint Usage/Research Center on Mathematics and Theoretical Physics JPMXP0619217849).}

\address{School of Mathematics and Statistics,
Nanjing University of Information Science \& Technology, Nanjing, Jiangsu 210044, P.\,R. China}
\email{hafezi@nuist.edu.cn}

\address{Department of Mathematics, Faculty of Science, Shizuoka University,
836 Ohya, Suruga-ku, Shizuoka, 422-8529, Japan;}
\address{
Institute for Advanced Study, KUIAS, Kyoto University,
Yoshida Ushinomiya-cho, Sakyo-ku, Kyoto 606-8501, 
Japan; and}
\address{Osaka Central Advanced Mathematical Institute,
3-3-138 Sugimoto, Sumiyoshi-ku,
Osaka, 558-8585, Japan.}
\email{asashiba.hideto@shizuoka.ac.jp}

\address{School of Mathematical Sciences, Shanghai Key Laboratory of PMMP, East China Normal University,  Shanghai 200241, P. R. China}
\email{keshavarz@ipm.ir}

\subjclass[2010]{18A25, 16G70, 16G10}

\keywords{Galois  (pre)covering;  morphism category;  locally bounded category; }

\begin{abstract}
Throughout this paper $G$ is a fixed group, and $k$ is a fixed field.  All categories are assumed to be $k$-linear.
First we give a systematic way to induce $G$-precoverings by adjoint functors using a 2-categorical machinery, which unifies many similar constructions of
$G$-precoverings.

Now let $\mathcal{C}$ be a skeletally small category with a $G$-action, $\mathcal{C}/G$ the orbit category of $\mathcal{C}$,
$(P, \phi) : \mathcal{C} \rightarrow \mathcal{C}/G$ the canonical $G$-covering, and $\mmod\CC$, $\mmod (\mathcal{C}/G)$ the categories
of finitely generated modules over $\mathcal{C}, \mathcal{C}/G$, respectively.
Then it is well known that there exists a canonical G-precovering $(P., \phi.) : \mathrm{mod}\mbox{-} \mathcal{C} \rightarrow \mmod(\mathcal{C}/G)$.
By applying the machinery above to this $(P., \phi.)$, new $G$-precoverings $(\mathrm{mod}\mbox{-} \mathcal{C}) / S \rightarrow (\mathrm{mod}\mbox{-} \mathcal{C}/G)/S'$ are induced
between the factor categories or localizations of $\mathrm{mod}\mbox{-} \mathcal{C}$ and $\mathrm{mod}\mbox{-} \mathcal{C}/G$, respectively.

This is further applied to the morphism category $\mathrm{H}(\mathrm{mod}\mbox{-}  \mathcal{C})$ of $\mathrm{mod}\mbox{-}  \mathcal{C}$ to have a $G$-precovering $\mathrm{fp}(\CK) \rightarrow \mathrm{fp}(\CK')$
between the categories of finitely presented modules over
suitable subcategories $\CK$ and $\CK'$ of 
$\mmod\CC$ and $ \mathrm{mod}\mbox{-} \mathcal{C}/G$, respectively.
\end{abstract} 

\maketitle

\tableofcontents

\section{Introduction}

Throughout this paper $G$ is a group, and $k$ is a field.
All categories are assumed to be $k$-linear.
A $G$-category is a $k$-category with a $G$-action.
We denote by $\kCAT$ (resp.\ $\GCAT$) the 2-category of light $k$-categories
(resp.\ $G$-categories),
where a category $\CC$ is said to be \emph{light}
if the objects form a subset of a fixed universe $\frkU$
and $\CC(x,y)$ is an element of $\frkU$ for all $x, y \in \CC$
(see the beginning of Sect.\ 2, and Definition \ref{dfn:GCat} for details).

\subsection{Examples of constructions of $G$-precoverings}

The notion of a Galois covering functor with group $G$ is introduced by Gabriel \cite{G}, which makes problems of the module category $\mmod A$ of a finite-dimensional algebra $A$ defined by a bound quiver $(Q,I)$
easier as follows.
There is a way to compute a locally finite-dimensional category $\calC$ with an action of a group $G$ such that the path-category $k[Q,I]$ of $(Q,I)$ is equivalent to the orbit category $\calC/G$
(note that $\mmod A$ is equivalent to $\mmod k[Q,I]$,
and hence to $\mmod(\calC/G)$), and in many cases
$\mmod \calC$ is easier to study than $\mmod A$,
and it is possible to control $\mmod(\calC/G) \simeq \mmod A$ by
$\mmod \calC$. 
Hence problems of $\mmod A$ can be reduced to those
of $\mmod \calC$ (see \cite{BG, Ri} for example).
Now, each Galois covering functor with group $G$ is isomorphic to the canonical functor $P : \calC \to \calC/G$ from a locally finite-dimensional category $\CC$ with a $G$-action to the orbit category $\CC/G$ of $\CC$ by $G$. 
In special cases, he induced from $P$ a Galois covering $P\down \colon \ind \CC \to \ind(\CC/G)$ with group $G$ from a full subcategory
$\ind \CC$ of the category $\mmod \CC$ of finitely generated
$\CC$-modules whose object set forms a $G$-stable complete set of representatives of isomorphism classes of indecomposable $\CC$-modules to a full subcategory $\ind (\CC/G)$
of the category $\mmod(\CC/G)$ of finitely generated
$\CC/G$-modules consisting of complete representatives of
isomorphism classes of indecomposable modules.

In \cite{A}, one of the authors extended the setting of covering functors
from locally finite-dimensional categories to any linear categories
by using a family $\phi$ of natural isomorphisms,
and introduced the notion of $G$-coverings and
$G$-precoverings $(P, \phi)\colon \calC \to \calC/G$ with a generalized version of the orbit category $\calC/G$ of a $G$-category $\calC$.
By the same idea as in \cite{G},
from a $G$-precovering $(P, \ph)$,
he also induced
the $G$-precoverings
$(P\down, \ph\down) \colon \mmod \calC \to \mmod(\calC/G)$
and 
$(P\down, \ph\down) \colon \Kb(\prj \calC) \to \Kb(\prj(\calC/G))$.
There are other constructions of $G$-precoverings
such as 
\begin{itemize}
\item
in \cite{As} (1997) for bounded homotopy categories of
finitely generated projective modules,
\item
in \cite{BL} (2014) for bounded derived categories,
\item
in \cite{AHV} (2018) for bounded derived categories,
singularity categories,
Gorenstein defect categories, and
stable Gorenstein projectives,
\item
in \cite{P} (2019) for finitely presented functors,
\item
in unpublised paper \cite{HM} (2020) for stable modules categories, and
\item
in \cite{HZhao} and \cite{HZou} (2023) for factor categories of modules categories.
\end{itemize}
See section \ref{ss:relations-other-works} for more details.

\subsection{Unification}

In this paper we will unify these constructions of $G$-precoverings by using a 2-categorical approach,
and give some further useful tools to inducing $G$-precoverings
as stated in Theorems \ref{thm:factor-case} and \ref{thm:general-case}.
To state the main theorems, we need the following
notion that generalizes Gabriel's idea to induce $G$-precoverings,
where $\cpt(\CA)$ is the full subcategory of a category $\CA$ consisting of the compact objects:

\newtheorem*{dfn-intro}{Definition}

\begin{definition}
[see Definition \ref{dfn:precov-adj}]
Let $\calA = (\calA, A)$ be a $G$-category having small coproducts, $\calB$ a $k$-category
and $(L, \ph) \colon \calA \to \calB$ a $G$-invariant functor, where $L$
has a right adjoint $R \colon \calB \to \calA$
with an adjunction isomorphism $\om \colon L \adj R$.
Then $(L, \ph)$ is said {\em to induce a $G$-precovering}
$(L, \ph) \colon \cpt(\calA) \to \cpt(\calB)$ {\em by the adjunction} $\om$
if the following conditions are satisfied:
\begin{enumerate}
\item[{\rm (PA1)}]
There exists a natural isomorphism $t \colon \Ds_{a\in G}A_a \To RL$.
\item[{\rm (PA2)}]
For any $x,y \in \calA$, the following diagram commutes:
\begin{equation}
\begin{tikzcd}
\Ds_{a \in G} \calA(x, ay) & \calA(x, \Ds_{a \in G} ay)\\
\calB(Lx, Ly) & \calA(x, RLy)
\Ar{1-1}{1-2}{"\nu", }
\Ar{1-2}{2-2}{"{\calA(x, t_y)}", "\rotatebox{90}{$\sim$}"'}
\Ar{1-1}{2-1}{"{(L, \ph)^{(2)}_{x,y}}"'}
\Ar{2-1}{2-2}{"\sim","\om_{x, Ly}"'}
\end{tikzcd},
\end{equation}
where $\nu$ is the canonical morphism.
Note that $\nu$ becomes an isomorphism if $x \in \cpt(\calA)$.
\item[{\rm (PA3)}]
$R$ preserves small coproducts.
\end{enumerate}
If in the above, the condition (PA3) is {\em not assumed}, then
$(L, \ph)$ is said {\em to induce a $G$-precovering}
$(L, \ph) \colon \cpt(\calA) \to \calB$
{\em by the adjunction} $\om$.
\end{definition}

Using this definition, for example,
the first theorem is stated as follows.

\newtheorem*{thm-A}{Theorem A}
\newtheorem*{prp-A}{Proposition A}
\newtheorem*{prp-B}{Proposition B}
\newtheorem*{prp-C}{Proposition C}
\newtheorem*{prp-D}{Proposition D}

\begin{theorem}
[see Theorem \ref{thm:factor-case}]
Consider the following diagram of $2$-categories, $2$-functors,
and a strict $2$-natural transformation on the left, and its extension
by Lemma \ref{lem:ext-2-nat} on the right:
\[
\begin{tikzcd}
\sfC & \makebox[1em][l]{$\kCAT$,}
\Ar{1-1}{1-2}{"\si", ""'{name=a}, bend left}
\Ar{1-1}{1-2}{"U"', ""{name=b}, bend right}
\Ar{a}{b}{"\pi", Rightarrow}
\end{tikzcd}\qquad\qquad
\begin{tikzcd}
\hat{\sfC} & \makebox[1em][l]{$\GCAT$,}
\Ar{1-1}{1-2}{"\hat{\si}", ""'{name=a}, bend left}
\Ar{1-1}{1-2}{"\hat{U}"', ""{name=b}, bend right}
\Ar{a}{b}{"\hat{\pi}", Rightarrow}
\end{tikzcd}
\]
where $\sfC$ is a $2$-subcategory of $\kCAT$ and $\si$ is the inclusion $2$-functor.
Assume the following conditions for all $\calC \in \sfC$:
\begin{enumerate}
\item[{\rm (0)}]
$\calC$ has small coproducts, and
$\pi_{\calC} \colon \calC \to U(\calC)$ preserves small coproducts.
%
\item[{\rm (i)}]
$\pi_{\calC}$ is the identity on the objects.
\item[{\rm (ii)}]
$\pi_{\calC}$ is epic among
representable right $U(\calC)$-modules
(see Definition \ref{dfn:epic-among}).
\end{enumerate}
Now let $(\calA, A)$ be a $G$-category in $\hat{\sfC}$,
$\calB$ a $k$-category in $\sfC$,
and $(L, \ph) \colon (\calA, A) \to \calB$ a $G$-invariant functor in $\hat{\sfC}$
 (cf. Remark 1.6),
where $L$ has a right adjoint $R \colon \calB \to \calA$ in $\sfC$
with an adjunction isomorphism $\om \colon L \adj R$ in $\sfC$.
Then the following hold.
\begin{enumerate}
\item
Assume that $(L, \ph) \colon \calA \to \calB$ induces a $G$-precovering $\cpt(\calA) \to \calB$ by the adjunction $\om$.
Then
\begin{enumerate}
\item
the $G$-invariant functor
$U(L, \ph) \colon U(\calA) \to U(\calB)$ induces a $G$-precovering
\[
U(L,\ph) \colon \cpt(U(\calA)) \to U(\calB)
\]
by the adjunction $U(\om) \colon U(L) \adj U(R)$, and
\item

if further $\cpt(\calA) = \pi_\calA(\cpt(\calA)) \subseteq \cpt(U(\calA))$
(this is the case, for example, 
if $\pi_{\calA}$ is 
full by Lemma \ref{lem:cpt-cpt-mod}),
then
$U(L, \ph)$ restricts to a $G$-precovering
\[
U(L,\ph) \colon U(\cpt(\calA)) \to U(\calB).
\]
\end{enumerate}

\item
Assume that $(L, \ph) \colon \calA \to \calB$ induces a $G$-precovering $\cpt(\calA) \to \cpt(\calB)$
by the adjunction $\om$.
Then
\begin{enumerate}
\item
the $G$-invariant functor
$U(L, \ph) \colon U(\calA) \to U(\calB)$ induces a $G$-precovering
\[
U(L,\ph) \colon \cpt(U(\calA)) \to \cpt(U(\calB))
\]
by the adjunction $U(\om) \colon U(L) \adj U(R)$, and
\item
if further $\cpt(\calA) \subseteq \cpt(U(\calA))$, then
$U(L, \ph)$ restricts to a $G$-precovering
\[
U(L,\ph) \colon U(\cpt(\calA)) \to U(\cpt(\calB)).
\]
\end{enumerate}
\end{enumerate}
\end{theorem}

These theorems are applied to constructions of 
stable module categories, or more generally,
factor categories (Proposition \ref{Prp:factor-covering}), and localizations (Proposition \ref{Prp:localization-covering}) as in
the following statements.

\begin{theorem}
[see Proposition \ref{Prp:factor-covering}]
Let $\calC$ be a skeletally small $G$-category, and 
$\calD$ a $G$-stable class of objects in $\Mod \calC$
closed under small coproducts with the property that
each object in $\calD$ is a small coproduct of finitely generated objects.
Denote by $S/G$ the ideal $\ang{P\down(\cpt(\calD))}$ of $\mmod (\calC/G)$. 
Then the following hold.
\begin{enumerate}
\item
The functor
\[
U(P\down, \ph\down) \colon (\mmod\CC)/S \to  (\mmod\CC/G)/(S/G)
\]
induced by the quotient $2$-functor $U = U_i$ in {\rm Lemma \ref{lem:factor-case}}
is a $G$-precovering.
\item
If $\CC$ is a locally support-finite locally bounded category,
then
\[
U(P\down, \ph\down) \colon (\mmod\CC)/S \to  (\mmod\CC/G)/(S/G)
\]
is a $G$-covering.
\end{enumerate}
\end{theorem}

\begin{theorem}
[see Proposition \ref{Prp:localization-covering}]
Let $\calC$ be a skeletally small $G$-category such that
$\Mod \calC$ is locally noetherian, 
both $\mmod \calC$ and $\mmod(\calC/G)$ are abelian,
and $\calS'$ a $G$-stable localizing subcategory of $\Mod\calC$.
Denote by $\calS$, $\calS/G$ the dense subcategories of $\mmod\calC$, $\mmod (\calC/G)$
induced from $\calS'$ as in Lemma \ref{lem:quotient-case}, respectively.
Then the following hold.
\begin{enumerate}
\item
The functor
\[
U(P\down, \ph\down) \colon (\mmod\CC)/\calS \to  (\mmod\CC/G)/(\calS/G)
\]
induced by the quotient $2$-functor $U = U_m$ in {\rm Example \ref{exm:quot-2-fun}}
is a $G$-precovering.
\item
If $\CC$ is a locally support-finite locally bounded category,
then
\[
U(P\down, \ph\down) \colon (\mmod\CC)/\calS \to  (\mmod\CC/G)/(\calS/G)
\]
is a $G$-covering.
\end{enumerate}
\end{theorem}

\subsection{Applications}

We apply these tools to the morphism category
$\rmH(\Mod \CC)$
of the module category $\Mod \CC$
over a skeletally small category $\CC$
to have the following theorem:

\begin{theorem}[see Proposition \ref{Prop-pushdown-adj}]
The $G$-invariant functor
$\rmH(P\down, \ph\down):\rmH(\Mod \CC) $ $\to \rmH(\Mod (\CC/G))$ induces a $G$-precovering 
$$\rmH(P\down, \ph\down):\rmH(\mmod \CC)\to \rmH(\mmod (\CC/G))$$
by the adjunction $\rmH(\theta)$, where
$\th \colon P\down \adj P\up$ is the adjunction
between the pushdown $P\down$ and the pullup $P\up$ (see Subsection \ref{Pushdownfunctor} for details).
\end{theorem}

For a full subcategory $\CK$ of an additive category $\CC$,
we denote by $\add\CK$ (resp.\ $\Add\CK$),  the full subcategory of $\CC$ consisting of the direct summands of finite (resp.\ small) coproducts of objects in $\CK$.
The theorem above is generalized to Proposition \ref{Proposition3last}
on subcategories $\CK$ and $\CK'$ of $\mmod \CC$
and $\mmod (\CC/G)$, respectively satisfying that
$P\down$ sends $\CK$ into $\CK'$ and
$P\up$ sends $\CK'$ into $\Add \CK$.

\newcommand{\Fp}{\mathrm{Fp}}

Now since the category $\CF(\CM)$ of finitely presented functors from a linear category $\CM$ to $\Mod k$
is shown to be equivalent to a factor category
$\frac{\rmH(\CM)}{\ang{\CU_\CM}}$
of the morphism category $\rmH(\CM)$ of $\CM$
by a suitable ideal $\ang{\CU_\CM}$,
the functor $\rmH(P\down, \ph\down)$ above yields $G$-precoverings
$\Fp(P\down) \colon \CF(\Mod \CC) \to \CF(\Mod(\CC/G))$
and
$\CF(P\down):\CF(\mmod \CC)\rt \CF(\mmod (\CC/G))$
as in the following theorem, where in the statement (2)
the vertical line as in $P\down|$ denotes the 
restriction of $P\down$ to a suitable subcategory.
\begin{theorem}[see Theorem \ref{thm:G-Fp-precovering}]
The following statements hold. 
\begin{enumerate}
\item
The functor  ${\rm Fp}(P\down)$ is  a $G$-precovering,
and has a right adjoint ${\rm Fp}(P\up)$
with an adjucntion $\Theta \colon {\rm Fp}(P\down) \adj {\rm Fp}(P\up)$.

Moreover, there is the following commutative diagram
\begin{equation}
\begin{tikzcd}
\CF(\mmod \CC) & \CF(\Mod \CC)\\
\CF(\mmod (\CC/G))& \CF(\Mod (\CC/G))
\Ar{1-1}{1-2}{"", }
\Ar{1-2}{2-2}{"{\rm Fp}(P\down)", ""'}
\Ar{1-1}{2-1}{"{{\rm fp}(P\down)}"'}
\Ar{2-1}{2-2}{"",""'}
\end{tikzcd},
\end{equation}
where the horizontal functors are embeddings.
\item 
The functor  ${\rm Fp}(P\down |)$ is  a $G$-precovering,
and has a right adjoint ${\rm Fp}(P\up|)$
with an adjucntion $\Theta \colon {\rm Fp}(P\down|) \adj {\rm Fp}(P\up|)$. 
Moreover, there is the following commutative diagram
\begin{equation}
\begin{tikzcd}
\CF(\CK) & \CF({\rm Add}\mbox{-}\CK)\\
\CF(\CK')& \CF({\rm Add}\mbox{-}\CK')
\Ar{1-1}{1-2}{"", }
\Ar{1-2}{2-2}{"{\rm Fp}(P\down|)", ""'}
\Ar{1-1}{2-1}{"{{\rm fp}(P\down|)}"'}
\Ar{2-1}{2-2}{"",""'}
\end{tikzcd},
\end{equation}
where the horizontal functors are embeddings.
\end{enumerate}
\end{theorem}

\subsection{Relations with other works}
\label{ss:relations-other-works}

It is possible
to recover some results from \cite{HM} using our general setting. For example, Proposition \ref{Prop 2.6} corresponds to \cite[Proposition 2.6]{HM}.
In addition, there have been recent papers that discuss the construction of G-precoverings by taking quotients. In \cite{HZhao}, the authors investigate G-coverings between relative stable categories in Theorem 5.5. Then, in \cite[Corollary 5.6]{HZhao}, they recover \cite[Proposition 2.6]{HM}. This result is a direct application of Theorem \ref{thm:general-case} using Example \ref{ex-G-stable} and the stable 2-functor $U_p$ defined in Example \ref{exm:stable-2-fun}. Another paper, \cite{HZou}, introduces the concept of G-liftable ideals and proceeds to construct Galois G-coverings of quotient categories in \cite[Theorem 3.9]{HZou} associated with these G-liftable ideals. Again, their results can be obtained by applying Theorem \ref{thm:general-case} using the 2-functor $U_i$ defined in Example \ref{exm:factor-2-fun}.

A functor $\Phi:\CF(\mmod \CC)\rt \CF(\mmod (\CC/G))$ is defined in  \cite[Section 5]{P} and proved in \cite[Theorem 5.5]{P} to be  $G$-precovering and  satisfies  properties similar to the pushdown $P\down$.  This functor plays an essential role in \cite{P} in  proving  the stability of the Krull-Gabriel dimension under Galois coverings. Using \cite[Proposition 5.3]{P}, we can deduce the existence of a natural isomorphism $\CF( P \down)\cong \Phi$. Therefore, using our morphism method, we can recover the functor $\Phi$ and its fundamental properties established in \cite[Theorem 5.5]{P}. Moreover, it will be proven that the functor $\Phi$ satisfies  similar conditions as those given in Definition \ref{dfn:precov-adj}, see Theorem \ref{thm:G-Fp-precovering} for more details.

We should remark that Z. Leszczy\'nski and A. Skowro{\'n}ski in \cite[Theorem~5]{LSk} showed when the (upper) triangular matrix algebra $T_2(\La)$ of a finite-dimensional algebra $\La$ is of polynomial growth, then  there is a  Galois functor $F^{(2)}_{\la}:\mmod T_2(\tilde{\La})\rt \mmod T_2(\La)$ which is induced by a Galois covering $F:\tilde{\La}\rt \La$
with $\tilde{\La}$ a simply connected locally bounded $k$-category. 
Next, they applied the  Galois functor  $F^{(2)}_{\la}$ to classify algebras $\La$ for which the triangular matrix algebra $T_2(\La)$ is of tame representation type.
Note also that \cite[Theorem 5]{LSk} together with \cite[Theorem 2.5]{Le94} is used to show that
a commutative 2-dimensional grid%
\footnote{This is a bound quiver defined as the product quiver of the equi-oriented Dynkin quiver of types $A_m$ and $A_n$ for $m,n \ge 1$
with full commutativity relations.} $G_{m,n}$ with $m \ge n \ge 2$
is wild 
iff $mn \ge 12$, and is representation-finite iff
$n = 2$ and $2 \le m \le 4$, which is used in TDA (topological data analysis).
What we have done in our work is to reformulate a covering theory for triangular matrix algebras in terms of functor categories and morphism categories
using the fact that
the module category $\mmod T_2(\La)$ is equivalent to
the morphism category $\rmH(\mmod \La)$ of $\mmod \La$.

\section{Preliminary}\label{Perliminaries}

In this section, for the convenience of the reader, we present some definitions and results that will be used throughout the paper.

Throughout this paper $G$ denotes a group, $k$ is a field,
and the category of $k$-vector spaces are denoted by $\Mod k$.
A category $\CC$ is called a {\em $k $-linear} category,
or a {\em $k$-category} for short,
if the morphism sets are $k$-vector spaces and the compositions of morphisms are $k$-bilinear.
A functor $F \colon \CC \to \CC'$ between $k$-categories are said to be {\em $k$-linear}
if the induced map $F_{(x,y)} \colon \CC(x,y) \to \CC'(Fx,Fy)$ is $k$-linear
for all objects $x, y$ of $\CC$.
We always assume that functors between $k$-categories are $k$-linear.

We fix a (Grothendieck) universe $\frkU$ containing the set $\bbN$ of natural numbers once for all.
A set is called a {\em small} set (resp.\ a {\em $1$-class}) if it is an element (resp.\ a subset) of $\frkU$.
A {\em small} (resp.\ {\em light}) category is a category $\calC$ whose objects form a small set (resp.\ a $1$-class), and
$\calC(x,y)$ is a small set for all objects $x, y$ of $\calC$.
We denote by $\kCAT$ the 2-category whose objects are the light $k$-categories, 
whose 1-morphisms are the $k$-functors between objects,
and whose 2-morphisms are the natural transformations between 1-morphisms,
also by $\kCat$ the full 2-subcategory of $\kCAT$ whose objects are the small $k$-categories.
Note that if $\calC$ is an object of $\kCat$, then $\Mod \calC$ (see the subsection \ref{ssec:fun-cat} for the definition) is an object of $\kCAT$
(see, for instance, \cite[Chapter 4, Appendix A]{Asa-book} for details on $2$-categories and a set theoretical foundation, respectively).

A \emph{skeleton} of a $k$-category $\calC$ is a full subcategory of $\calC$ whose objects form a complete set of representatives of the isomorphism classes of objects in $\calC$, which is clearly equivalent to $\calC$. Therefore a skeleton of $\calC$ is uniquely determined up to equivalences, and is denoted by $\sk(\calC)$.  A category is said to be \emph{skeletally small}, if its skeleton is a small category.
Now let $\calC$ be a skeletally small category.
Then since $\sk(\calC)$ is equivalent to $\calC$, 
their module categories are equivalent: $\Mod\sk(\calC) \simeq \Mod\calC$, the former is again skeletally small and a light category.
Hence $\Mod\calC$ is also skeletally small, and is equivalent to
a light category, although $\Mod\calC$ itself is not a light category
(it is 2-moderate in a terminology defined in \cite[Appendix A]{Asa-book}).
Therefore, we may consider $\Mod(\Mod\calC)$ without any set theoretical problems, and is again skeletally small and equivalent to a light category $\Mod(\sk(\Mod\calC))$ ($\Mod(\Mod\calC)$ itself is 3-moderate).
Hence also there are no problems to consider subcategories of $\Mod(\Mod\calC)$ such as the subcategory $\CF(\Mod\calC)$ consisting of finitely presented objects in $\Mod(\Mod\calC)$,
which we are interested in.
By this reason, we usually assume that a $k$-category $\calC$ is
skeletally small when we consider its module category $\Mod\calC$.

\subsection{Orbit categories}

\begin{definition}[$G$-Categories]
A {\em category with a $G$-action} (or a {\em $G$-category} for short) is a pair
$(\CC, A)$ of a $k$-category $\CC$ and a group homomorphism $A\colon G \to \Aut(\CC)$,
where $\Aut(\CC)$ is the group of automorphisms of $\CC$.
We set $ A_a:=A(a) $  for all $ a \in G $.

If there seems to be no confusion,
we simply write $\CC = (\CC, A)$ by omitting $A$,
and we denote $G$-actions by the same letter $A$, also
we usually write $ax:= A_{a}(x)$ and $af:= A_a(f)$
for all $a \in G$, $x \in \CC$, and all morphisms $f$ in $\CC$.

Recall that the action of $G$ is said to be {\em free}
if $ax \neq x$ for all $x \in \CC$ and $a \in G \setminus \{1_G\}$,
or equivalently if the map  $\rho_x:G\rt Gx\, (:= \{ax \mid a \in G\})$, $a\mapsto ax$ is a bijection
for all $x \in \CC$. 
\end{definition}

\begin{definition}[$G$-Invariant Functors]
\label{Def-G-invariant}
Let $\calC$ be a $G$-category, and $\calC'$ a $k$-category.
Then a {\em $G$-invariant functor} from $\calC$ to $\calC'$ is a pair $(F, \ph)$
of a functor $F: \CC \to \CC'$ and
a family $\ph:=(\ph_{a })_{a \in G}$
of natural isomorphisms $\ph_{a}: F \To FA_{a}$ such that
the following diagram of natural isomorphisms commutes for all $a, b \in G$:
\[
\begin{tikzcd}
F & FA_a\\
FA_{ba} & FA_bA_a
\Ar{1-1}{1-2}{"\ph_a", Rightarrow}
\Ar{2-1}{2-2}{"", equal}
\Ar{1-1}{2-1}{"\ph_{ba}"', Rightarrow}
\Ar{1-2}{2-2}{"\ph_b A_a", Rightarrow}.
\end{tikzcd}
\]
This family $\ph$ was called an {\em invariant adjuster} of $F$ in \cite[Definition 1.1]{A},
but here we call it a {\em structure} of the $G$-invariant functor.
\end{definition}
In the above, we denote the value of $\ph_a$ at an object $x \in \CC$ by
$\ph_{a,x} \colon Fx \to FA_a x$.
Then the diagram above is commutative if and only if
the following diagram of morphisms in $\CC'$ is commutative for all $x \in \CC$:
\[
\begin{tikzcd}
Fx & FA_a x\\
FA_{ba} x & FA_bA_a x
\Ar{1-1}{1-2}{"\ph_{a,x}"}
\Ar{2-1}{2-2}{"", equal}
\Ar{1-1}{2-1}{"\ph_{ba, x}"'}
\Ar{1-2}{2-2}{"\ph_{b, ax}"}.
\end{tikzcd}
\]
Note that for any $a \in G$ and $x \in \CC$,
the commutativity of the diagram shows that
$\ph_{1,x} = 1_{Fx}$ and $\ph\inv_{a , x}= \ph_{a\inv, ax}$ (\cite[Remark 1.2]{A}).
The definition above is also given in \cite[Defintion 2.3]{BL}
under the name of a {\em $G$-stable functor}.
 
\begin{definition}
[$G$-equivariant functors]

Let $\CC =(\calC, A)$ and $\CC'= (\calC', A')$ be $G$-categories.
Then a $G$-{\em equivariant} functor from $\calC$ to $\calC'$ is a pair $(F,\ph)$
of a functor $F\colon \calC \to \calC'$ and a family 
$\ph = (\ph_{a})_{a\in G}$ of natural isomorphisms
$\ph_{a}\colon A'_{a}F \To FA_{a}$ ($a \in G$)
such that the following diagram of natural isomorphisms commutes for all  $a,b \in G$:
\[
\begin{tikzcd}
A'_b A'_a F & A'_bFA_a & FA_b A_a\\
A'_{ba}F && FA_{ba}
\Ar{1-1}{1-2}{"A'_b\ph_a", Rightarrow}
\Ar{1-2}{1-3}{"\ph_b A_a", Rightarrow}
\Ar{2-1}{2-3}{"\ph_{ba}"', Rightarrow}
\Ar{1-1}{2-1}{"", equal}
\Ar{1-3}{2-3}{"", equal}.
\end{tikzcd}
\]
This $\ph$ was called an {\em equivariance adjuster} of $F$ in \cite[Definition 4.8]{A},
but here we call it a {\em structure} of the $G$-equivariant functor.
\end{definition}

\begin{definition}[Morphisms between $G$-equivariant functors]
Let $\CC =(\calC, A)$ and $\CC'= (\calC', A')$ be $G$-categories, and
$(F,\ph),(F',\ph')  \colon \calC \to \calC'$ $G$-equivariant functors.
Then a morphism from $(F,\ph)$ to $(F',\ph')$ is a natural transformation
$\al \colon F \To F'$ such that the following diagram of natural transformations commutes:
\[
\begin{tikzcd}
A'_aF & FA_a\\
A'_aF' & F'A_a
\Ar{1-1}{1-2}{"\ph_a", Rightarrow}
\Ar{2-1}{2-2}{"\ph'_a", Rightarrow}
\Ar{1-1}{2-1}{"A'_a \al", Rightarrow}
\Ar{1-2}{2-2}{"\al A_a", Rightarrow}
\end{tikzcd}.
\]
\end{definition}

\begin{definition}
\label{dfn:GCat}
By $\GCat$ we denote the $2$-category whose objects are the small $G$-categories,
whose $1$-morphisms are the $G$-equivariant functors between objects,
whose $2$-morphisms are the morphisms between $1$-morphisms.
$\GCAT$ is similarly defined, where the objects are the light $G$-categories 
(see Definition \ref{dfn:2-category} for details).
\end{definition}

\begin{remark}
We regard every $k$-category $\calC$ as a $G$-category $(\calC, 1)$, where
$1$ is the trivial action $1 \colon G \to \Aut(\calC)$.
This defines an embedding $\kCAT \to \GCAT$.
Then a $G$-invariant functor in Definition \ref{Def-G-invariant} is nothing
but a $G$-equivarinat functor.
\end{remark}

We freely apply the remark above throughout the paper.
 
\begin{definition}[$G$-coverings]\label{Def.G-precovering}
Let $(F,\ph): \CC \to \CC'$ be a $G$-invariant functor from a $G$-category
$\CC = (\CC, A)$ to a $k$-category $\CC'$.
Then $(F, \ph)$ is called a {\em $G$-precovering} if for any $x, y \in\CC$
the following two $k$-homomorphisms are isomorphisms.
\[
\begin{aligned}
(F, \ph)^{(1)}_{x,y}&: \bigoplus_{a \in G} \CC(ax, y) \to \CC'(Fx, Fy), \quad
    (f_{a })_{a \in G} \mapsto \sum_{a  \in G} F(f_{a }) \ph_{a , x};\\
(F, \ph)^{(2)}_{x, y} &: \bigoplus_{b \in G} \CC(x, b y)\to \CC'(Fx, Fy), \quad
    (f_{b})_{b \in G} \mapsto \sum_{b \in G}\ph_{b^{-1}, b y} F(f_{b}).
\end{aligned}
\]
In addition, if $F$ is dense, then $(F,\ph)$ is called a {\em $G$-covering} {\cite[Definition 1.7]{A}}.
\end{definition}

As proved in \cite[Proposition 1.6]{A}, in the above definition it is enough to check that all $F^{(1)}_{x,y}$ are isomorphisms, or all $F^{(2)}_{x,y}$ are isomorphisms.

\begin{definition}[Orbit Categories]
Let $\CC$ be a $G$-category.
Then the {\em orbit category} $\CC/G$ of $\CC$ by $G$ is defined as follows:
\begin{enumerate}
\item
The class of objects is the same as that of $\CC$.
\item
For any $x,y \in \CC/G$, the morphism set $\CC/G(x,y)$ is given by
\[
\left\{f=(f_{b,a })_{(a , b)} \in \prod_{(a , b)\in G\times G} \CC(a  x, b y) \right.\left.\ \rule[-3ex]{0.5pt}{7ex}\ \begin{array}{l} f \text{ is row-finite, column-finite}\\
\text{and } f_{c b, c a }= c (f_{b, a }), \forall c \in G\end{array}\right\},
\]
%
where $f$ is said to be {\em row-finite} (resp. {\rm column-finite}) if for any $a \in G$ the set $ \{ b \in G \ \vert \ f_{a,b} \neq 0 \}$ (resp. $ \{ b \in G \ \vert \ f_{b,a} \neq 0 \} $) is finite.
\item
For two composable  morphisms $x \st{f}\to y \st{g}\to z$ in $\CC/G$, we set
$$gf:= (\sum_{c \in G}g_{b, c} f_{c, a })_{(a , b) \in G\times G}.$$
\end{enumerate}

There is a canonical functor $P : \CC \to \CC/G$ given by $P(x)=x$ and $P(f)= (\delta_{a , b} a  f)_{(a , b)}$ for all $x, y \in \CC$ and $f \in \CC(x, y)$, which is extended to a $G$-invariant functor $P=(P, \ph): \CC \to \CC/G$
by defining $\ph$ 
as follows:
Set 
$\ph_{c, x}:=(\delta_{a  , b c}1_{a  x})_{(a , b)}
\in \CC/G(Px, P c x)$
for all $ c \in G$ and $x \in \CC$, $\ph_{c}:=(\ph_{c,x})_{x \in \CC}: P \to PA_{c}$ and $\ph:=(\ph_{c})_{c \in G}$ {\cite[Definitions 2.4 and 2.5]{A}}.
Then $(P, \ph)$ turns out to be a $G$-covering,
and is called the {\em canonical} $G$-covering associated to the orbit category $\calC/G$ \cite[Proposition 2.6]{A}.
\end{definition}

\subsection{Functor categories}
\label{ssec:fun-cat}

Let $\CC$ be a skeletally small category.
A contravariant functor from $\CC$ to the category $\Mod k$ is
called a (right) {\em $\CC$-module}.
Any object $x$ of $\CC$ serves us a contravariant functor
$\calC(\blank, x) \colon \CC \to \Mod k$
{\em represented by} $x$ called a {\em representable functor},
which is a typical example of a right $\CC$-module.
All the $\CC$-modules and the natural transformations between them form an abelian category
as the functor category from $\CC^\op$ to $\Mod k$, which is
denoted by $\Mod\CC$. The hom sets in $\Mod \CC$ is sometimes denoted by $\Hom_{\CC}(\blank, \blank)$.
It is well known that 
$\Mod\CC$ {\em has small coproducts}, i.e., has arbitrary direct sums with small index sets.
%

An object $F$ of $\Mod\CC$ is said to be {\em finitely generated}
if there exists an epimorphism from a finite direct sum of representable functors to $F$, that is,
if we have an epimorphism
\begin{equation}\label{eq:fg}
\bigoplus_{i=1}^n{\CC}(\blank, x_i) \to F
\end{equation}
for some finitely many objects $x_1,\dots, x_n$ of $\CC$.
The full subcategory of $\Mod\calC$ consisting of all finitely generated $\CC$-modules is denoted by $\mmod\CC$.

A finitely generated $\CC$-module $F$ is projective if and only if it is isomorphic to a direct summand of a finite direct sum of representable functors.
We denote by $\prj \CC$ the full subcategory of $\mmod \CC$ consisting of all finitely generated projective $\CC$-modules.

Assume further that $\CC$ is additive.
Then note that a $\CC$-module $F$ is finitely generated if and only if
there exists an exact sequence
\[
\CC(\blank, x) \to F \to 0
\]
for some object $x$ of $\CC$ because
$\Ds_{i=1}^n \CC(\blank, x_i) \cong \CC(\blank, \Ds_{i=1}^n x_i)$
in \eqref{eq:fg}.
A $\CC$-module $F$ is called {\em finitely presented} if there exists an exact sequence
\[
\CC(\blank, x') \to \CC(\blank, x) \to F \to 0
\]
for some $x, x' \in \CC$.
In this paper, the full subcategory of $\Mod \CC$ consisting of all
finitely presented $\CC$-modules is denoted by $\CF(\CC)$
\footnote{This was denoted by $\calF(\CC)$ in \cite{P}, but
we did not use this notation in this paper
because $\calF$ is sometimes used to denote the subcategory of modules having a filtration
with factors in a set of prescribed modules.
In the literature, the notation $\mmod \CC$ is usually used,
but we have already adapted it for the subcategory of finitely generated modules.
The notation $\CF$ is also used in papers (Henning's paper and Bill's paper).}.

It is proved by Auslander that $\CF(\CC)$ is an abelian category if and only if $\CC$ admits weak kernels \cite[Proposition 2.1]{Au1}.
This is the case, for instance, when $\CC$ is a contravariantly finite subcategory of an abelian category. Let $\CA$ be an abelian category.
For $\CC \subseteq \CA$, we denote by $\underline{\CC}$ the stable category of $\CC.$ The canonical functor $\pi:\CC\rt \underline{\CC}$ induces a fully faithful functor from $\CF(\underline{\CC})$ to $\CF(\CC),$ where the image is the subcategory consisting of those functors that vanish on projective objects.
Therefore, $\CF(\udl{\CC})$ can be considered as a subcategory of $\CF(\CC)$.

For a full subcategory $\CK$ of an additive category $\CC$,
we denote by $\add\CK$ (resp.\ $\Add\CK$),  the full subcategory of $\CC$ consisting of the direct summands of finite (resp.\ small) coproducts of objects in $\CK$.

Let $\DD$ be a class of objects of a category $\CX$.
Then the ideal of $\CX$ formed by all morphisms factoring through
a finite direct sum of objects in $\DD$ is denoted by $\ang{\DD}$.
%
The factor category $\frac{\CX}{\ang{\DD}}$ has the same objects as  $ \CX$ and for any  $X, Y \in \CX$
\[
\frac{\CX}{\ang{\DD}}(X, Y ) := \frac{\mathcal{X}(X, Y )} {\ang{\DD}(X, Y )}.
\]

\subsection{Pushdown functor}
\label{Pushdownfunctor}
Let $\CC = (\calC, A)$ be a $G$-category.
Then the $k$-category $\Mod\calC$ turns out to be a $G$-category $\Mod\calC = (\Mod\calC, \Mod A)$ by defining a $G$-action $\Mod A$ on $\Mod \CC$ as follows: For each $a  \in G$, $(\Mod A)_{a}:= \ovl{A}_a \colon \Mod \CC \to \Mod \CC$ is an isomorphism sending each morphisms $f \colon X \to Y$ in $\Mod \CC$ to ${}^a f \colon {}^a X \to {}^a Y$, where
\[
{}^{a}u := u \circ A_{a }^{-1}
\]
for each $u = X, f, Y$.
Note that ${}^{a }\CC(\blank,x) = \CC(a ^{-1}(\blank),x) \cong \CC(\blank, a  x)$ for all $x \in \CC$.
Then we see that the $G$-action $\Mod A$ on $\Mod\calC$ restricts to the $G$-action (denoted by $\prj A$) on the category $\prj \calC$, and hence to the $G$-action (denoted by $\mmod A$) on the category $\mmod\calC$ because the epimorphism in \eqref{eq:fg} induces an epimorphism
$\bigoplus_{i=1}^n{\CC}(\blank, a x_i) \to {}^a F$.
In this way we regard $\mmod \calC = (\mmod \calC, \mmod A)$ and $\prj \calC = (\prj \calC, \prj A)$ as $G$-categories.

\begin{definition}
\label{dfn:ind-C}
In the following, we denote by $\ind\CC$ a full subcategory of $\mmod \CC$ whose objects form a complete set of representatives of isoclasses of indecomposable modules in $\mmod \CC$ {\em that is closed under the $G$-action} $\Mod A$.
This makes it possible to restrict the $G$-action $\Mod A$ to the $G$-action (denoted by $\ind A$) on $\ind \calC$, and we regard $\ind \calC =(\ind \calC, \ind A)$ as a $G$-category.
\end{definition}

If there seems to be no confusion, we write the $G$-action on $\Mod A, \mmod A, \prj A$, and $ \ind A$ simply by $\ovl{A}$.

Let $\CC$ be a skeletally small $G$-category.
The canonical functor $P: \CC \to \CC/G$ induces a functor
$P\up: \Mod (\CC/G) \to \Mod \CC$, given by $P\up X =X \circ P^\op$
for all $X \in \Mod (\CC/G)$, where $P^\op \colon \calC^\op \to (\calC/G)^\op$ is a functor defined by $P$ in an obvious way.
This functor is called the {\it pullup} of $P$.
It is well known that $P\up$ has a left adjoint $P\down$
called the {\em pushdown} of $P$,
the precise form of which was given in the proof of \cite[Theorem 4.3]{A}.
In fact, $P\down: \Mod \CC \to \Mod (\CC/G)$ is defined as follows. Let $X \in \Mod \CC$.
Then for an object $x \in \obj(\CC/G)= \obj(\CC)$,
$(P_{\centerdot}X)(x):= \bigoplus_{a  \in G}X(a  x)$.
For a morphism $f:x \rt y$ in $\CC/G$, where $f=(f_{b, a })_{(a , b) \in G \times G} \in \CC/G(x,y)$, $(P_{\centerdot}X)(f)$ is given by the following commutative diagram
\[
\begin{tikzcd}[column sep=80pt]
(P\down X)(y) & (P\down X)(x)\\
\Ds_{b\in G} X(by) & \Ds_{a\in G} X(ax)
\Ar{1-1}{1-2}{"(P\down X)(f)"}
\Ar{2-1}{2-2}{"{(X(f_{b, a}))_{(a, b) \in G\times G}}"'}
\Ar{1-1}{2-1}{equal}
\Ar{1-2}{2-2}{equal}
\end{tikzcd}.
\]
Also, if $u: X \to X'$ is a morphism in $\Mod \CC$, then for every $x \in {\rm Obj}(\CC/G)$, $(P_{\centerdot}u)_{x}$ is defined by the following commutative diagram
\[
\begin{tikzcd}[column sep=60pt]
(P\down X)(x) & (P\down X')(x)\\
\Ds_{a\in G} X(ax) & \Ds_{a\in G} X'(ax)
\Ar{1-1}{1-2}{"(P\down u)_x"}
\Ar{2-1}{2-2}{"\Ds_{a\in G}u_{ax}"'}
\Ar{1-1}{2-1}{equal}
\Ar{1-2}{2-2}{equal}
\end{tikzcd}.
\]
For each $X \in \Mod \CC$ and $Y \in \Mod \CC/G$ the adjunction
\[
\th_{X, Y} \colon
\Hom_{\CC/G}(P\down X, Y) \to \Hom_{\CC}(X, P\up Y)
\]
is given by
$(\th_{X, Y}\, \al)_x:=\al_{x,1}:X(x) \to Y(x)=Y(Px)=(P\up Y)(x)$ for all $x \in \CC$ and $\al = (\al_x)_{x \in \CC/G} \in \Hom_{\CC/G}(P\down X, Y)$, where
\[
\al_x = (\al_{x,a})_{a\in G} \colon (P\down X)(x) = \Ds_{a \in G}X(ax) \to Y(x)
\]
for all $x \in \CC/G$;
and its inverse
\[
\th\inv_{X, Y} \colon
\Hom_{\CC}(X, P\up Y) \to \Hom_{\CC/G}(P\down X, Y)
\]
is given by
$(\th\inv_{X,Y}\, f)_x:=(Y(\phi_{a, x})f_{ax})_{a \in G}
\colon (P\down X)(x) = \Ds_{a \in G}X(ax) \to Y(x)$ 
for all $f \in \Hom_{\CC}(X, P\up Y)$ and $x\in \CC/G$,
where $Y(\phi_{a, x})f_{ax} \colon X(ax) \to Y(x)$ is the composite
\[
X(ax) \ya{f_{ax}} Y(ax) \ya{Y(\phi_{a, x})} Y(x),
\]
see \cite[proof of Theorem 4.3]{A}.
Moreover, 
$P\down$ is extended to a $G$-invariant functor as follows.
For each $c \in G$, define a morphism
${\ph\down}_c \colon P\down \to P\down \circ \ovl{A}_{c}$ by ${\ph\down}_{c}:=({\ph\down}_{c, X})_{X \in \Mod \CC}$,
where ${\ph\down}_{c, X}$ is defined by the commutative diagram
\[
\begin{tikzcd}[column sep=110pt]
(P\down X)(x) & (P\down\,{}^{c}\!X)(x)\\
\Ds_{a \in G} X(ax) & \Ds_{b \in G}X(c\inv bx)
\Ar{1-1}{1-2}{"{{\ph\down}_{c,X,x}}"}
\Ar{2-1}{2-2}{"{(\delta_{a, c^{-1}b}\idty_{X(ax)})_{(a, b ) \in G \times G}}"'}
\Ar{1-1}{2-1}{equal}
\Ar{1-2}{2-2}{equal}
\end{tikzcd}
\]
for all $x \in \CC$.
Then ${\ph\down}_{c}$ turns out to be a natural isomorphism for all $c \in G$, and then by setting
$\ph\down:=({\ph\down}_{c})_{c \in G}$,
we obtain a $G$-invariant functor $(P\down, \ph\down)$.

Finally, as mentioned in the proof, we have the equalities
\[
(P^{\centerdot}P_{\centerdot}X)(x)=\Ds_{a \in G}X(ax)=(\Ds_{a \in G}{}^{a^{-1}}X)(x)
\]
for all $X \in \Mod \CC$
and $x \in \CC$.
Let $\si_b = (\de_{a,b}\idty_{\,{}^b\!X})_{a\in G} \colon {}^bX \to \Ds_{a\in G} {}^aX$ and
$\si'_b = (\de_{a,b^{-1}}\idty_{\,{}^b\!X})_{a\in G} \colon {}^bX \to \Ds_{a\in G} {}^{a\inv}X$
be the canonical injections for all $b \in G$.
Then the universality of the coproduct gives us an isomorphism
$t_X \colon \Ds_{a\in G}{}^aX \to \Ds_{a\in G}{}^{a\inv}X$ such that the diagram
\[
\begin{tikzcd}
{}^bX & \Ds_{a\in G}{}^{a\inv}X\\
\Ds_{a\in G}{}^aX
\Ar{1-1}{1-2}{"\si'_{b}"}
\Ar{1-1}{2-1}{"\si_{b}"'}
\Ar{2-1}{1-2}{"t_X"'}
\end{tikzcd}
\]
commutes. 
Thus the family $t = (t_X)_{X \in \Mod \calC}$
turns out to be a natural isomorphism
\[
t \colon \Ds_{a \in G} {}^a(\blank) \To P\up P\down,
\]
where
$\Ds_{a \in G} {}^a(\blank) \colon \Mod \CC \to \Mod \CC$ is a functor defined by
sending $f \colon X \to Y$ to
$\Ds_{a\in G}{}^a f \colon \Ds_{a\in G}{}^a X \to \Ds_{a\in G}{}^a Y$.
%
%
By using $(P\down, \phi\down)$
and the natural isomorphism
$t$ 
above,
we have the following vital commutative diagram.

\begin{proposition}[{\cite[Proof of Theorem 4.3]{A}}]
\label{Vitalcommutativediagram}
For any $X, Y \in \mmod \CC$, 
we have a commutative diagram
\[
\begin{tikzcd}
\Ds_{a \in G} \mmod\CC(X, {}^{a} Y) & \Mod\CC(X, \Ds_{a \in G} {}^{a} Y)\\
\mmod \CC/G(P\down X, P\down Y) & \Mod\CC(X, P\up P\down Y)
\Ar{1-1}{1-2}{"\nu", "\sim"'}
\Ar{1-2}{2-2}{"{\Mod\CC(X, t_Y)}", "\rotatebox{90}{$\sim$}"'}
\Ar{1-1}{2-1}{"{(P\down, \ph\down)^{(2)}_{X,Y}}"'}
\Ar{2-1}{2-2}{"\sim","\th_{X, P\down Y}"'}
\end{tikzcd},
\]
where $\nu$ is the canonical isomorphism, and $\th$ is an adjunction
of the adjoint $P\down \adj P\up$.
\end{proposition}

\subsection{Locally bounded categories}

A locally bounded category is a $k$-category $\CC$ satisfying the following conditions:
\begin{itemize}
\item The endomorphism $\CC(x,x)$ is local for every $x \in \CC$.
\item $\CC$ is basic, i.e., if $x\neq y$, then $x\ncong y$.
\item For each $x \in \CC$, $\sum_{y \in \CC} [\CC(x,y): k] <\infty$ and $\sum_{y \in \CC} [\CC(y,x):k] <\infty$.
\end{itemize}
Observe that if $\CC$ is a locally bounded $k$-category with a free $G$-action, then the orbit category $\CC/G$ is equivalent to the classical orbit category in Gabriel's sense.

\begin{remark}
\label{rmk:fg-lb}
Let $\CC$ be a locally bounded $k$-category.
Then a $\CC$-module $F$ is finitely generated and projective
if and only if
$F$ is isomorphic to the direct sum of a finite number of modules of the form $\CC(-,x)$ for some $x \in \CC$,
which satisfies $\sum_{y \in \CC} [\CC(y,x):k] <\infty$ by definition.
Therefore, in particular,  $F$ is finitely generated if and only if it is finite-dimensional, i.e., $\sum_{x \in \CC} [F(x): k]< \infty$, see \cite[2.2]{BG}.
\end{remark}

We need the following result throughout the paper.

\begin{lemma} \cite[Lemma 3.5]{G}\label{Inde-Pushdown}
Let $\CC$ be a locally bounded $G$-category and $G$ acts  freely on $\ind\CC$. Then the pushdown functor $P_{\centerdot}:\mmod \CC \to \mmod (\CC/G)$ maps any indecomposable $\CC$-module to an indecomposable $\CC/G$-module.
\end{lemma}

Let $\CC$ be a $G$-category and $M$ a $\CC$-module. 
We denote by $\supp M$
the {\em support} of $M$, i.e., 
the full subcategory of $\CC$ consisting of 
all objects $x$ of $\CC$ such that $M(x)\neq 0$.
Let $x$ be an object of $\CC$.
We denote by $\CC_x$ the full subcategory of $\CC$ formed by
the vertices of all $\supp M$, where $M \in \ind \CC$ and $M(x)\neq 0$.
A locally bounded $k$-category $\CC$ is called {\em locally support-finite}
if for every $x \in \CC$, $\CC_x$ is finite.

\newcommand\DLSthm{\cite[Theorem]{DLS}}
\begin{proposition}[\DLSthm]
\label{loc-supp-fin}
Let $\CC$ be a locally support-finite $G$-category and $G$ acts freely on $\ind \CC$. Then the  pushdown functor induces a $G$-covering
$$P_{\centerdot}: \mmod \CC \to \mmod (\CC/G).$$
\end{proposition}

\section{Extension of a 2-endofunctor of $\kCAT$ to $\GCAT$}

To unify the theorems on the existence of a $G$-precovering,
we use $2$-categorical approach.
In this section, we prepare necessary facts for the later use.
Let us start with the following definitions.

\begin{definition}[{\cite[Definition 4.1.1]{Asa-book}}]
\label{dfn:2-category}
A $2$-category is a sequence of the following data that satisfies the axioms below.

Data:
\begin{itemize}
    \item A non-empty set $\sfC_0$,
    \item A family of categories $(\sfC(x,y))_{x,y \in \sfC_0}$,
    \item A family of functors $\circ:= (\circ_{x,y,z}:\sfC(y,z) \times \sfC(x,y) \longrightarrow \sfC(x,z))_{x,y,z \in \sfC_0}$,
    \item A family of functors $(u_x: 1 \longrightarrow \sfC(x,y))_{x \in \sfC_0}$, where $1$ denotes the category having only a single object $*$ and only a single morphism $\id_{*}$.
\end{itemize}
Axioms:
\begin{itemize}
    \item  (associativity) The following diagram is commutative for all $x,y,z,w \in \sfC_0$;
    \begin{equation}
\begin{tikzcd}
	{\sfC(z,w)\times \sfC(y,z) \times \sfC(x,y)} &&& {\sfC(y,w) \times \sfC(x,y)} \\
	\\
	{\sfC(z,w) \times \sfC(x,z)} &&& {\sfC(x,w)}
	\arrow["{1\times \circ}"', from=1-1, to=3-1]
	\arrow["{\circ \times 1}", from=1-1, to=1-4]
	\arrow["\circ", from=1-4, to=3-4]
	\arrow["\circ"', from=3-1, to=3-4]
\end{tikzcd}
\end{equation}
    \item (unitality) The following diagram is commutative for all $x,y \in \sfC_0$:
    \begin{equation}
\begin{tikzcd}
	{1 \times \sfC(x,y)} && {\sfC(x,y) \times 1} \\
	& {\sfC(x,y)} \\
	{\sfC(y,y) \times  \sfC(x,y)} && {\sfC(x,y) \times \sfC(x,x)}
	\arrow["{u_y \times \sfC(x,y)}"', from=1-1, to=3-1]
	\arrow["{\sfC(x,y) \times u_x}", from=1-3, to=3-3]
	\arrow["{\pr_2}", from=1-1, to=2-2]
	\arrow["{\pr_1}"', from=1-3, to=2-2]
	\arrow["\circ"', from=3-1, to=2-2]
	\arrow["\circ", from=3-3, to=2-2]
\end{tikzcd},
\end{equation}
\end{itemize}
where $\pr_i$ denotes the $i$-th projections for all $i = 1,2$.
Elements of $\sfC_0$ are called objects of $\sfC$, objects (resp. morphisms, compositions) of $\sfC(x,y)$ (with $x,y \in \sfC_0$) are called $1$-morphisms ($2$-morphisms, vertical compositions) of $\sfC$. Furthermore, $\circ_{x,y,z}$ (with $x,y,z \in \sfC$) is called horizantal compositions of $\sfC$,
and we set $\idty_{x}:=u_x(*)$ and $\idty_{\idty_x}:= u_x(\idty_*)$ for all $x \in \sfC_0$.
We denote the vertical compositions (i.e., compositions of the categories $\sfC(x,y)$)  by a symbol $\bullet$ to distinguish it from the horizontal compositions $\circ$. Sometimes we omit the symbol of the horizontal composition to write $gf$ for $g\circ f$.
In a $2$-category $\sfC$ we display objects $x,y \in \sfC_0$, $1$-morphisms $f,g \in \sfC(x,y)$, and  $2$-morphisms $\alpha \in \sfC(x,y)(f,g)$ as follows:

\[
\begin{tikzcd}[row sep=50pt, column sep=60pt]
x & y
\Ar{1-1}{1-2}{"f", ""'{name=a}, bend left}
\Ar{1-1}{1-2}{"g"', ""{name=b}, bend right}
\Ar{a}{b}{"\al", Rightarrow}
\end{tikzcd}.
\]
\end{definition}

\begin{definition}[{\cite[Definition 4.1.11]{Asa-book}}]\label{2-functor}
Let $\sfC$ and  $\sfD$ be $2$-categories. A pair of the following data that satisfies the axioms below is called a $2$-functor from $\sfC$ to $\sfD$, and is denoted by $\sfC \rightarrow \sfD$;

Data:
\begin{itemize}
    \item A map $X: \sfC_0 \rightarrow \sfD_0$,
    \item A family of functors $(X_{(x,y)}: \sfC(x,y) \rightarrow \sfD(X(x),X(y)))_{(x,y)\in \sfC_0 \times \sfC_0}$, where we shortly write $X(f)$ for $X_{(x,y)}(f)$ for all $f \in \sfC(x,y)_0 \cup \sfC(x,y)_1$.
\end{itemize}

Axioms:
\begin{itemize}
    \item For each $x,y,z \in \sfC_0$ the following is commutative:
    \[\begin{tikzcd}
	{\sfC(y,z) \times \sfC(x,y) } && {\sfC(x,z)} \\
	\\
	{\sfD(X(y),X(z)) \times \sfD(X(x),X(y))} && {\sfD(X(x),X(z))}
	\arrow["\circ", from=1-1, to=1-3]
	\arrow["\circ"', from=3-1, to=3-3]
	\arrow["{X_{(y,z)} \times X_{(x,y)}}"', from=1-1, to=3-1]
	\arrow["{X_{(x,z)}}", from=1-3, to=3-3]
\end{tikzcd}\]
    \item For each $x\in \sfC_0$ the following is commutative:
    \[\begin{tikzcd}
	{\mathbf{1}} && {\sfC(x,x)} \\
	\\
	&& {\sfD(X(x),X(x))}
	\arrow["{u_x}", from=1-1, to=1-3]
	\arrow["{X_{(x,x)}}", from=1-3, to=3-3]
	\arrow["{u_{X(x)}}"', from=1-1, to=3-3]
\end{tikzcd}\]
\end{itemize}
\end{definition}

\begin{lemma}
\label{lem:ext-2-fun}
Let $U \colon \kCAT \to \kCAT$ be a $2$-functor.
Then the following hold.
\begin{enumerate}
\item
Let $\calC = (\calC, A)$ be a $G$-category.
Define a pair $\hat{U}(\calC) = (\hat{U}(\calC), \hat{U}(A))$ by setting
$\hat{U}(\calC):= U(\calC)$, and $\hat{U}(A)_a:= U(A_a)$
for all $a \in G$.
Then it turns out to be a $G$-category.

\item
Let $(F, \ph) \colon \calC \to \calD$ be a $G$-equivariant functor of $G$-categories.
Define a pair $\hat{U}(F, \ph):= (\hat{U}(F), \hat{U}(\ph))$
by setting $\hat{U}(F):= U(F)$, and $\hat{U}(\ph)_a:= U(\ph_a)$
for all $a\in G$.
Then it turns out to be a $G$-equivariant functor
$\hat{U}(F, \ph) \colon \hat{U}(\calC) \to \hat{U}(\calD)$.

\item
Let $(F, \ph), (F', \ph') \colon \calC \to \calD$ be $G$-equivariant functors of $G$-categories,
and $\al \colon (F, \ph) \To (F', \ph')$ a morphism between them.
Set $\hat{U}(\al):= U(\al)$.
Then it turns out to be a morphism $\hat{U}(\al) \colon \hat{U}(F, \ph) \To \hat{U}(F', \ph')$
of $G$-equivariant functors.

\item By the above,
$U$ is extended to a $2$-functor $\hat{U} \colon \GCAT \to \GCAT$.
\end{enumerate}
\end{lemma}

\begin{proof}
Since the proofs of (1), (2) and (3) are similar, we only show the statement (2) among them.

Set $\calC = (\calC, A)$ and $\calD = (\calD, A')$.
Then in the diagrams below,
$(F,\ph)$ is given by the left one, and it makes the right one commutative
for all $a,b \in G$:
\[
\begin{tikzcd}
\calC & \calD\\
\calC & \calD
\Ar{1-1}{1-2}{"F"}
\Ar{2-1}{2-2}{"F"'}
\Ar{1-1}{2-1}{"A_a"'}
\Ar{1-2}{2-2}{"A'_a"}
\Ar{1-2}{2-1}{"\ph_a"', Rightarrow}
\end{tikzcd},
\quad
\begin{tikzcd}
A'_b A'_a F & A'_bFA_a & FA_b A_a\\
A'_{ba}F && FA_{ba}
\Ar{1-1}{1-2}{"A'_b\ph_a", Rightarrow}
\Ar{1-2}{1-3}{"\ph_b A_a", Rightarrow}
\Ar{2-1}{2-3}{"\ph_{ba}"', Rightarrow}
\Ar{1-1}{2-1}{"", equal}
\Ar{1-3}{2-3}{"", equal}
\end{tikzcd}.
\]
Apply $U$ to these diagrams to have the following:
\[\tiny
\begin{tikzcd}[column sep=30pt]
U(\calC) & U(\calD)\\
U(\calC) & U(\calD)
\Ar{1-1}{1-2}{"U(F)"}
\Ar{2-1}{2-2}{"U(F)"'}
\Ar{1-1}{2-1}{"U(A_a)"'}
\Ar{1-2}{2-2}{"U(A'_a)"}
\Ar{1-2}{2-1}{"U(\ph_a)"', Rightarrow}
\end{tikzcd},
\ 
\begin{tikzcd}[column sep=40pt]
U(A'_b) U(A'_a) U(F) & U(A'_b)U(F)U(A_a) & U(F)U(A_b) U(A_a)\\
U(A'_{ba})U(F) && U(F)U(A_{ba})
\Ar{1-1}{1-2}{"U(A'_b)U(\ph_a)", Rightarrow}
\Ar{1-2}{1-3}{"U(\ph_b) U(A_a)", Rightarrow}
\Ar{2-1}{2-3}{"U(\ph_{ba})"', Rightarrow}
\Ar{1-1}{2-1}{"", equal}
\Ar{1-3}{2-3}{"", equal}
\end{tikzcd}.
\]
The left diagram gives the pair $(\hat{U}(F), \hat{U}(\ph))$,
and the commutativity of the right one shows that it is a $G$-equivariant functor
$\hat{U}(\calC) \to \hat{U}(\calD)$.

(4) Let $(\calC, A)$ be a $G$-category.
Then by definition we have
\[
\hat{U}(\idty_{(\calC, A)}) = \hat{U}(\idty_{\calC}, (\idty_{A_a})_{a\in G})
= (\idty_{U(\calC)}, (\idty_{U(A_a)})_{a\in G}) = \idty_{\hat{U}(\calC, A)}.
\]
Let $(\calC, A) \ya{(F,\ph)} (\calC',A') \ya{(F',\ph')} (\calC'',A'')$ be
$G$-equivariant functors of $G$-categories.
Recall that $(F',\ph') (F, \ph) = (F'F, ((F'\ph_a) \bullet (\ph'_a F))_{a\in G})$,
where $\bullet$ denotes the vertical composition.
Then we have
\[
\begin{aligned}
\hat{U}((F',\ph') (F, \ph)) &= (U(F')U(F), ((U(F')U(\ph_a))\bullet (U(\ph'_a) U(F)))_{a\in G})\\
&= (U(F'),(U(\ph'_a))_{a\in G}) (U(F), (U(\ph_a))_{a\in G})
= \hat{U}(F',\ph') \hat{U}(F, \ph)).
\end{aligned}
\]
Finally, by the definition of $\hat{U}$ on $2$-morphisms, it is obvious that
$\hat{U}$ preserves both vertical and horizontal compositions of $2$-morphisms, and the identities of 2-morphisms.
\end{proof}

To deal with stable module categories,
we need a slight modification as follows.
\begin{definition}
\label{dfn:2ext}
Let $\sfC$ be a $2$-subcategory of $\kCAT$.
Then we define a $2$-subcategory $\hat{\sfC}$ of $\GCAT$ as follows.

First, objects are the $G$-categories $(\calC, A)$ with $\calC \in \sfC$.
Next, for any objects $(\calC, A)$ and $(\calC', A')$, $1$-morphisms
from $(\calC, A)$ to $(\calC', A')$ are
the $G$-equivariant functors $(F,\ph) \colon (\calC, A) \to (\calC', A')$
with $F \in \sfC(\calC, \calC')$.
Finally, for any $1$-morphisms $(F,\ph), (F',\ph') \colon (\calC, A) \to (\calC', A')$
in $\hat{\sfC}$,
$2$-morphisms from $(F,\ph)$ to $(F',\ph')$ are the $2$-morphisms
$(F,\ph) \To (F',\ph')$ in $\GCAT$.
\end{definition}

Using this, Lemma \ref{lem:ext-2-fun} is slightly generalized as follows.

\begin{lemma}
\label{lem:ext-2-fun-gen}
Let $\sfC$ be a $2$-subcategory of $\kCAT$, and
$U \colon \sfC \to \kCAT$ a $2$-functor.
For each object (resp.\ $1$-morphism, $2$-morphism) $x$ in $\sfC$,
set the value of $\hat{U}(x)$ as in Lemma \ref{lem:ext-2-fun}.
Then this defines a $2$-functor $\hat{U} \colon \hat{\sfC} \to \GCAT$.
\end{lemma}

\begin{proof}
The same proof of Lemma \ref{lem:ext-2-fun} works.
\end{proof}

\begin{example}
\label{exm:stable-2-fun}
We give an example of a $2$-functor from a $2$-subcategory of $\kCAT$ to $\kCAT$.
We denote by $\kCAT_p$ the $2$-subcategory of $\kCAT$
whose objects are the light $k$-categories,
whose $1$-morphisms are the $k$-functors $\calA \to \calB$ sending
projective objects of $\calA$ to projective objects of $\calB$,
and whose $2$-morphisms are the natural transformations between $1$-morphisms.
Now we define a $2$-functor
\[
U:= U_p \colon \kCAT_p \to \kCAT
\]
called the {\em stable $2$-functor} as follows.

{\bf On objecs:}
Let $\calA \in \kCAT_p$.
Then $U(\calA):= \udl{\calA}:= \calA/P_\calA$,
where $P_\calA$ is the ideal of $\calA$ consisting of all morphisms in $\calA$
factoring through some projective object of $\calA$.

{\bf On $1$-morphisms:}
Let $F \colon \calA \to \calB$ be in $\kCAT_p$.
Then since $F(P_\calA) \subseteq P_\calB$, it induces a $k$-functor
$U(F):= \udl{F} \colon \udl{\calA} \to \udl{\calB}$.
We denote by $\pi_\calA$ the canonical functor $\calA \to \udl{\calA}$,
and set $\udl{f}:= \pi_\calA(f)$ for all morphisms $f$ in $\calA$.

{\bf On $2$-morphisms:}
Let $F, F' \colon \calA \to \calB$ be in $\kCAT_p$, and
$\al \colon F \To F'$ a $2$-morphism in $\kCAT_p$.
Thus $\al = (\al_x)_{x \in \calA}$ is a natural transformation from $F$ to $F'$
with $\al_x \in \calB(Fx, F'x)$.
Set $U(\al):= \udl{\al}:= (\udl{\al_x})_{x \in \calA}$.
Then $U(\al)$ turns out to be a natural transformation $\udl{F} \To \udl{F'}$.

Then as easily seen, $U \colon \kCAT_p \to \kCAT$ is a $2$-functor.

For $\sfC = \kCAT_p$, we set $\GCAT_p:= \hat{\sfC}$.
Namely, objects of $\GCAT_p$ are the objects of $\GCAT$,
$1$-morphisms are the $G$-equivariant functor $(\calC, A) \to (\calC', A')$
sending projective objects of $\calC$ to projective objectis of $\calC'$,
and $2$-morphisms are the $2$-morphisms in $\GCAT$ between $1$-morphisms
in $\GCAT_p$.
Then by Lemma \ref{lem:ext-2-fun-gen}, $U$ is extended to a $2$-functor
\[
\hat{U}:= \hat{U}_p \colon \GCAT_p \to \GCAT.
\]
\end{example}

To deal with also factor categories and localizations of categories,
we need a further extension as follows.

\begin{definition}\label{dfn:G-stable-2cats}
We define the following two $2$-categories.
\begin{enumerate}
\item 
A $2$-category $\kCAT_s$ is defined as follows:\\
Objects are the pairs $(\calC, \calP)$ of a $k$-category $\calC$ and a subclass
$\calP$ of the morphisms class of $\calC$, which are called
{\em structured $k$-categories}.\\
For any objects $(\calC, \calP)$ and $(\calC', \calP')$,
$1$-morphisms from the former to the latter are the
functors $F \colon \calC \to \calC'$ such that $F(\calP) \subseteq \calP'$.
In particular, an automorphism of $(\calC, \calP)$ is
an automorphism $A$ of $\calC$
satisfying $A(\calP) \subseteq \calP$.\\
For any $1$-morphisms $F, F' \colon (\calC, \calP) \to (\calC', \calP')$,
$2$-morphisms from $F$ to $F'$ are the natural transformations $F \To F'$.

\item
A $2$-category $\GCAT_s$ is defined as follows:\\
Objects are the pairs $((\calC, A), \calP)$ (or triples $(\calC, A, \calP)$),
where $(\calC,\calP)$ is a structured
$k$-category and $A \colon G \to \Aut(\calC, \calP)$ is a group homomorphism,
or equivalently, $(\calC, A)$ is a $G$-category and $\calP$ is a {\em $G$-stable} subclass
of the morphism class of $\calC$ in the sense that
$A_a(\calP) \subseteq \calP$ for all $a \in G$.
$(\calC, A, \calP)$ are called
{\em structured $G$-categories}.\\
\hspace{1mm}For any objects $(\calC, A, \calP)$ and $(\calC', A', \calP')$,
$1$-morphisms from the former to the latter are the
$1$-morphisms $(F, \ph) \colon (\calC, A) \to (\calC',A')$ in $\GCAT$
such that $F(\calP) \subseteq \calP'$.\\
\hspace{1mm}For any $1$-morphisms $(F,\ph), (F', \ph') \colon (\calC, A, \calP) \to (\calC', A', \calP')$,
$2$-morphisms from the former to the latter are the $2$-morphisms $(F,\ph) \To (F',\ph')$
in $\GCAT$.
\end{enumerate}
\end{definition}

\begin{definition}
\label{dfn:2ext-gen}
Let $\sfC$ be a $2$-subcategory of $\kCAT_s$.
We define a $2$-subcategory $\hat{\sfC}$ of $\GCAT_s$ as follows.

First, objects are the objects $(\calC, A, \calP)$ in $\GCAT_s$ with $(\calC, \calP) \in \sfC$.
Next, for any objects $(\calC, A, \calP)$ and $(\calC', A', \calP')$, $1$-morphisms
from the former to the latter are
the $G$-equivariant functors $(F,\ph) \colon (\calC, A) \to (\calC', A')$
with $F \in \sfC((\calC, \calP), (\calC', \calP'))$.

Finally, for any $1$-morphisms
$(F,\ph), (F',\ph') \colon (\calC, \calP, A) \to (\calC', \calP', A')$
in $\hat{\sfC}$,
$2$-morphisms from $(F,\ph)$ to $(F',\ph')$ are the $2$-morphisms
$(F,\ph) \To (F',\ph')$ in $\GCAT$.
\end{definition}

\begin{lemma}
\label{lem:ext-2-fun-gen+}
Let $\sfC$ be a $2$-subcategory of $\kCAT_s$, and
$U \colon \sfC \to \kCAT$ a $2$-functor.
For each object (resp.\ $1$-morphism, $2$-morphism) $x$ in $\sfC$,
set the value of $\hat{U}(x)$ as in Lemma \ref{lem:ext-2-fun}.
Then this defines a $2$-functor $\hat{U} \colon \hat{\sfC} \to \GCAT$.
\end{lemma}

\begin{proof}
The same proof of Lemma \ref{lem:ext-2-fun} works.
\end{proof}

If there seems to be no confusion, we denote $\hat{U}$ just by $U$.
Then we have $U(A)_a = U(A_a)$ and $U(\ph)_a = U(\ph_a)$
for all $G$-categories $(\calC, A)$ and $G$-equivariant functors $(F, \ph)$.

We extend the construction in Lemma \ref{lem:ext-2-fun-gen}
to {\em strict $2$-natural transformations}
(see \cite[Definition 4.1.14]{Asa-book}).

\begin{lemma}
\label{lem:ext-2-nat}
Consider the following diagram of $2$-categories, $2$-functors,
and a strict $2$-natural transformation:
\[
\begin{tikzcd}
\sfC & \makebox[1em][l]{$\kCAT$}
\Ar{1-1}{1-2}{"V", ""'{name=a}, bend left}
\Ar{1-1}{1-2}{"U"', ""{name=b}, bend right}
\Ar{a}{b}{"\ro", Rightarrow}
\end{tikzcd}\qquad,
\]
where $\sfC$ is a $2$-subcategory of $\kCAT$
(resp.\ $\kCAT_s$).
Then it is extended to the following diagram of $2$-categories, $2$-functors,
and a strict $2$-natural transformation:
\[
\begin{tikzcd}
\hat{\sfC} & \makebox[1em][l]{$\GCAT$}
\Ar{1-1}{1-2}{"\hat{V}", ""'{name=a}, bend left}
\Ar{1-1}{1-2}{"\hat{U}"', ""{name=b}, bend right}
\Ar{a}{b}{"\hat{\ro}", Rightarrow}
\end{tikzcd}\qquad,
\]
where $\hat{\sfC}, \hat{U}, \hat{V}$ are defined as in Lemma \ref{lem:ext-2-fun-gen}.
\end{lemma}

\begin{proof}
We show the assertion in the case that $\sfC$ is a $2$-subcategory of $\kCAT$.
The assertion in the remaining case is shown similarly.
Since $\ro$ is a strict $2$-natural transformation,
any diagram in $\sfC$ on the left yields a commutative diagram on the right:
\[
\begin{tikzcd}[row sep=50pt, column sep=60pt]
\calC & \calC'
\Ar{1-1}{1-2}{"F", ""'{name=a}, bend left}
\Ar{1-1}{1-2}{"F'"', ""{name=b}, bend right}
\Ar{a}{b}{"\al", Rightarrow}
\end{tikzcd},
\quad
\begin{tikzcd}[row sep=50pt, column sep=60pt]
V(\calC) & V(\calC')\\
U(\calC) & U(\calC')
\Ar{1-1}{1-2}{"V(F)", ""'{name=a}, bend left}
\Ar{1-1}{1-2}{"V(F')"', ""{name=b}, bend right}
\Ar{a}{b}{"V(\al)", Rightarrow}
\Ar{2-1}{2-2}{"U(F)", ""'{name=aa}, bend left}
\Ar{2-1}{2-2}{"U(F')"', ""{name=bb}, bend right}
\Ar{aa}{bb}{"U(\al)", Rightarrow}
\Ar{1-1}{2-1}{"\ro_\calC"'}
\Ar{1-2}{2-2}{"\ro_{\calC'}"}
\end{tikzcd}.
\]

Thus we have
\begin{align}
\label{eq:strict-2-nat-1}
\ro_{\calC'}V(F) &= U(F) \ro_\calC, \text{ and}\\
\label{eq:strict-2-nat-2}
\ro_{\calC'}V(\al) &= U(\al) \ro_\calC.
\end{align}
First, we define a candidate of a strict $2$-natural transformation
$\hat{\ro} \colon \hat{V} \To \hat{U}$, which should be a family
$\hat{\ro}:= (\hat{\ro}_{(\calC, A)})_{(\calC, A)\in \hat{\sfC}}$
with
\[
\hat{\ro}_{(\calC, A)} \colon (V(\calC), (V(A_a))_{a\in G}) \to (U(\calC), (U(A_a))_{a\in G})
\]
a $G$-equivariant functor.
Thus we are looking for a $2$-morphism in the diagram
\[
\begin{tikzcd}
V(\calC) & U(\calC)\\
V(\calC) & U(\calC)
\Ar{1-1}{1-2}{"\ro_\calC"}
\Ar{2-1}{2-2}{"\ro_\calC"'}
\Ar{1-1}{2-1}{"V(A_a)"'}
\Ar{1-2}{2-2}{"U(A_a)"}
\Ar{1-2}{2-1}{Rightarrow, dashed}
\end{tikzcd}.
\]
Apply \eqref{eq:strict-2-nat-1} to $A_a \colon \calC \to \calC$ to have
$\ro_{\calC}V(A_a) = U(A_a) \ro_\calC$ for all $a \in G$, 
which makes it possible to define $\hat{\ro}$
by 
\begin{equation}
\label{eq:ext-nat-tr}
\hat{\ro}_{(\calC, A)}:= (\ro_\calC, (\idty_{U(A_a)\ro_\calC})_{a\in G}).
\end{equation}

We now show that $\hat{\ro}$ is a strict $2$-natural transformation,
namely that
any diagram in $\hat{\sfC}$ on the left yields a commutative diagram on the right:
\[
\begin{tikzcd}[row sep=50pt, column sep=60pt]
(\calC,A) & (\calC',A')
\Ar{1-1}{1-2}{"{(F,\ph)}", ""'{name=a}, bend left}
\Ar{1-1}{1-2}{"{(F',\ph')}"', ""{name=b}, bend right}
\Ar{a}{b}{"\al", Rightarrow}
\end{tikzcd},
\quad
\begin{tikzcd}[row sep=50pt, column sep=60pt]
\hat{V}(\calC,A) & \hat{V}(\calC',A')\\
\hat{U}(\calC,A) & \hat{U}(\calC',A')
\Ar{1-1}{1-2}{"{\hat{V}(F,\ph)}", ""'{name=a}, bend left}
\Ar{1-1}{1-2}{"{\hat{V}(F',\ph')}"', ""{name=b}, bend right}
\Ar{a}{b}{"\hat{V}(\al)", Rightarrow}
\Ar{2-1}{2-2}{"{\hat{U}(F,\ph)}", ""'{name=aa}, bend left}
\Ar{2-1}{2-2}{"{\hat{U}(F',\ph')}"', ""{name=bb}, bend right}
\Ar{aa}{bb}{"\hat{U}(\al)", Rightarrow}
\Ar{1-1}{2-1}{"{\hat{\ro}_{(\calC,A)}}"'}
\Ar{1-2}{2-2}{"{\hat{\ro}_{(\calC',A')}}"}
\end{tikzcd}.
\]
Thus we have to show that
\begin{align}
\label{eq:strict-2-nat-1a}
\hat{\ro}_{(\calC',A')}\hat{V}(F,\ph) &= \hat{U}(F,\ph) \hat{\ro}_{(\calC,A)}, \text{ and}\\
\label{eq:strict-2-nat-2a}
\hat{\ro}_{(\calC',A')}\hat{V}(\al) &= \hat{U}(\al) \hat{\ro}_{(\calC,A)}.
\end{align}
By definition, \eqref{eq:strict-2-nat-2a}
coincides with \eqref{eq:strict-2-nat-2}.
The left hand side of \eqref{eq:strict-2-nat-1a} is equal to the following:
\[
\begin{aligned}
(\ro_{\calC'}, (\idty_{U(A_a)\ro_{\calC'}})_{a}) (V(F), (V(\ph_a))_{a})
&= (\ro_{\calC'}V(F), ((\ro_{\calC'}V(\ph_a))\bullet (\idty_{U(A_a)\ro_{\calC'}}))_{a})\\
&= (\ro_{\calC'}V(F), (\ro_{\calC'}V(\ph_a))_{a}).
\end{aligned}
\]
The right hand side of \eqref{eq:strict-2-nat-1a} is equal to the following:
\[
\begin{aligned}
(U(F), (U(\ph_a))_{a})(\ro_\calC, (\idty_{U(A_a)\ro_{\calC}})_{a})
&= (U(F)\ro_\calC, ((U(F)\idty_{U(A_a)\ro_{\calC}})\bullet (U(\ph_a)\ro_{\calC}))_{a})\\
&= (U(F)\ro_\calC, (U(\ph_a)\ro_{\calC})_{a}).
\end{aligned}
\]
We now compare these last values.
By \eqref{eq:strict-2-nat-1}, the first components coincide.
By applying \eqref{eq:strict-2-nat-2} to $\al:= \ph_a$ for all $a \in G$, we see
that the second components coincide as well.
Hence \eqref{eq:strict-2-nat-1a} holds.
\end{proof}

\begin{example}
\label{exm:stable-2-nat}
In Example \ref{exm:stable-2-fun}, we note that the family
$\pi:= (\pi_\calC)_{\calC \in \kCAT_p}$
defines a strict $2$-natural transformation
\[
\begin{tikzcd}
\makebox[1em][r]{$\kCAT_p$} & \makebox[1em][l]{$\kCAT$}
\Ar{1-1}{1-2}{"\si", ""'{name=a}, bend left}
\Ar{1-1}{1-2}{"U"', ""{name=b}, bend right}
\Ar{a}{b}{"\pi", Rightarrow}
\end{tikzcd}\qquad,
\]
where $\si$ is the inclusion $2$-functor.
Then by Example \ref{exm:stable-2-fun} and Lemma \ref{lem:ext-2-nat}, the family
\[
\hat{\pi}:= (\pi_\calC, (\idty_{\pi_\calC A_a})_{a\in G})_{(\calC,A) \in \GCAT_p}
\]
defines a strict $2$-natural transformation
\begin{equation}
\label{eq:stable-2-nat}
\begin{tikzcd}
\makebox[1em][r]{$\GCAT_p$} & \makebox[1em][l]{$\GCAT$}
\Ar{1-1}{1-2}{"\hat{\si}", ""'{name=a}, bend left}
\Ar{1-1}{1-2}{"\hat{U}"', ""{name=b}, bend right}
\Ar{a}{b}{"\hat{\pi}", Rightarrow}
\end{tikzcd}\qquad.
\end{equation}
Note that $\hat{\si}$ is nothing but the inclusion $2$-functor.
\end{example}

For applications, we now add the following two examples of $2$-functors.

\begin{example}
\label{exm:factor-2-fun}

We define a $2$-category $\kCAT_i$ as a $2$-subcategory of $\kCAT_s$ as follows.

{\bf Objects:}
Objects are the pairs $(\calC, S)$, where $\calC \in \kCAT$,
and $S$ is an ideal of $\calC$.

{\bf 1-morphims:}
Let $(\calA, S), (\calB, T)$ be objects.
Then $1$-morphisms from $(\calA, S)$ to $(\calB, T)$
are the $k$-functors $F \colon \calA \to \calB$ such that
$F(S) \subseteq T$.

{\bf $2$-morphims:}
Let $(\calA, S), (\calB, T)$ be objects, and
$E, F \colon (\calA, S) \to (\calB, T)$ $1$-morphisms.
Then $2$-morphisms from $E$ to $F$ are the natural transformations
from $E$ to $F$.

Now we define a $2$-functor
\[
U:= U_i \colon \kCAT_i \to \kCAT
\]
called the {\em factor $2$-functor} as follows.

{\bf On objecs:}
Let $(\calA,S) \in \kCAT_i$.
Then $U(\calA,S):= \calA/S$.

{\bf On $1$-morphisms:}
Let $F \colon (\calA, S) \to (\calB, T)$ be in $\kCAT_m$.
Then since $F(S) \subseteq T$, it induces a $k$-functor
$U(F) \colon \calA/S \to \calB/T$.
We denote by $\pi_{(\calA,S)}$ the canonical functor $\calA \to \calA/S$,
and set $\ovl{f}:= \pi_{(\calA,S)}(f)$ for all morphisms $f$ in $\calA$.
Note here that $\pi_{(\calA,S)}\colon \calA \to \calA/S$
is a full functor for all $(\calA,S) \in \kCAT_i$.

{\bf On $2$-morphisms:}
Let $E, F \colon (\calA,S) \to (\calB,T)$ be in $\kCAT_i$, and
$\al \colon E \To F$ a $2$-morphism in $\kCAT_i$.
Thus $\al = (\al_x)_{x \in \calA}$ is a natural transformation from $E$ to $F$
with $\al_x \in \calB(Ex, Fx)$.
Set $U(\al):= (\ovl{\al_x})_{x \in \calA}$.
Then $U(\al)$ turns out to be a natural transformation $U(E) \To U(F)$.

Then as easily seen, $U \colon \kCAT_i \to \kCAT$ is a $2$-functor.

For $\sfC = \kCAT_i$, we set $\GCAT_i:= \hat{\sfC}$.
Namely, objects of $\GCAT_i$ are the pairs $(\calC, S)$,
where $\calC$ is an object of $\GCAT$, and $S$ is an ideal of $\calC$,
$1$-morphisms $((\calA, A),S) \to ((\calB, B),T)$
are the $G$-equivariant functors $(F, \ph)\colon (\calA, A) \to (\calB, B)$ with  $F(S) \subseteq T$,
and $2$-morphisms are the $2$-morphisms in $\GCAT$ between the $1$-morphisms in $\GCAT_i$.
Then by Lemma \ref{lem:ext-2-fun-gen+}, $U$ is extended to a $2$-functor
\[
\hat{U}:= \hat{U}_i \colon \GCAT_i \to \GCAT.
\]
As in Example \ref{exm:stable-2-nat},
the family
$\pi_i:= \pi:= (\pi_{(\calA,S)}\colon \calA \to \calA/S)_{(\calA,S) \in \kCAT_i}$
defines a strict $2$-natural transformation
\[
\begin{tikzcd}
\makebox[1em][r]{$\kCAT_i$} & \makebox[1em][l]{$\kCAT$}
\Ar{1-1}{1-2}{"\si_i", ""'{name=a}, bend left}
\Ar{1-1}{1-2}{"U_i"', ""{name=b}, bend right}
\Ar{a}{b}{"\pi_i", Rightarrow}
\end{tikzcd}\qquad,
\]
where $\si_i$ is the forgetful $2$-functor
$((\calA, S)\ya{F} (\calB,T)) \mapsto (\calA \ya{F} \calB)$.
Then by Lemma \ref{lem:ext-2-nat}, the family
\[
\hat{\pi}_i:= \hat{\pi}:= (\pi_{(\calA,S)}, (\idty_{(\pi_{(\calA,S)} A_a)})_{a\in G})_{(\calA,A,S) \in \GCAT_i}
\]
defines a strict $2$-natural transformation
\begin{equation}
\label{eq:factor-2-nat}
\begin{tikzcd}
\makebox[1em][r]{$\GCAT_i$} & \makebox[1em][l]{$\GCAT$}
\Ar{1-1}{1-2}{"\hat{\si}_i", ""'{name=a}, bend left}
\Ar{1-1}{1-2}{"\hat{U}_i"', ""{name=b}, bend right}
\Ar{a}{b}{"\hat{\pi}_i", Rightarrow}
\end{tikzcd}\qquad,
\end{equation}
where $\hat{\si}_i$ is the forgetful $2$-functor
\[
((\calA, A, S)\ya{(F,\ph)}(\calB, B, T)) \mapsto ((\calA, A) \ya{(F,\ph)} (\calB, B)).
\]
\end{example}

\begin{remark}
We remark that $\kCAT_p$ is related to a 2-subcategory of $\kCAT_i$.
Let $\kCAT_P$ be a $2$-subcategory of $\kCAT_i$
whose objects are the pairs $(\calA, P_\calA)$ of an object
$\calA$ of $\kCAT$ and the ideal $P_\calA$ of $\calA$
consisting of all morphisms in $\calA$
factoring through some projective object of $\calA$,
and whose morphisms are given by
$\kCAT_P((\calA, P_\calA), (\calB, P_\calB)):=
\kCAT_i((\calA, P_\calA), (\calB, P_\calB))$
for all objects $(\calA, P_\calA), (\calB, P_\calB)$
of $\kCAT_P$.
Then we have
\[
\kCAT_p(\calA, \calB) \subseteq \kCAT_P((\calA, P_\calA), (\calB, P_\calB)),
\]
more precisely,
\[
\footnotesize
\kCAT_p(\calA, \calB) = \{F \in \kCAT_P((\calA, P_\calA), (\calB, P_\calB))
\mid F \text{sends projectives to projectives}\}
\]
and
$U_p(\calA) = \calA/P_\calA = U_i(\calA, P_\calA)$
for all objects $\calA, \calB$ of $\kCAT$.
\end{remark}

\begin{example}[Quotient case]
\label{exm:quot-2-fun}
We define a $2$-category $\kCAT_m$ as a $2$-subcategory of $\kCAT_s$ as follows.

{\bf Objects:}
Objects are the pairs $(\calC, S)$, where $\calC \in \kCAT$,
and $S$ is a bicalculable multiplicative system of $\calC$.

{\bf 1-morphims:}
Let $(\calA, S), (\calB, T)$ be objects.
Then $1$-morphisms from $(\calA, S)$ to $(\calB, T)$
are the $k$-functors $F \colon \calA \to \calB$ such that
$F(S) \subseteq T$.

{\bf $2$-morphims:}
Let $(\calA, S), (\calB, T)$ be objects, and
$E, F \colon (\calA, S) \to (\calB, T)$ $1$-morphisms.
Then $2$-morphisms from $E$ to $F$ are the natural transformations
from $E$ to $F$.

Now we define a $2$-functor
\[
U:= U_m \colon \kCAT_m \to \kCAT
\]
called the {\em quotient $2$-functor} (or {\em localization $2$-functor}) as follows.

{\bf On objecs:}
Let $(\calA,S) \in \kCAT_m$.
Then $U(\calA,S):= \calA[S\inv]$.

{\bf On $1$-morphisms:}
Let $F \colon (\calA, S) \to (\calB, T)$ be in $\kCAT_m$.
Then since $F(S) \subseteq T$, it induces a $k$-functor
$U(F) \colon \calA[S\inv] \to \calB[T\inv]$.
We denote by $\pi_{(\calA,S)}$ the quotient functor $\calA \to \calA[S\inv]$,
and set $\ovl{f}:= \pi_{(\calA,S)}(f)$ for all morphisms $f$ in $\calA$.

{\bf On $2$-morphisms:}
Let $E, F \colon (\calA,S) \to (\calB,T)$ be in $\kCAT_m$, and
$\al \colon E \To F$ a $2$-morphism in $\kCAT_m$.
Thus $\al = (\al_x)_{x \in \calA}$ is a natural transformation from $E$ to $F$
with $\al_x \in \calB(Ex, Fx)$.
Set $U(\al):= (\ovl{\al_x})_{x \in \calA}$.
Then $U(\al)$ turns out to be a natural transformation $U(E) \To U(F)$.

Then as easily seen, $U \colon \kCAT_m \to \kCAT$ is a $2$-functor.

For $\sfC = \kCAT_m$, we set $\GCAT_m:= \hat{\sfC}$.
Namely, objects of $\GCAT_m$ are the pairs $(\calC, S)$,
where $\calC$ is an object of $\GCAT$, and $S$ is a bicalculable multiplicative system of $\calC$,
$1$-morphisms $(\calA, A,S) \to (\calB, B,T)$
are the $G$-equivariant functors $(F,\ph) \colon (\calA, A) \to (\calB, B)$ with $F(S) \subseteq T$,
and $2$-morphisms are the $2$-morphisms in $\GCAT$ between $1$-morphisms in $\GCAT_m$.
Then by Lemma \ref{lem:ext-2-fun-gen+}, $U$ is extended to a $2$-functor
\[
\hat{U}:= \hat{U}_m \colon \GCAT_m \to \GCAT.
\]
As in Example \ref{exm:factor-2-fun},
the family
$\pi_m:= \pi:= (\pi_{(\calA,S)}\colon \calA \to \calA[S\inv])_{(\calA,S) \in \kCAT_m}$
defines a strict $2$-natural transformation
\[
\begin{tikzcd}
\makebox[1em][r]{$\kCAT_m$} & \makebox[1em][l]{$\kCAT$}
\Ar{1-1}{1-2}{"\si_m", ""'{name=a}, bend left}
\Ar{1-1}{1-2}{"U_m"', ""{name=b}, bend right}
\Ar{a}{b}{"\pi_m", Rightarrow}
\end{tikzcd}\qquad,
\]
where $\si_m$ is the forgetful $2$-functor
$((\calA, S)\ya{F} (\calB,T)) \mapsto (\calA \ya{F} \calB)$.
Then by Lemma \ref{lem:ext-2-nat}, the family
\[
\hat{\pi}_m:= \hat{\pi}:= (\pi_{(\calA,S)}, (\idty_{(\pi_{(\calA,S)} A_a)})_{a\in G})_{(\calA,S,A) \in \GCAT_m}
\]
defines a strict $2$-natural transformation
\begin{equation}
\label{eq:stable-2-nat-m}
\begin{tikzcd}
\makebox[1em][r]{$\GCAT_m$} & \makebox[1em][l]{$\GCAT$}
\Ar{1-1}{1-2}{"\hat{\si}_m", ""'{name=a}, bend left}
\Ar{1-1}{1-2}{"\hat{U}_m"', ""{name=b}, bend right}
\Ar{a}{b}{"\hat{\pi}_m", Rightarrow}
\end{tikzcd}\qquad,
\end{equation}
where $\hat{\si}_m$ is the forgetful $2$-functor
$((\calA, A, S)\ya{(F,\ph)}(\calB, B, T)) \mapsto ((\calA, A) \ya{(F,\ph)} (\calB, B))$.
\end{example}

\begin{remark}[Dense/localizing subcategories]
\label{rmk:dense-localizing}
For instance, bicalculable multiplicative systems can be obtained as follows.
Let $\calC$ be an abelian category.
Recall that a full subcategory $\calS$ of $\calC$ is said to be {\em dense} (or a {\em Serre subcategory}) if it is closed under subobjects, factor objects, and extensions.
If this is the case, then
$S:= S_\calC(\calS):= \{f: \text{ a morphism in } \calC \mid \Ker f, \Coker f \in \calS\}$ becomes a bicalculable multiplicative system of $\calC$.
In this case, we set
$\calC/\calS:= \calC[S\inv]$
(see e.g.\ \cite[Ch.\ 4, Sect.\ 4.1]{Po73}).
Recall also that $\calS$ is called {\em localyzing}
if the canonical functor $Q \colon \calC \to \calC/\calS$ has
a right adjoint.
This is equivalent to saying that $\calS$ is closed under small coproducts
if $\calC$ has small coproducts and injective envelopes
(see e.g.\ \cite[Ch.\ 4, Proposition 6.3]{Po73}).
\end{remark}

\section{$G$-precoverings induced by adjoint}

In this section, we prove a theorem that produces new $G$-precoverings from old.
First, we review the notion of the adjoint pair in 2-categories.

\begin{definition}
\label{dfn:adj-sys}
Let $\sfC$ be a $2$-category, and consider the diagram
\[\begin{tikzcd}[row sep= 10pt, column sep=10pt]
\calC &&& \calD \\
{}&{}&{}&{}\\
\calC &&& \calD
\Ar{1-1}{1-4}{"L"}
\Ar{1-4}{3-1}{"R" description}
\Ar{3-1}{3-4}{"L"'}
\Ar{1-1}{3-1}{"\idty_\calC"'}
\Ar{1-4}{3-4}{"\idty_\calD"}
\Ar{2-1}{2-2}{"\et", Rightarrow}
\Ar{2-3}{2-4}{"\ep"', Rightarrow}
\end{tikzcd}
\]
in $\sfC$.
Then the quadruple $(L,R, \et, \ep)$ is called an {\em adjoint system}
from $\calC$ to $\calD$ if the following diagrams are commutative:
\[
\begin{tikzcd}
L & L\circ R \circ L\\
& L
\Ar{1-1}{1-2}{"\idty_L \circ \et", Rightarrow}
\Ar{1-2}{2-2}{"\ep \circ \idty_L", Rightarrow}
\Ar{1-1}{2-2}{equal}
\end{tikzcd} \text{ and }
\begin{tikzcd}
R & R\circ L \circ R\\
& R
\Ar{1-1}{1-2}{"\et\circ \idty_R", Rightarrow}
\Ar{1-2}{2-2}{"\idty_R\circ \ep ", Rightarrow}
\Ar{1-1}{2-2}{equal}
\end{tikzcd}.
\]
(Here, $L\et = \idty_L \circ \et$ and $\ep_L = \ep \circ \idty_L$ are the horizontal composites.)
When this is the case, $L$ is called a {\em left adjoint} to $R$,
$R$ is called a {\em right adjoint} to $L$, and $(L,R)$ is called an {\em adjoint pair}.
The $2$-morphism $\et$ (resp.\ $\ep$) is called a \emph{unit}
(resp.\ a \emph{counit}) of this adjoint system.

In the case that $\sfC =  \GCAT$, the family
\[
\om = (\om_{x,y} \colon \calD(Lx,y) \to \calC(x,Ry))_{x \in \calC, y \in \calD}
\]
defined by $f \mapsto Rf \circ \et_x$ is well-known to be a natural isomorphism
(see e.g., \cite[Theorem 4.4.7]{Asa-book}), 
which is called the {\em adjunction} of this adjoint system,
and we denote this by $\om \colon L \adj R$.
\end{definition}
Since a $2$-functor preserves identity $2$-morphisms and
both vertical and horizontal compositions of $2$-morphisms,
the following is obvious.

\begin{lemma}
\label{lem:2-fun-pres-adj}
Let $U \colon \sfC \to \sfD$ be a $2$-functor between $2$-categories,
and $(L, R, \et, \ep)$ an adjoint system from $\calC$ to $\calD$ in $\sfC$.
Then $(U(L), U(R), U(\et), U(\ep))$ turns out to be
an adjoint system from $U(\calC)$ to $U(\calD)$ in $\sfD$.

In particular, when $\sfC = \GCAT$, if
$\om = (\om_{x,y} \colon \calD(Lx,y) \to \calC(x,Ry))_{x \in \calC, y \in \calD}$
is the adjunction of $(L, R, \et, \ep)$, then the adjunction
\[
U(\om) = (U(\om)_{x,y}\colon U(\calD)(U(L)x,y) \to U(\calC)(x,U(R)y))_{x \in U(\calC), y \in U(\calD)}
\]
of $(U(L), U(R), U(\et), U(\ep))$ is given by
$f \mapsto U(R)f \circ U(\et)_x$.
\qed
\end{lemma}

\subsection{Preservation of compactness}
\label{ssec:compact}

Let $\calA$ be a $k$-category and $x \in \calA$.
Consider a cocone $(y, (s_i \colon y_i \to y)_{i\in I})$ in $\calA$.
Then we have a cocone
\[
(\calA(x, y),\ (\calA(x, s_i) \colon \calA(x, y_i) \to \calA(x, y))_{i\in I})
\]
in $\Mod k$.
Since $\Mod k$ has small coproducts, we have the following commutative diagram:
\[
\begin{tikzcd}
\calA(x, y_i) && \calA(x, y)\\
&\Ds_{i\in I} \calA(x, y_i)
\Ar{1-1}{1-3}{"{\calA(x,s_i)}"}
\Ar{1-1}{2-2}{"u_i"'}
\Ar{2-2}{1-3}{"\nu"'}
\end{tikzcd},
\]
where $u_i$ is the canonical injection defined by $f \mapsto (\de_{j,i}f)_{j \in I}$
for all $i\in I$, and $\nu$ is the canonical morphism
defined by
\begin{equation}
\label{eq:compact}
\nu((f_i)_{i \in I}):= \sum_{i \in I}s_i\circ f_i
\quad\text{for all}\ (f_i)_{i \in I} \in \Ds_{i \in I}\calA(x, y_i).
\end{equation}
We set $[\calA(x,s_i)]_{i \in I}:= \nu$.

With this notation, let us recall the following definition.

\begin{definition}
\label{dfn:compact}
An object $x$ in a $k$-category $\calA$ is said to be {\em compact}
if for any small set $I$ and any family $(y_i)_{i\in I}$ of objects of $\calA$,
the functor $\calA(x, \blank)$ {\em preserves} the coproduct of $(y_i)_{i\in I}$,
in the sense that
for any coproduct $(y, (s_i \colon y_i \to y)_{i\in I})$ of $(y_i)_{i\in I}$,
the canonical map $[\calA(x,s_i)]_{i \in I} \colon \Ds_{i\in I} \calA(x, y_i) \to \calA(x, y)$
is an isomorphism of $k$-vector spaces.

We denote by $\cpt(\calA)$ the full subcategory of $\calA$
consisting of the compact objects.
For example, if $\calA = \Mod \calC$ for some $k$-category $\calC$,
then $\cpt(\calA) = \mmod \calC$.
\end{definition}

The following remark is useful when we show the compactness of an object.

\begin{remark}
\label{rmk:compact}
In the above definition,
it is well-known that $[\calA(x,s_i)]_{i \in I}$ is an injection in general (see e.g \cite[Lemma 1.5.17 (3)]{Asa-book}).
Note also that a coproduct is uniquely determined up to natural isomorphism.
Hence to prove that $x$ is compact, it is enough to show that there exists a coproduct
$(y, (s_i \colon y_i \to y)_{i\in I})$ of $(y_i)_{i\in I}$
such that the canonical map $[\calA(x,s_i)]_{i \in I}$ is a surjection.
\end{remark}

\begin{lemma}
\label{lem:left-adj-cpt}
Let $(L, R, \et, \ep)$ be an adjoint system from $\calA$ to $\calB$ in a $2$-category $\sfC$,
which is a $2$-subcategory of $\kCAT$.
If $R$ preserves small coproducts (for instance, if $R$ has a right adjoint), then
\[
L(\cpt(\calA)) \subseteq \cpt(\calB).
\]
\end{lemma}

\begin{proof}
Set $\om$ to be the adjunction isomorphism
of the adjoint system, and
let $x \in \cpt(\calA)$.
To show that $L(x)$ is compact in $\calB$,
let $I$ be  a small set, and $(y, (\si_i \colon y_i \to y)_{i \in I})$ a
coproduct of a family $(y_i)_{i\in I}$ of objects of $\calB$.
Then by the universality of coproducts, there exists a unique map
$\th \colon \Ds_{i\in I}\calB(F(x), y_i) \to \calB(L(x),y)$
and a unique map
$\et \colon \Ds_{i\in I}\calA(x, R(y_i)) \to \calA(x,R(y))$
such that the upper and lower triangles in the diagram
\[
\begin{tikzcd}
\calB(L(x), y_i) && \calB(L(x), y)\\
& \Ds_{i\in I}\calB(L(x), y_i)\\
\calA(x, R(y_i)) && \calA(x, R(y))\\
& \Ds_{i\in I}\calA(x, R(y_i))
\Ar{1-1}{1-3}{"{\calB(F(x), \si_i)}"}
\Ar{1-1}{2-2}{"v_i"'}
\Ar{2-2}{1-3}{"\th"', dashed}
\Ar{3-1}{3-3}{"{\calA(x, R(\si_i))}", near start}
\Ar{3-1}{4-2}{"u_i"'}
\Ar{4-2}{3-3}{"\et"', dashed}
\Ar{1-1}{3-1}{"\om_{x,y_i}"', "\rotatebox{90}{$\sim$}"}
\Ar{1-3}{3-3}{"\om_{x,y}", "\rotatebox{90}{$\sim$}"'}
\Ar{2-2}{4-2}{"\Ds_{i\in I}\om_{x,y_i}", "\rotatebox{90}{$\sim$}"', near start, crossing over}
\end{tikzcd}
\]
commutes, where $u_i$ and $v_i$ are the canonical injections.
To show that $L(x) \in \cpt(\calB)$, it is enough to show that $\th$ is bijective.
Here, the back square commutes by the naturality of $\om$,
and so does the left square by definitions of $u_i, v_i$.
Set $\al$ (resp.\ $\be$) to be the clockwise (resp.\ counter-clockwise)
composite of the right square.
The commutativities of these squares and two triangles
show that $\al v_i = \be v_i$ for all $i \in I$.
Therefore by the universality of coproduct, we have $\al = \be$,
that is, the right square commutes.
Now since $R$ preserves coproducts, $(R(y), R(\si_i)_{i\in I})$ is a coproduct of
$(R(y_i))_{i\in I}$.
Therefore, since $x \in \cpt(\calA)$, $\et$ is bijective,
and hence so is $\th$.
\end{proof}

\begin{remark}
\label{rmk:autoeq-cpt}
Let $\calA$ and $\calB$ be in a $2$-category $\sfC$,
which is a $2$-subcategory of $\kCAT$.
If $L \colon \calA \to \calB$ is an equivalence with a quasi-inverse $R$,
then it is well known that there exists an adjoint system $(L, R, \et, \ep)$
from $\calA$ to $\calB$ (see e.g., \cite[Corollary 4.4.5]{Asa-book}).
Hence by Lemma \ref{lem:left-adj-cpt}, $L(\cpt(\calA)) \subseteq \cpt(\calB)$ because $R$ has a right (and a left) adjoint $L$.
In particular, if $L$ is an auto-equivalence of $\calA$,
then $L(\cpt(\calA)) \subseteq \cpt(\calA)$.
\end{remark}


\subsection{$G$-precoverings induced by adjoint}

By referring to Proposition \ref{Vitalcommutativediagram}, a property of $(P\down, \ph\down)$,
we define the following.

\begin{definition}
\label{dfn:precov-adj}
Let $\calA = (\calA, A)$ be a $G$-category having small coproducts%
\footnote{
This means that $\calA$ has small coproducts.
Note that we use the direct sum $\Ds_{a\in G}ay$ in the diagram \eqref{eq:precov-adj}, which requires this assumption.}
, $\calB$ a $k$-category
and $(L, \ph) \colon \calA \to \calB$ a $G$-invariant functor, where $L$
has a right adjoint $R \colon \calB \to \calA$
with an adjunction isomorphism $\om \colon L \adj R$.
Then $(L, \ph)$ is said {\em to induce a $G$-precovering}
$(L, \ph) \colon \cpt(\calA) \to \cpt(\calB)$ {\em by the adjunction} $\om$
if the following conditions are satisfied:
\begin{enumerate}
\item[{\rm (PA1)}]
There exists a natural isomorphism $t \colon \Ds_{a\in G}A_a \To RL$.
\item[{\rm (PA2)}]
For any $x,y \in \calA$, the following diagram commutes:
\begin{equation}\label{eq:precov-adj}
\begin{tikzcd}
\Ds_{a \in G} \calA(x, ay) & \calA(x, \Ds_{a \in G} ay)\\
\calB(Lx, Ly) & \calA(x, RLy)
\Ar{1-1}{1-2}{"\nu", }
\Ar{1-2}{2-2}{"{\calA(x, t_y)}", "\rotatebox{90}{$\sim$}"'}
\Ar{1-1}{2-1}{"{(L, \ph)^{(2)}_{x,y}}"'}
\Ar{2-1}{2-2}{"\sim","\om_{x, Ly}"'}
\end{tikzcd},
\end{equation}
where $\nu$ is the canonical morphism.
Note that $\nu$ becomes an isomorphism if $x \in \cpt(\calA)$.
\item[{\rm (PA3)}]
$R$ preserves small coproducts.
\end{enumerate}
If in the above, the condition (PA3) is {\em not assumed}, then
$(L, \ph)$ is said {\em to induce a $G$-precovering}
$(L, \ph) \colon \cpt(\calA) \to \calB$
{\em by the adjunction} $\om$.
%
%
Note that by the condition (PA2), the induced functor $(L, \ph) \colon \cpt(\calA) \to \calB$ is actually a $G$-precovering, and
the condition (PA3) guarantees that 
the functor $(L, \ph) \colon \cpt(\calA) \to \cpt(\calB)$
is induced by Lemma \ref{lem:left-adj-cpt}, and is a $G$-precovering.
Note the difference of two types.
For the former, its codomain is $\cpt(\calB)$ and (PA3) is assumed,
while for the latter, its codomain is $\calB$ and (PA3) is not assumed.
\end{definition}

By noting that ${\rm cpt}(\Mod \CC)=\mmod \CC$ and ${\rm cpt}(\Mod \CC/G)=\mmod \CC/G$
for a small $G$-category $\calC$,
Definition \ref{dfn:precov-adj} rephrases Proposition \ref{Vitalcommutativediagram}
as follows.

\begin{proposition}
\label{prp:pushdown-precov}
Let $\calC$ be a skeletally small $G$-category, $(P, \ph) \colon \calC \to \calC/G$ the canonical functor.
Then the pushdown functor
$(P\down, \ph\down) \colon \Mod \calC \to \Mod \calC/G$ induces a $G$-precovering
$(P\down, \ph\down) \colon \mmod\calC \to \mmod \calC/G$ by the adjunction $\th \colon P\down \adj P\up$. \qed
\end{proposition}

Let $\CC$ be a $G$-category, and let $\mathcal{K}$ be a full subcategory of $\CC$ that is $G$-stable, i.e., $A_a(\mathcal{K}) \subseteq \mathcal{K}$ for every $a \in G$. Assume that $\CC$ is idempotent complete.

Under these conditions, we can make the following observations:

\begin{itemize}
\item [$(1)$] The subcategories ${\rm Add}\mbox{-}\CK$ and ${\rm add}\mbox{-}\CK$ are also $G$-stable;
    \item [$(2)$] If $\CK \subseteq \cpt(\CC)$, then $\cpt({\rm Add}\mbox{-}\CK)={\rm add}\mbox{-}\CK.$
\end{itemize}
In particular, for a skeletally small category $\CC$, we have
\[
\cpt(\Mod \CC)=\cpt({\rm Add}\mbox{-}(\mmod \CC))=\mmod \CC.
\]

\begin{example}
\label{ex-G-stable}
Keep the notations used in subsection \ref{Pushdownfunctor}.
%
Let $\CK$ be $G$-stable subcategory of $\mmod \CC$,
and $\CK'$ a subcategory of $\mmod(\CC/G)$.
Assume $P\down$ sends $\CK$ to $\CK'$ and $P\up$ sends $\CK'$ to ${\rm Add}\mbox{-}\CK$.
Then, we have  $P\down ({\rm Add}\mbox{-}\CK) \subseteq {\rm Add}\mbox{-}\CK'$
and $P\up({\rm Add}\mbox{-}\CK')\subseteq {\rm Add}\mbox{-}\CK$
because both $P\down$ and $P\up$ preserve small coproducts.
Hence, according to the observation made earlier (before the example),
we can restrict the $G$-invariant functor 
$(P\down, \ph\down)\colon \Mod \CC \to \Mod (\CC/G)$
to obtain the $G$-invariant functor 
$(P\down|, \ph\down|)\colon {\rm  Add}\mbox{-}\CK\rt {\rm Add}\mbox{-}\CK'$,
which induces a $G$-precovering
$(P\down|, \ph\down |)\colon {\rm  add}\mbox{-}\CK\rt {\rm add}\mbox{-}\CK'$
by the adjuction $\theta| \colon (P\down|) \adj (P\up|)$.
\end{example}

\begin{remark}
The definition above is extended to the case that $\calA, \calB \in \kCat_s$ in a canonical way, which is used in Theorem \ref{thm:general-case}.
To see the difference, note the required conditions on $L$, $R$ and $\om$ in this theorem.
\end{remark}

To deal with localization cases, we introduce the following terminology.

\begin{definition}
\label{dfn:epic-among}
Let $\pi \colon \calA \to \calB$ be a $k$-functor of $k$-categories
that is the identity on the objects, that is, $\obj(\calA) = \obj(\calB)$, and $\pi(x) = x$ for all $x \in \calA$.


$\pi$ is said to be {\em epic among representable right $\calB$-modules}
if for any $y \in \cpt(\calA)$, $z \in \calA$, and any morphims
$u, v \colon \calB(\blank, y) \to \calB(\blank, z)$ of right $\calB$-modules,
the equality $u\pi = v\pi$ on $\cpt(\calA)(\blank, y)$, 
(namely, $u \circ \pi_{(\blank,y)}
= v \circ \pi_{(\blank,y)}$,
see the diagram below)
implies $u = v$:
\[
\begin{tikzcd}
\cpt(\calA)(\blank,y) & \calB(\blank, y) & \calB(\blank, z)
\Ar{1-1}{1-2}{"\pi_{(\blank,y)}"}
\Ar{1-2}{1-3}{"u", yshift=5pt}
\Ar{1-2}{1-3}{"v"', yshift=-5pt}
\end{tikzcd}.
\]
By Yoneda's lemma, this condition is equivalent to the following:
For any $y \in \cpt(\calA)$, $z \in \calA$, and any morphims
$u, v \in \calB(y, z)$,
the equality $u \circ \pi(g) = v \circ \pi(g)$ for all $x \in \cpt(\calA)$
and $g \in \calA(x,y)$ implies $u = v$.
\end{definition}




\begin{example}
\label{exm:ess-surj-repbl}
Let $\pi \colon \calA \to \calB$ be a $k$-functor of $k$-categories.
In the following cases $\pi$ is
epic among representable right modules:

\begin{enumerate}
\item
The case where $\pi$ is a full functor.
\item
The case where $S$ is a bicalculable multiplicative system of $\calA$,
$\calB:= \calA[S\inv]$, and $\pi$ is the localization functor.
\end{enumerate}
In the case (1), the assertion is obvious.
In the case (2), note that $\pi$ is the identity on the objects.
To show the assertion, 
let $x \in \calA$, $y \in \calB$.
%
and $u, v \colon \calB(\blank, x) \to \calB(\blank, y)$ morphisms of right $\calB$-modules,
and assume that the equality $u\pi = v\pi$ holds on $\calA(\blank, x)$.
Take any $z \in \calA$ and $g \in \calB(z,x)$.
Then $g$ has the form $g = \pi(f)\pi(s)\inv$ for some
$f \in \calA(w,x)$, $s \in S(w,z)$ and $w \in \calA$.
By assumption, we have $u(\pi(f)) = v(\pi(f))$.
Since both $u$ and $v$ are morphisms of right $\calB$-modules, we have
$u(g) = u(\pi(f)\pi(s)\inv) = u(\pi(f))\pi(s)\inv =  v(\pi(f))\pi(s)\inv = v(\pi(f)\pi(s)\inv) = v(g)$.
Thus $u = v$.
\end{example}

\begin{lemma}
\label{lem:compact}
Let $\pi \colon \calA \to \calB$ be a $k$-functor of $k$-categories.
Let $I$ be a small set, and $(y_i)_{i\in I}$ any family
of objects in $\calA$, and
assume that $(y_i)_{i\in I}$ has a coproduct
$(y, (s_i \colon y_i \to y)_{i\in I})$ in $\calA$.
Then the following diagram is commutative for all $x \in  \calA$:
\begin{equation}
\label{eq:compact-preserved}
\begin{tikzcd}
\Ds_{i\in I}\calA(x, y_i) & \calA(x, y)\\
\Ds_{i\in I}\calB(\pi(x), \pi(y_i)) & \calB(\pi(x), \pi(y))
\Ar{1-1}{1-2}{"\nu_{x}"}
\Ar{2-1}{2-2}{"\nu'_{x}"'}
\Ar{1-1}{2-1}{"{\Ds_{i\in I}\pi_{(x,y_i)}}"'}
\Ar{1-2}{2-2}{"{\pi_{(x,y)}}"}
\end{tikzcd},
\end{equation}
where $\nu_x$ and $\nu'_x$ are the canonical maps defined 
as in \eqref{eq:compact}.
Note here that $\nu'_x$ is defined even if the family
$(\pi(y), (\pi(s_i) \colon \pi(y_i) \to \pi(y))_{i\in I})$ is not a coproduct
of $(\pi(y_i))_{i\in I}$ because this family is a cocone.
\end{lemma}

\begin{proof}
Take any $(f_i)_{i\in I} \in \Ds_{i\in I}\calA(z, y_i)$.
Then the clockwise composite of morphisms send it to
$\pi(\sum_{i\in I} s_i \circ f_i)$, and the counter-clockwise composite of morphisms send it to
$\sum_{i\in I} \pi(s_i) \circ \pi(f_i)$,
and hence they coincide.
\end{proof}

\begin{lemma}
\label{lem:cpt-cpt}
Let $\pi \colon \calA \to \calB$ be a $k$-functor of $k$-categories, and assume that
\begin{enumerate}
\item
$\calA$ has small coproducts,
\item
$\pi$ preserves small coproducts,
\item
$\pi$ is surjective on objects, and
\item
$\pi$ is full.
\end{enumerate}
Then $\pi(\cpt(\calA)) \subseteq \cpt(\calB)$.
\end{lemma}

\begin{proof}
Let $x \in \cpt(\calA)$.
To show that $\pi(x)$ is in $\cpt(\calB)$,
take any small set $I$ and any family $(y'_i)_{i\in I}$ of objects of $\calB$.
By (3), there exist objects $y_i$ of $\calA$
such that $y'_i = \pi(y_i)$.
By (1), there exists a coproduct
$(y, (s_i \colon y_i \to y)_{i\in I})$ of $(y_i)_{i\in I}$ in $\calA$.
By (2),
$(\pi(y), (\pi(s_i) \colon y'_i \to \pi(y))_{i\in I})$
turns out to be a coproduct of $(y'_i)_{i\in I}$ in $\calB$.
Then we have the commutative diagram \eqref{eq:compact-preserved}.
It is enough to show that $\nu'_x$ is surjective by Remark \ref{rmk:compact}.
But this follows from the facts that
$\pi_{(x,y)}$ is surjective by (4)
and that $\nu_x$ is an isomorphism because $x \in \cpt(\calA)$.
As a consequence, $\pi(x) \in \cpt(\calB)$.
\end{proof}

\begin{lemma}
\label{lem:cpt-cpt-mod}
Let $\pi \colon \calA \to \calB$ be a full $k$-functor of $k$-categories,
which is identity on objects.
If $\calA = \Mod \calC$ for some skeletally small $k$-category $\calC$,
then $\pi$ preserves small coproducts.
Hence by Lemma \ref{lem:cpt-cpt}, we have $\pi(\cpt(\calA)) \subseteq \cpt(\calB)$.
\end{lemma}

\begin{proof}
Assume that $\calA = \Mod \calC$ for some skeletally small $k$-category $\calC$.
Then it is well know that $\calA$ has small coproducts that is given as direct sums.
To show that $\pi$ preserves small coproducts, let $I$ be a small set,
and let $M_i \in \calA$ for all $i \in I$.
Set $M:= \Ds_{i\in I} M_i$, and let $q_i \colon M_i \to M$ be the canonical injection,
and $p_i \colon M \to M_i$ the composite of the inclusion
$\si \colon M \to \prod_{j\in I} M_j$ followed by
the canonical projection $\prod_{j\in I}M_j \to M_i$ for all $i\in I$.
Then we have
\begin{align}
&p_i q_j = \de_{i,j}\idty_{M_i} \ \text{ for all }i,j \in I; \text{ and}\\
\label{eq:sum-coprod}
&\sum_{i \in I}q_i p_i = \idty_M,
\end{align}
where the left hand side of the second equality is summable in the sense that
for each $x \in \calC$, the set
$\{i \in I \mid q_{i,x}p_{i,x} \ne 0 \colon M(x) \to M(x)\}$ is finite.
To complete the proof,
it is enough to show that $(M, (q'_i)_{i\in I})$
is a coproduct of $(M_i)_{i\in I}$ in $\calB$
because $(M, (q_i)_{i\in I})$ is a coproduct of $(M_i)_{i\in I}$ in $\calA$,
where we set $q'_i:= \pi(q_i)$ for all $i\in I$ for short.
To show this, take any $N \in \calB$ and
let $f'_i$ be in $\calB(M_i, N)$ for all $i \in I$.
Then since $\pi$ is full, for each $i\in I$, there exists $f_i \in \calA(M_i, M)$
such that $\pi(f_i) = f'_i$.
By the universality of the coproduct $(M, (q_i)_{i\in I})$,
there exists a unique $f \in \calA(M, N)$ such that
$f_i = f q_i$ for all $i\in I$.
By setting $f':= \pi(f) \in \calB(M,N)$, we have
$f'_i = f' q'_i$ for all $i\in I$.
It remains to show the uniqueness of this $f'$.
To show this,
assume that $g' \in \calB(M,N)$ satisfies $f'_i = g'q'_i$ for all $i\in I$.
We have only to show that $f' = g'$.
By assumption, there exists some $g \in \calA(M,N)$ such that $g' = \pi(g)$.
Set $u:= f - g \in \calA(M,N)$.
Then $\pi(uq_i) = \pi(fq_i -gq_i) = f'q'_i - g'q'_i = f'_i -f'_i = 0$
for all $i \in I$.
Hence by \eqref{eq:sum-coprod},
we have\[
\pi(u) = \pi(\sum_{i\in I}uq_ip_i) = \sum_{i\in I}\pi(uq_i)\pi(p_i) = 0.
\]
Therefore, $0 = \pi(f-g) = f' - g'$, and thus we have $f' = g'$.
\end{proof}

\subsection{Main theorems}

In the following theorems, we give a general tool to induce
$G$-precoverings in terms of 2-categories.

\begin{theorem}
\label{thm:factor-case}
Consider the following diagram of $2$-categories, $2$-functors,
and a strict $2$-natural transformation on the left, and its extension

by Lemma \ref{lem:ext-2-nat} on the right:
\[
\begin{tikzcd}
\sfC & \makebox[1em][l]{$\kCAT$,}
\Ar{1-1}{1-2}{"\si", ""'{name=a}, bend left}
\Ar{1-1}{1-2}{"U"', ""{name=b}, bend right}
\Ar{a}{b}{"\pi", Rightarrow}
\end{tikzcd}\qquad\qquad
\begin{tikzcd}
\hat{\sfC} & \makebox[1em][l]{$\GCAT$,}
\Ar{1-1}{1-2}{"\hat{\si}", ""'{name=a}, bend left}
\Ar{1-1}{1-2}{"\hat{U}"', ""{name=b}, bend right}
\Ar{a}{b}{"\hat{\pi}", Rightarrow}
\end{tikzcd}
\]
where $\sfC$ is a $2$-subcategory of $\kCAT$ and $\si$ is the inclusion $2$-functor.
Assume the following conditions for all $\calC \in \sfC$:
\begin{enumerate}
\item[{\rm (0)}]
$\calC$ has small coproducts, and
$\pi_{\calC} \colon \calC \to U(\calC)$ preserves small coproducts.
%
\item[{\rm (i)}]
$\pi_{\calC}$ is the identity on the objects.
\item[{\rm (ii)}]
$\pi_{\calC}$ is epic among
representable right $U(\calC)$-modules.
\end{enumerate}
Now let $(\calA, A)$ be a $G$-category in $\hat{\sfC}$,
where $\calA$ has small coproducts,
$\calB$ a $k$-category in $\sfC$,
and $(L, \ph) \colon (\calA, A) \to \calB$ a $G$-invariant functor in $\hat{\sfC}$
(cf. Remark 1.6),
where $L$ has a right adjoint $R \colon \calB \to \calA$ in $\sfC$
with an adjunction isomorphism $\om \colon L \adj R$ in $\sfC$.
Then the following hold.
\begin{enumerate}
\item
Assume that $(L, \ph) \colon \calA \to \calB$ induces a $G$-precovering $\cpt(\calA) \to \calB$ by the adjunction $\om$.
Then
\begin{enumerate}
\item
the $G$-invariant functor
$U(L, \ph) \colon U(\calA) \to U(\calB)$ induces a $G$-precovering
\[
U(L,\ph) \colon \cpt(U(\calA)) \to U(\calB)
\]
by the adjunction $U(\om) \colon U(L) \adj U(R)$, and
\item

if further $\cpt(\calA) = \pi_\calA(\cpt(\calA)) \subseteq \cpt(U(\calA))$
(this is the case, for example, 
if $\pi_{\calA}$ is 
full by Lemma \ref{lem:cpt-cpt-mod}),
then
$U(L, \ph)$ restricts to a $G$-precovering
\[
U(L,\ph) \colon U(\cpt(\calA)) \to U(\calB).
\]
\end{enumerate}

\item
Assume that $(L, \ph) \colon \calA \to \calB$ induces a $G$-precovering $\cpt(\calA) \to \cpt(\calB)$
by the adjunction $\om$.
Then
\begin{enumerate}
\item
the $G$-invariant functor
$U(L, \ph) \colon U(\calA) \to U(\calB)$ induces a $G$-precovering
\[
U(L,\ph) \colon \cpt(U(\calA)) \to \cpt(U(\calB))
\]
by the adjunction $U(\om) \colon U(L) \adj U(R)$, and
\item

if further $\cpt(\calA) \subseteq \cpt(U(\calA))$, then
$U(L, \ph)$ restricts to a $G$-precovering
\[
U(L,\ph) \colon U(\cpt(\calA)) \to U(\cpt(\calB)).
\]
\end{enumerate}
%
%
%
\end{enumerate}
\end{theorem}

\begin{proof}
First of all, we remark that for each $\calA \in \hat{\sfC}$,
$\hat{\pi}_\calA$ satisfies the conditions corresponding to both (i) and (ii)
by \eqref{eq:ext-nat-tr} as well.
By abuse of notation, we denote $\hat{U}$ and $\hat{\pi}$ just by $U$ and $\pi$,
respectively.

(1)(a)
Assume that $(L, \ph) \colon \calA \to \calB$ in $\sfC$ induces
a $G$-precovering $\cpt(\calA) \to \calB$ by the adjunction $\om \colon L \adj R$.
Then $(L, \ph)$ satisfies the conditions (PA1) and (PA2) in Definition \ref{dfn:precov-adj}.
Since $U$ is a 2-functor, we have an adjunction $U(\om) \colon U(L) \adj U(R)$ by Lemma \ref{lem:2-fun-pres-adj}.
By (PA1) for $(L, \ph)$, we have a natural isomorpism
\[
U(t) \colon \Ds_{a\in G}U(A_a) \To U(R)U(L).
\]
Thus $U(L, \ph)$ satisfies (PA1).

To show the statement (a),
it remains to show the condition (PA2) for $U(L, \ph)$,
namely the commutativity of the diagram
\begin{equation}\label{eq:new-precov-adj}
\begin{tikzcd}[column sep = 50pt]
\Ds_{a \in G} U(\calA)(x, ay) & U(\calA)(x, \Ds_{a \in G} ay)\\
U(\calB)(U(L)X, U(L)y) & U(\calA)(x, U(R)U(L)y)
\Ar{1-1}{1-2}{"\nu'"}
\Ar{1-2}{2-2}{"{U(\calA)(x, U(t_y))}"}
\Ar{1-1}{2-1}{"{(U(L, \ph))^{(2)}_{x,y}}"'}
\Ar{2-1}{2-2}{"U(\om)_{x, U(L)y}"'}
\end{tikzcd}
\end{equation}
for all $x,y \in U(\calA)$, where $\nu'$ is the canonical morphism.
Construct the following diagram by combining the diagrams \eqref{eq:precov-adj} and \eqref{eq:new-precov-adj} with $\pi$:
\begin{equation}
\label{eq:cubic-diagram}
\footnotesize
\begin{tikzcd}[column sep=40pt]
	{\Ds_{a\in G}\calA(x,ay)} & {\calA(x,\Ds_{a\in G} ay)}\\
	{\calB(Lx, Ly)} & {\calA(x,RLy)}\\
	\\
	& {\Ds_{a\in G}U(\calA)(x,ay)} & {U(\calA)(x,\Ds_{a\in G} ay)} \\
	& {U(\calB)(U(L)x,U(L)y)} & {U(\calA)(x,U(R)U(L)y)}
	\arrow["{(L, \ph)^{(2)}_{x,y}}"', from=1-1, to=2-1]
	\arrow["{\calA(x,t_y)}"', from=1-2, to=2-2]
	\arrow["\nu", from=1-1, to=1-2]
	\arrow["{\om_{x,Ly}}", from=2-1, to=2-2]
	\arrow["{(U(L),U(\ph))^{(2)}_{x,y}}", from=4-2, to=5-2]
	\arrow["{U(\calA)(x,U(t)_y})", from=4-3, to=5-3]
	\arrow["{U(\om)_{x,U(L)y}}"', from=5-2, to=5-3]
	\arrow["{\Ds_{a\in G}\pi_{\calA,(x,ay)}}"{description}, from=1-1, to=4-2, crossing over]
	\arrow["{\pi_{\calB, (Lx,Ly)}}"{description}, from=2-1, to=5-2]
	\arrow["{\pi_{\calA, (x, \Ds_{a\in G}ay)}}", from=1-2, to=4-3]
	\arrow["{\pi_{\calA, (x, RLy)}}"{description}, from=2-2, to=5-3, near start]
	\arrow["{\nu'}"', from=4-2, to=4-3, crossing over]
\end{tikzcd}.
\end{equation}
Then the back square commutes by the condition (PA2), 
the commutativity of \eqref{eq:precov-adj}.
We now verify the commutativity of the left, right, upper, and lower squares.

\begin{clm}
The left square of \eqref{eq:cubic-diagram} is commutative.
\end{clm}
Indeed, we first show the commutativity of the following diagaram:
\begin{equation}
\label{eq:left-sq}
\begin{tikzcd}[column sep=50pt]
	{\calA(x,ay)} & {U(\calA)(x,ay)} \\
	{\calB(Lx, Lay)} & {U(\calB)(U(L)x,U(L)ay)} \\
	{\calB(Lx, Ly)} & {U(\calB)(U(L)x,U(L)y)}
	\arrow["{\pi_{\calA, (x,ay)}}", from=1-1, to=1-2]
	\arrow["{\pi_{\calB,(Lx, Lay)}}", from=2-1, to=2-2]
	\arrow["{\pi_{\calB, (Lx, Ly)}}"', from=3-1, to=3-2]
	\arrow["{L_{(x,ay)}}"', from=1-1, to=2-1]
	\arrow["{U(L)_{(x,ay)}}", from=1-2, to=2-2]
	\arrow["{\calB(Lx,\phi\inv_{a,y})}"', from=2-1, to=3-1]
	\arrow["{U(\calB)(U(L)x,U(\phi)\inv_{a,y})}", from=2-2, to=3-2]
\end{tikzcd}.
\end{equation}
Consider the diagram
\[
\begin{tikzcd}[row sep=50pt, column sep=60pt]
\calA & \calB\\
U(\calA) & U(\calB)
\Ar{1-1}{1-2}{"LA_a", ""'{name=a}, bend left}
\Ar{1-1}{1-2}{"L"', ""{name=b}, bend right}
\Ar{a}{b}{"\ph_a\inv", Rightarrow}
\Ar{2-1}{2-2}{"U(L)U(A_a)", ""'{name=aa}, bend left}
\Ar{2-1}{2-2}{"U(L)"', ""{name=bb}, bend right}
\Ar{aa}{bb}{"U(\ph_a\inv)", Rightarrow}
\Ar{1-1}{2-1}{"\pi_\calA"'}
\Ar{1-2}{2-2}{"\pi_\calB"}
\end{tikzcd}.
\]

Then by the strict $1$-naturality of $\pi$, we have $\pi_\calB L = U(L) \pi_\calA$, 
which shows the commutativity of the upper square of \eqref{eq:left-sq}, and
by the strict $2$-naturality of $\pi$, we have $\pi_\calB \ph_a\inv = U(\ph_a\inv)\pi_\calA$,
which shows the commutativity of the lower square of \eqref{eq:left-sq}
because any $g \in \calB(Lx, Lay)$ is sent by the counter-clockwise composite to
$\pi_\calB(\ph\inv_{a,y}\circ g)$, by the clockwise composite to
$U(\ph)\inv_{a,y}\circ \pi_\calB(g)$, and they coincide by
\[
\pi_\calB(\ph\inv_{a,y}\circ g) = \pi_\calB(\ph\inv_{a,y})\circ \pi_\calB(g) = U(\ph\inv_{a})_{\pi_\calA(y)}\circ \pi_\calB(g)=U(\ph)\inv_{a,y}\circ \pi_\calB(g).
\]
Now to show the claim, take any $(f_a)_{a\in G} \in \Ds_{a\in G}\calA(x,ay)$.
Then the counter-clockwise composite sends it to
$\pi_\calB(\sum_{a\in G}\ph\inv_{a,y}L(f_a)) = \sum_{a\in G}\pi_\calB(\ph\inv_{a,y}L(f_a))$,
and the clockwise composite sends it to
$\sum_{a\in G}U(\ph\inv_a)_y U(L)(\pi_\calA(f_a))$, and they coincide by the commutativity
of \eqref{eq:left-sq}.
\begin{clm}
The right square of \eqref{eq:cubic-diagram} is commutative.
\end{clm}
indeed,
by applying the strict $2$-naturality of $\pi$ to the diagram
\[
\begin{tikzcd}[row sep=50pt, column sep=60pt]
\calA & \calA\\
U(\calA) & U(\calA)
\Ar{1-1}{1-2}{"\Ds_{a\in G}A_a", ""'{name=a}, bend left}
\Ar{1-1}{1-2}{"RL"', ""{name=b}, bend right}
\Ar{a}{b}{"t", Rightarrow}
\Ar{2-1}{2-2}{"U(\Ds_{a\in G}A_a)", ""'{name=aa}, bend left}
\Ar{2-1}{2-2}{"U(R)U(L)"', ""{name=bb}, bend right}
\Ar{aa}{bb}{"U(t)", Rightarrow}
\Ar{1-1}{2-1}{"\pi_\calA"'}
\Ar{1-2}{2-2}{"\pi_\calA"}
\end{tikzcd},
\]
we have $\pi_\calA \circ t = U(t) \circ \pi_\calA$.
By evaluating it at $y$, we have
\begin{equation}
\label{eq:pi-U}
\pi_\calA(t_y) = U(t)_{\pi_{\calA}(y)} = U(t)_y.
\end{equation}
To show the claim, take any $f \in \calA(x, \Ds_{a\in G}ay)$.
Then the counter-clockwise composite sends it to $\pi_\calA(t_y \circ f) = \pi_\calA(t_y) \circ \pi_\calA(f)$, while the clockwise composite sends it to
$U(t)_y \circ \pi_\calA(f)$, and they coincide by \eqref{eq:pi-U}.
\begin{clm}
The upper square of \eqref{eq:cubic-diagram} is commutative.
\end{clm}
This follows by Lemma \ref{lem:compact}.

%
\begin{clm}
The lower square of \eqref{eq:cubic-diagram} is commutative.
\end{clm}
Indeed, let $f \in \calB(Lx, Ly)$.
Then by the definition of $U(\om)$ 
in Lemma \ref{lem:2-fun-pres-adj},
the counter-clockwise composite sends $f$ to $U(R)(\pi_\calB(f)) \circ U(\et)_x$,
while the clockwise composite sends $f$ to $\pi_\calA(Rf \circ \et_x) = \pi_\calA(Rf) \circ \pi_\calA(\et_x)$.
Thus, it is enough to show that
\begin{equation}
\label{eq:lower-sq}
U(R)(\pi_\calB(f)) \circ U(\et)_x = \pi_\calA(Rf) \circ \pi_\calA(\et_x).
\end{equation}
By the commutative diagram
\begin{equation}
\label{eq:nat-pi-R}
\begin{tikzcd}
\calB & \calA\\
U(\calB) & U(\calA)
\Ar{1-1}{1-2}{"R"}
\Ar{2-1}{2-2}{"U(R)"'}
\Ar{1-1}{2-1}{"\pi_\calB"'}
\Ar{1-2}{2-2}{"\pi_\calA"}
\end{tikzcd},
\end{equation}
we have
\begin{equation}\label{eq:lower-eq-1-nat}
U(R)(\pi_\calB(f)) = \pi_\calA(Rf).
\end{equation}
By the strict $2$-naturality of $\pi$, we also have $\pi_\calA \circ \et = U(\et) \circ \pi_\calA$.
By evaluating it at $x \in \calA$, we have
\begin{equation}\label{eq:lower-eq-2-nat}
\pi_\calA(\et_x) = U(\et)_{\pi_\calA(x)} = U(\et)_x.
\end{equation}
\eqref{eq:lower-eq-1-nat} and \eqref{eq:lower-eq-2-nat} show the equality \eqref{eq:lower-sq}.

We now show the commutativity of \eqref{eq:new-precov-adj}.
By fixing $y \in \calA$, we set $u_x$ and $v_x$ to be the clockwise composite and
counter-clockwise composite of \eqref{eq:new-precov-adj}, respectively,
and show that $u_x = v_x$ for all $x \in \calA$.
By the claims above, we have
\[
u_x \circ \Ds_{a\in G}\pi_{\calA,(x,ay)}
= v_x \circ \Ds_{a\in G}\pi_{\calA,(x,ay)}.
\]
For each $a \in G$,
let $\si_a \colon U(\calA)(x, ay) \to \Ds_{a\in G}U(\calA)(x, ay)$ be the  
canonical injection,
and set $u_{x,a}:= u_x\si_a, v_{x.a}:= v_x\si_a$ so that $u_x = (u_{x,a})_{a\in G}$
and $v_x = (v_{x,a})_{a\in G}$.
Then the equality above gives us the following for all $a \in G$:
\begin{equation}
\label{eq:comm-with-pi}
u_{x,a} \circ \pi_{\calA,(x,ay)}
= v_{x,a} \circ \pi_{\calA,(x,ay)}.
\end{equation}
Let $a \in G$.
Here, regard $x$ as a variable, and consider the correspondences
\[
u_a, v_a \colon U(\calA)(\blank, ay) \to U(\calA)(\blank, U(R)U(L)y)
\]
given by $u_{x,a}$ and $v_{x,a}$, respectively.
Then it is enough to show that these are morphisms of right $U(\calA)$-modules
because if this is the case,
then since $\pi_\calA$ is epic among representable right $U(\calA)$-modules by assumption (ii),

we have $u_a = v_a$ for all $a \in G$,
and hence $u_x =v_x$ for all $x \in \calA$, as desired.

By definition, we see that
\begin{align}
\label{eq:u}
u_a(f) &= U(t)_y \circ s_a \circ f\\
\label{eq:v}
v_a(f) &= U(R)(\ph_{a,y}\inv) \circ  (U(R)U(L))(f) \circ U(\et)_x
\end{align}
for all $f \in U(\calA)(x, ay)$ and $x \in U(\calA)$,
where $s_a \colon ay \to \Ds_{a\in G} ay$ is the injection of coproduct.
By the form of $\eqref{eq:u}$, it is obvious that $u_a$ is a morphism of right $U(\calA)$-modules.
To show that $v_a$ is a morphism of right $U(\calA)$-modules, take any $g \in U(\calA)(z,x)$ with $z \in U(\calA)$.
We have only to show that $v_a(fg) = v_a(f)g$.
By the naturality of $U(\et)$, we have
\[
\begin{aligned}
v_a(f)g &= U(R)(\ph_{a,y}\inv) \circ  (U(R)U(L))(f) \circ U(\et)_x \circ g\\
&= U(R)(\ph_{a,y}\inv) \circ  (U(R)U(L))(f)\circ (U(R)U(L))(g) \circ U(\et)_z\\
&= U(R)(\ph_{a,y}\inv) \circ  (U(R)U(L))(fg) \circ U(\et)_z\\
&= v_a(fg).
\end{aligned}
\]
This finishes the proof of (1)(a).

(1)(b) Assume that $\cpt(\calA) \subseteq \cpt(U(\calA))$.
More precisely stated,
$\obj(\cpt(\calA)) \subseteq \obj(\cpt(U(\calA)))$.
Then since $\pi$ is the identity on objects, we have
\[
\obj(U(\cpt(\calA))) = \obj(\cpt(\calA))
\subseteq \obj(\cpt(U(\calA))).
\]
Thus roughly stated, $U(\cpt(\calA)) \subseteq \cpt(U(\calA))$.
Hence the assertion follows.

(2)(a) This follows from (1) by the following claim.
\begin{clm}
If $R$ preserves small coproduts, then so does $U(R)$.
\end{clm}
To show this, look at the commutative diagram \eqref{eq:nat-pi-R},
and let $I$ be a small set and $(x_i)_{i \in I}$ a family of objects of $U(\calB)$.
Then since $\calB$ has small coproducts, there exists a coproduct
 $(y, (s_i \colon x_i \to y)_{i \in I})$ of $(x_i)_{i\in I}$ in $\calB$,
and then
$(y, (\pi_\calB(s_i) \colon x_i \to y)_{i \in I})$ turns out to be 
a coproduct of $(x_i)_{i\in I}$ in $U(\calB)$ because
$\pi_\calB$ preserves coproducts.
By the universality of coproducts, any other coproduct of
$(x_i)_{i\in I}$ in $U(\calB)$ is isomorphic to this coproduct.
Then since both $R$ and $\pi_\calA$ preserve coproducts, we see that
$U(R)(y, (\pi_\calB(s_i))_{i\in I}) = (U(R) \circ \pi_\calB)(y, (s_i)_{i\in I})
= (\pi_\calA \circ R)(y, (s_i)_{i\in I})$
is a coproduct of $(U(R)(x_i))_{i\in I}$ in $U(\calA)$, which proves the claim.


(2)(b) 
By (PA3) and Lemma \ref{lem:left-adj-cpt}, $(L, \ph)$ restricts to
$(L', \ph') \colon \cpt(\calA) \to \cpt(\calB)$.
Then in the diagram below, the $2$-functor $U$ sends
the commutative outer square to
the commutative inner square, and other four trapezoids are commutative
by the naturality of $\pi$:
\[
\begin{tikzcd}
\cpt(\calA) &&& \cpt(\calB)\\
&U(\cpt(\calA)) & U(\cpt(\calB))\\
&U(\calA) & U(\calB)\\
\calA &&& \calB
\Ar{1-1}{1-4}{"{(L', \ph')}"}
\Ar{4-1}{4-4}{"{(L, \ph)}"'}
\Ar{1-1}{4-1}{"s"', hookrightarrow}
\Ar{1-4}{4-4}{"s'", hookrightarrow}
\Ar{2-2}{2-3}{"{U(L', \ph')}"}
\Ar{3-2}{3-3}{"{U(L, \ph)}"'}
\Ar{2-2}{3-2}{"U(s)"'}
\Ar{2-3}{3-3}{"U(s')"}
\Ar{1-1}{2-2}{"\pi_{\cpt(\calA)}"'}
\Ar{1-4}{2-3}{"\pi_{\cpt(\calB)}"}
\Ar{4-1}{3-2}{"\pi_{\calA}"'}
\Ar{4-4}{3-3}{"\pi_{\calB}"}
\end{tikzcd},
\]

where $s, s'$ are the inclusion functors.
Here, by assumption (i), both $U(s)$ and $U(s')$ are the inclusions on
the objects.
Thus $U(L,\ph)$ ristricts to $U(L',\ph')$ that sends objects of $U(\cpt(\calA))$ into those of $U(\cpt(\calB))$.
Then the assertion follows from (1)(b).
\end{proof}

%

To generalize the theorem above, we introduce the following notation.

\begin{notation}\label{Notat-com-general-thm}
Let $(\calA, S)$ (resp.\ $(\calA, A, S)$) be an object of $\kCAT_s$ (resp.\  $\GCAT_s$).
Then we set
\[
\cpt(\calA, S):= (\cpt(\calA), \cpt(S))\ \text{and}\ 
\cpt(\calA, A, S):= (\cpt(\calA), A|_{\cpt(\calA)}, \cpt(S)),
\]
where $\cpt(S):= S|_{\cpt(\calA)\times \cpt(\calA)} = \bigsqcup_{x,y \in \cpt(\calA)}S(x,y)$.
Note here that $A_a(\cpt(\calA)) \subseteq \cpt(\calA)$ for all $a \in G$
by Remark \ref{rmk:autoeq-cpt}, which enables us to define the restriction
$A|_{\cpt(\calA)} \colon G \to \Aut(\cpt(\calA)), a \mapsto A_a|_{\cpt(\calA)}$.
\end{notation}

Theorem \ref{thm:factor-case} is generalized to a $2$-subcategory of $\kCAT_s$
as follows.

\begin{theorem}
\label{thm:general-case}
Consider the following diagram of $2$-categories, $2$-functors,
and a strict $2$-natural transformation on the left, and its extension
by Lemma \ref{lem:ext-2-nat} on the right:
\[
\begin{tikzcd}
\sfC & \makebox[1em][l]{$\kCAT_s$,}
\Ar{1-1}{1-2}{"\si", ""'{name=a}, bend left}
\Ar{1-1}{1-2}{"U"', ""{name=b}, bend right}
\Ar{a}{b}{"\pi", Rightarrow}
\end{tikzcd}\qquad\qquad
\begin{tikzcd}
\hat{\sfC} & \makebox[1em][l]{$\GCAT_s$,}
\Ar{1-1}{1-2}{"\hat{\si}", ""'{name=a}, bend left}
\Ar{1-1}{1-2}{"\hat{U}"', ""{name=b}, bend right}
\Ar{a}{b}{"\hat{\pi}", Rightarrow}
\end{tikzcd}
\]
where $\sfC$ is a $2$-subcategory of $\kCAT_s$ and $\si$ is the forgetful $2$-functor.
Assume the following conditions for all $\calA \in \sfC$: 
\begin{enumerate}
\item[{\rm (0)}]
$\calA$ has small coproducts, and
$\pi_{\calA} \colon \calA \to U(\calA)$ preserves small coproducts.
%
\item[{\rm (i)}]
$\pi_{\calA}$ is the identity on the objects,
that is, $\obj(\calA) = \obj(U(\calA))$, and $\pi_\calA(x) = x$ for all $x \in \calA$.
\item[{\rm (ii)}]
$\pi_{\calA}$ is epic among representable right $U(\calA)$-modules.
\end{enumerate}
Let $(\calA, A, S) \in \hat{\sfC}$ 
and $(\calB, T) \in \sfC$, and assume that there exists 
a $G$-invariant functor $(L, \ph) \colon (\calA, A, S) \to (\calB,T)$
in\footnote{%
Here we regard $(\calB,T) \in \hat{\sfC}$ by the trivial $G$-action on $\calB$.} $\hat{\sfC}$ such that
$L$ has a right adjoint $R \colon (\calB, T) \to (\calA,S)$ in $\sfC$
with an adjunction isomorphism $\om \colon L \adj R$ in $\sfC$.
Then the following hold.
\begin{enumerate}
\item
Assume that $(L, \ph) \colon (\calA,A, S) \to (\calB,T)$ induces a $G$-precovering
$\cpt(\calA,A,S) $ $ \to (\calB,T)$
by the adjunction $\om$. Then
\begin{enumerate}
\item
the $G$-invariant functor
$U(L, \ph) \colon U(\calA,A,S) \to U(\calB,T)$ induces a $G$-precovering
\[
U(L,\ph) \colon \cpt(U(\calA,A,S)) \to U(\calB,T)
\]
by the adjunction $U(\om) \colon U(L) \adj U(R)$, and
\item
if $\cpt(\calA) \subseteq \cpt(U(\calA,S))$, then
$U(L, \ph)$ restricts to a $G$-precovering
\[
\begin{aligned}
U(L,\ph) \colon U(\cpt(\calA,A,S)) \to U(\calB,T).
\end{aligned}
\]
\end{enumerate}

\item
Assume that $(L, \ph) \colon (\calA, A, S) \to (\calB,T)$ induces a $G$-precovering $\cpt(\calA,A,S) $ $ \to \cpt(\calB,T)$
by the adjunction $\om$.  Then
\begin{enumerate}
\item
the $G$-invariant functor
$U(L, \ph) \colon U(\calA,A,S) \to U(\calB,T)$ induces a $G$-precovering
\[
U(L,\ph) \colon \cpt(U(\calA,A,S)) \to \cpt(U(\calB,T))
\]
by the adjunction $U(\om) \colon U(L) \adj U(R)$, and
\item
if $\cpt(\calA) \subseteq \cpt(U(\calA,S))$, then
$U(L, \ph)$ restricts to a $G$-precovering
\[
U(L,\ph) \colon U(\cpt(\calA,A,S)) \to U(\cpt(\calB,T)).
\]
\end{enumerate}
\end{enumerate}
%
\end{theorem}

\begin{proof}
Just by taking care of well-definedness of morphisms,
the same proof as Theorem \ref{thm:factor-case} works.
\end{proof}

\subsection{Applications}

\begin{example}
\label{exm:stable-cat}
We can apply Theorem \ref{thm:factor-case} to 
the diagram \eqref{eq:stable-2-nat} given in Example \ref{exm:stable-2-nat},
especially in the case that $\calA = \Mod\calC$ for a skeletally small $G$-category $\calC$.
Indeed, the assumptions (i) and (ii) are obviously satisfied.
By Lemma \ref{lem:cpt-cpt-mod}, the assumption (0) is
satisfied, and we also have
$\cpt(\calA) \subseteq \cpt(U(\calA))$.
Here, note that $\cpt(\calA) = \mmod\calC$,
$U(\calA) = \uMod\calC, U(\cpt(\calA)) = \umod \calC$.
\end{example}

We obtain the following by Theorem \ref{thm:factor-case} (2)(b).

\begin{proposition}\label{Prop 2.6}
Let $\calC$ be a skeletally small $G$-category. Then the following hold.
\begin{enumerate}
\item
The functor $\udl{(P\down, \ph\down)} \colon \umod\CC \to  \umod\CC/G$
induced by $(P\down, \ph\down)$ is a $G$-precovering.
\item
If $\CC$ is a locally support-finite locally bounded category,
then  $\udl{(P\down, \ph\down)} \colon \umod\CC \to  \umod\CC/G$ is a $G$-covering.
\end{enumerate}
\end{proposition}

\begin{proof}
(1)
As stated in Proposition \ref{prp:pushdown-precov},
$(P\down, \ph\down) \colon \Mod \calC \to \Mod \calC/G$ induces a $G$-precovering
$(P\down, \ph\down) \colon \mmod\calC \to \mod \calC/G$
by the adjunction $\th \colon P\down \adj P\up$.
As in Example \ref{exm:stable-cat}, we can apply
Theorem \ref{thm:factor-case}(2)(b) to $(P\down, \ph\down)$ above
to have the assertion (1).
(Recall that the right adjoint $P\up$ to $P\down$ has a right adjoint,
and hence Lemma \ref{lem:left-adj-cpt}
also applies.)

(2) By Proposition \ref{loc-supp-fin}, 
$(P\down, \ph\down) \colon \mmod\calC \to \mod \calC/G$ is dense.
Hence so is
$\udl{(P\down, \ph\down)} \colon \umod\CC \to  \umod\CC/G$.
\end{proof}

In the factor $2$-functor case and 
the quotient $2$-functor case, 
we obtain the following two propositions by
Theorem \ref{thm:general-case} (2)(b).
Before each proposition, we prepare a lemma to apply the theorem in each setting.

\begin{remark}
Let $\calA$ be a $k$-category, and $S$ a class of morphisms of $\calA$
having the form $S = \ang{\calD}$ for some class of objects of $\calA$.
Then 

(1) $S$ is $G$-stable if $\calD$ is $G$-stable.

Indeed, if $f \colon X \to Y$ is in $S$, then $f$ factors through
some $D \in \calD$.
Let $a \in G$.
Then $af$ factors through $aD$ that is in $\calD$ because
$\calD$ is $G$-stable.  Hence $af$ is in $S$, as required.

Now, let $D'$ be the isomorphism closure of $D$,
i.e., $D':= \{X \in \calA \mid X \cong D\ \text{for some}\ D \in \calD\}$.
Then clearly we have

(2) $\ang{D} = \ang{D'}$.
\end{remark}



\begin{lemma}
\label{lem:factor-case}
Consider the factor $2$-functor setting:
\[
\begin{tikzcd}
\makebox[1em][r]{$\GCAT_i$} & \makebox[1em][l]{$\GCAT$}
\Ar{1-1}{1-2}{"\si_i", ""'{name=a}, bend left}
\Ar{1-1}{1-2}{"U_i"', ""{name=b}, bend right}
\Ar{a}{b}{"\pi_i", Rightarrow}
\end{tikzcd}\qquad.
\]

Let $(\calA,S') \in \GCAT_i$, and assume that $\calA$ is given in the form
$\calA = \Mod\calC$ for some skeletally small $G$-category $\calC$
with the canonical $G$-covering $(P, \ph) \colon \calC \to \calC/G$, 
and that $S'$ is an ideal of the form $S' = \ang{\calD}$
for some $G$-stable class $\calD$ of objects of $\calA$
closed under small coproducts.
Then by the remark above $S'$ is $G$-stable (see Definition \ref{dfn:G-stable-2cats}(2)).
Set $S:= \ang{\cpt(\calD)}$ be the ideal of $\mmod \calC$.
Then $S$ is $G$-stable because so is $\cpt{\calD}$ by Remark \ref{rmk:autoeq-cpt}.
Finally, assume that each object in $\calD$ is decomposed
into a small coproduct of compact objects.
%
%
In this case, the following statements hold.
\begin{enumerate}
\item 
We have $\cpt(S') = S$.
Hence
\[
U_i(\cpt(\calA, S')) = U_i(\mmod\calC, S) = (\mmod \calC)/S.
\]
\item
Set
$S'/G:= \ang{P\down(\calD)}$ and 
$S/G:= \ang{P\down(\cpt(\calD))}$
to be the ideal of $\Mod(\calC/G)$ and $\mmod(\calC/G)$), respectively.
Then we have
\[
\cpt(S'/G) = S/G.
\]
Hence
\[
U_i(\cpt(\Mod (\calC/G), S'/G)) = U_i(\mmod (\calC/G), S/G) = (\mmod (\calC/G))/(S/G).
\]
\item
We have $P\down(S') \subseteq S'/G$ and $P\up(S'/G) \subseteq S'$.
Hence we see that $(P\down, \ph\down)$ induces a $G$-invariant functors
$(\Mod\calC, S') \to (\Mod \calC/G, S'/G)$ and
$(\mmod\calC, S) \to (\mmod \calC/G, S/G)$ and that
$P\up$ induces a right adjoint
$(\Mod\calC/G, S'/G) \to (\Mod\calC, S')$
to the former.
\end{enumerate}
\end{lemma}

\begin{proof}

(1)
Since $\cpt(\calD) \subseteq \calD$,
we have $S \subseteq S'$.
Let $f \colon X \to Y$ in $S'$, and assume that both $X$ and $Y$
are in $\mmod \calC$.
Then $f$ factors through an object $D \in \calD$,
which is isomorphic to the dirct sum $\Ds_{i \in I} D_i$
for some $D_i \in \cpt(\calD)$ and some small set $I$.
Since $X$ is compact, there exists some finite subset $J$ of $I$
such that $f$ factors through $\Ds_{i \in J} D_i$, which is in $\cpt(\calD)$.
Hence $f \in S$.

(2)
Since $P\down$ sends finitely generated $\calC$-modules
to finitely generated $\calC/G$-modules,
each object in $P\down(\calD)$ is decomposed into a small coproduct
of compact objects.
Then the same argument as in (1) works in this case to
show the assertion.

(3)
By construction, it is obvious that $P\down(S') \subseteq S'/G$.
Let $D \in \calD$.
Then $P\up P\down (D) \cong \Ds_{a \in G} {}^aD$.
Since $\calD$ is $G$-stable,
${}^aD$ are in $\calD$ for all $a \in G$.
Moreover, since $\calD$ is closed under small coproducts,
we have $\Ds_{a \in G} {}^aD \in \calD$.
Thus $P\up P\down (D) \in \calD$.
Now let $f \in S'/G$.
Then $f$ factors through an object of the form $P\down(D)$
for some $D \in \calD$.
Hence $P\up(f)$ factors through $P\up P\down(D)$, which is
in $\calD$ by the above.
Thus, $P\up(f) \in S'$, and we have 
$P\up(S'/G) \subseteq S'$.
\end{proof}

\begin{proposition}
\label{Prp:factor-covering}
Let $\calC$ be a skeletally small $G$-category, and 
$\calD$ a $G$-stable class of objects in $\Mod \calC$
closed under small coproducts with the property that
each object in $\calD$ is a small coproduct of finitely generated objects.
Denote by $S/G$ the ideal $\ang{P\down(\cpt(\calD))}$ of $\mmod (\calC/G)$. 
Then the following hold.
\begin{enumerate}
\item
The functor
\[
U(P\down, \ph\down) \colon (\mmod\CC)/S \to  (\mmod\CC/G)/(S/G)
\]
induced by the quotient $2$-functor $U = U_i$ in {\rm Lemma \ref{lem:factor-case}}
is a $G$-precovering.
\item
If $\CC$ is a locally support-finite locally bounded category,
then
\[
U(P\down, \ph\down) \colon (\mmod\CC)/S \to  (\mmod\CC/G)/(S/G)
\]
is a $G$-covering.
\end{enumerate}
\end{proposition}

\begin{proof}
By Lemma \ref{lem:factor-case}, we can apply
Theorem \ref{thm:factor-case} (2)(a),
and by Lemma \ref{lem:cpt-cpt-mod}, we have
$\cpt(\calA) \subseteq \cpt(U_i(\calA))$.
Hence we can apply Theorem \ref{thm:factor-case} (2)(b)
to have the statement (1) above.
The statement (2) follows from (1).
\end{proof}

We next consider the quotient case.

\begin{lemma}
\label{lem:quotient-case}
Consider the quotient $2$-functor setting:
\[
\begin{tikzcd}
\makebox[1em][r]{$\GCAT_m$} & \makebox[1em][l]{$\GCAT$}
\Ar{1-1}{1-2}{"\si", ""'{name=a}, bend left}
\Ar{1-1}{1-2}{"U_m"', ""{name=b}, bend right}
\Ar{a}{b}{"\pi", Rightarrow}
\end{tikzcd}\qquad.
\]

Let $(\calA,S') = (\calA, A, S') \in \kCAT_m$.
Then $\pi_{(\calA,S')} \colon \calA \to U_m(\calA, S')$
preserves small coproducts because it has right adjoint.
Hence by Lemma \ref{lem:cpt-cpt} , we have
$\cpt(\calA) \subseteq \cpt(U_m(\calA, S'))$ if $\calA$ has small coproducts.

Assume now that $\calA$ is given in the form
$\calA = \Mod\calC$ for some skeletally small $G$-category $\calC$
such that both $\mmod\calC$ and $\mmod(\calC/G)$ are abelian,
with the canonical $G$-covering $(P, \ph) \colon \calC \to \calC/G$.
Then $\calA$ has small coproducts as required above.
Let $\calS'$ be a $G$-stable localizing subcategory of $\Mod \calC$,
and set $\calS'/G:= \Loc(\{P\down(X)\mid X \in \calS\})$ to be the smallest localizing subcategory of $\Mod \calC/G$
containing $P\down(X)$ for all $X \in \calS$
(see Remark \ref{rmk:dense-localizing} for remarks on localizing
subcategories).
Then the canonical functors $Q \colon \Mod\calC \to (\Mod\calC)/\calS$
and $Q/G \colon \Mod(\calC/G) \to (\Mod(\calC/G))/(\calS/G)$
turn out to be localizations, i.e.,
they are exact and have fully faithful right adjoints $J$ and $J/G$, respectively.
Thus we have the following diagram with four adjoint pairs:
\begin{equation}
\label{eq:4adjoints}
\begin{tikzcd}[row sep=50pt, column sep=70pt]
\Mod\calC & (\Mod\calC)/\calS'\\
\Mod(\calC/G) & (\Mod(\calC/G))/(\calS'/G)
\Ar{1-1}{1-2}{"Q", bend left=15pt}
\Ar{1-2}{1-1}{"J", bend left=15pt}
\Ar{1-1}{1-2}{"\rotatebox{-90}{$\adj$}", phantom}
\Ar{1-1}{2-1}{"P\down"', bend right}
\Ar{2-1}{1-1}{"P\up"', bend right}
\Ar{1-1}{2-1}{"\adj"', phantom}
\Ar{2-1}{2-2}{"Q/G", bend left=15pt}
\Ar{2-2}{2-1}{"J/G", bend left=15pt}
\Ar{2-1}{2-2}{"\rotatebox{-90}{$\adj$}" pos=0.65, phantom}
\Ar{1-2}{2-2}{"U_m(P\down)"', bend right}
\Ar{2-2}{1-2}{"U_m(P\up)"', bend right}
\Ar{1-2}{2-2}{"\adj"', phantom}
\end{tikzcd},
\end{equation}
where by definitions of $U_m(P\down)$ and $U_m(P\up)$, the following are commutative:
{\small
\[
\begin{tikzcd}
\Mod\calC & (\Mod\calC)/\calS'\\
\Mod(\calC/G) & (\Mod(\calC/G))/(\calS'/G)
\Ar{1-1}{1-2}{"Q"}
\Ar{1-1}{2-1}{"P\down"'}
\Ar{2-1}{2-2}{"Q/G"}
\Ar{1-2}{2-2}{"U_m(P\down)"'}
\end{tikzcd},
\ 
\begin{tikzcd}
\Mod\calC & (\Mod\calC)/\calS'\\
\Mod(\calC/G) & (\Mod(\calC/G))/(\calS'/G)
\Ar{1-1}{1-2}{"Q"}
\Ar{2-1}{1-1}{"P\up"'}
\Ar{2-1}{2-2}{"Q/G"}
\Ar{2-2}{1-2}{"U_m(P\up)"'}
\end{tikzcd}.
\]}
Then the following statements hold.
\begin{enumerate}
\item
Set $\calS:= \calS' \cap \mmod\calC$ and $\calS/G:= \calS'/G \cap \mmod(\calC/G)$.
Then they turn out to be dense subcategories of $\mmod\calC$ and $\mmod(\calC/G)$,
respectively.
\item
Set $S':=S_{\Mod\calC}(\calS')$ and $S:= S_{\mmod\calC}(\calS)$ to
be the corresponding bicalculable multiplicative systems of
$\Mod\calC$ and $\mmod\calC$.
Then we have
\[
\cpt(S') = S.
\]
Hence
\[
\begin{aligned}
U_m(\Mod\calC, S')&= (\Mod \calC)/\calS',\ \text{and}\\
U_m(\cpt(\Mod\calC, S'))&= U_m(\mmod \calC, S) = (\mmod \calC)/\calS.
\end{aligned}
\]
\item
Set $S'/G:= S_{\Mod(\calC/G)}(\calS'/G)$ and $S/G:= S_{\mmod(\calC/G)}(\calS/G)$
to be the corresponding bicalculable multiplicative systems of
$\Mod(\calC/G)$ and $\mmod(\calC/G)$, respectively.
Then we have
\[
\cpt(S'/G) = S/G.
\]
Hence
\[
\begin{aligned}
U_m(\Mod(\calC/G), S'/G))&= \Mod(\calC/G)/(\calS'/G),\ \text{and}\\
U_m(\cpt(\Mod(\calC/G), S'/G))&= U_m(\mmod \calC/G, S/G) = (\mmod \calC/G)/(\calS/G).
\end{aligned}
\]
\item
We have $P\down(S') \subseteq S'/G$, $P\down(S) \subseteq S/G$ and $P\up(S'/G) \subseteq S'$.
Hence we see that $(P\down, \ph\down)$ induces a $G$-invariant functors
$(\Mod\calC, S') \to (\Mod \calC/G, S'/G)$ and
$(\mmod\calC, S) \to (\mmod \calC/G, S/G)$ and that
$P\up$ induces the right adjoint
\[
P\up \colon (\Mod\calC/G, S'/G) \to (\Mod\calC, S')
\]
to the former.
\item
The functors $U_m(P\down)$ and $U_m(P\up)$ are naturally isomorphic to the functors
$Q/G\circ P\down \circ J$ and $Q \circ P\up \circ J/G$, respectively.
\end{enumerate}
\end{lemma}

\begin{proof}
(1) Since $\mmod\calC$ is abelian by assumption, it is a dense subcategory of $\Mod\calC$.
Hence the intersection $\calS:= \calS' \cap \mmod\calC$ turns out to be a dense
subcategory of $\Mod\calC$ and hence of $\mmod\calC$.
Similarly, $\calS/G$ is shown to be a dense subcategory of $\mmod(\calC/G)$.

(2) Since $\mmod\calC$ is abelian, for any $f \in \mmod \calC(M,N)$ with $M,N \in \mmod\calC$,
we have $\Ker f, \Coker f \in \mmod\calC$.  Hence
\[
\begin{aligned}
\cpt(S')&= \{f \in \Mod\calC(M,N) \mid M,N \in \mmod\calC,\Ker f, \Coker f \in \calS'\}\\
&= \{f \in \mmod\calC(M,N) \mid M,N \in \mmod\calC,\Ker f, \Coker f \in \calS\}\\
&= S.
\end{aligned}
\]

(3)
Similarly since $\mmod(\calC/G)$ is abelian,
for any $f \in \mmod(\calC/G)(M,N)$ with $M,N \in \mmod(\calC/G)$,
we have $\Ker f, \Coker f \in \mmod(\calC/G)$.  Hence
\[
\begin{aligned}
\cpt(S'/G) &=\{f \in \Mod(\calC/G)(M,N) \mid M,N \in \mmod(\calC/G), \Ker f, \Coker f \in \calS'/G\}\\
&=\{f \in \mmod(\calC/G)(M,N) \mid M,N \in \mmod(\calC/G), \Ker f, \Coker f \in \calS/G\}\\
&=S/G.
\end{aligned}
\]

(4)
Take any $f \in S'$.
Then $f \in \Mod\calC(M,N)$ for some $M,N \in \Mod\calC$ such that
$\Ker f, \Coker f \in \calS$.
Consider the exact sequence
\[
0 \to \Ker f \ya{\si} M \ya{f} N \ya{\pi} \Coker f \to 0.
\]
Since $P\down$ is exact, we have an exact sequence
\[
0 \to P\down(\Ker f) \ya{P\down(\si)} P\down(M) \ya{P\down(f)} P\down(N) \ya{P\down(\pi)} P\down(\Coker f) \to 0.
\]
Thus we have $\Ker P\down(f) \cong P\down(\Ker f)$ and
$\Coker P\down(f) \cong P\down(\Coker f)$.
Since $\Ker f, \Coker f \linebreak[3] \in \calS$,
we have $P\down(\Ker f), P\down(\Coker f) \in \calS/G$.
Therefore, $P\down(f) \in S'/G$.
Hence we have $P\down(S') \subseteq S'/G$.
In the above, by changing $\Mod\calC$ to $\mmod\calC$,
the same argument shows that $P\down(S) \subseteq S/G$.

Next let $f \in S'/G$.
Then $f \in \Mod(\calC/G)(M, N)$ for some $M, N \in \Mod(\calC/G)$
such that $\Ker f, \Coker f \in \calS'/G$.
Consider the exact sequence
\[
0 \to \Ker f \ya{\si} M \ya{f} N \ya{\pi} \Coker f \to 0.
\]
Since $P\up$ has both a left adjoint and a right adjoint,
$P\up$ preserves small limits and small colimits,
in particular, is exact.
Therefore, $P\up$ sends the exact sequence above to the exact sequence
\[
0 \to P\up(\Ker f) \ya{P\up(\si)} P\up(M) \ya{P\up(f)} P\up(N) \ya{P\up(\pi)} P\up(\Coker f) \to 0.
\]
Thus we have $\Ker P\up(f) \cong P\up(\Ker f)$ and
$\Coker P\up(f) \cong P\up(\Coker f)$.

We now show that $P\up(\calS'/G) \subseteq \calS'$.
To show this let $X \in \calS'$.
Then $P\up(P\down(X)) \cong \DS_{a \in G} {}^aX$.
Since $\calS'$ is $G$-stable, ${}^aX \in \calS'$ for all $a \in G$.
Further since $\calS'$ is localizing, it is closed under small coproducts,
and hence,
$\DS_{a \in G} {}^aX \in \calS'$.
As a consequence, we have
\[
P\up(\{P\down(X) \mid X \in \calS'\}) \subseteq \calS'.
\]
Since $P\up$ is exact and preserves small coproducts,
$P\up$ commutes with $\Loc$, and we have
\[
\begin{aligned}
P\up(\calS'/G)&= P\up(\Loc(\{P\down(X) \mid X \in \calS'\}))\\
&= \Loc(P\up(\{P\down(X) \mid X \in \calS'\}))\\
&\subseteq \Loc(\calS') = \calS',
\end{aligned}
\]
as desired.
Therefore, since $\Ker f, \Coker f  \in \calS'/G$,
we have 
$P\up(\Ker f), P\up(\Coker f) \in \calS'$.
Thus $P\up(f) \in S'$.
Hence we have $P\up(S'/G) \subseteq S'$.

(5)
Let $\et \colon \idty_{\Mod\calC} \To J\circ Q$ be the unit of the adjoint $Q \adj J$,
and let $X \in \Mod \calC$.
Then since $\calS'$ is localizing,
$\et_X \colon X \to J(Q(X))$ is in $S'$ (see \cite[Ch.\ 4, Proposition 4.3(2)]{Po73}).
Since $P\down(S') \subseteq S'/G$ by (4) above,
\[
(Q/G\circ P\down)(\et_X) \colon (Q/G\circ P\down)(X) \to (Q/G\circ P\down \circ J)(Q(X))
\]
turns out to be an isomorphism in $\Mod(\calC/G)/(\calS'/G)$.
This yields a natural isomorphism
$((Q/G\circ P\down)(\et_X))_{X \in (\Mod\calC)/\calS'} \colon U_m(P\down) \To Q/G\circ P\down \circ J$.
The rest is proved similarly.
\end{proof}

\begin{proposition}
\label{Prp:localization-covering}
Let $\calC$ be a skeletally small $G$-category such that
$\Mod \calC$ is locally noetherian, 
both $\mmod \calC$ and $\mmod(\calC/G)$ are abelian,
and $\calS'$ a $G$-stable localizing subcategory of $\Mod\calC$.
Denote by $\calS$, $\calS/G$ the dense subcategories of $\mmod\calC$, $\mmod (\calC/G)$
induced from $\calS'$ as in Lemma \ref{lem:quotient-case}, respectively.
Then the following hold.
\begin{enumerate}
\item
The functor
\[
U(P\down, \ph\down) \colon (\mmod\CC)/\calS \to  (\mmod\CC/G)/(\calS/G)
\]
induced by the quotient $2$-functor $U = U_m$ in {\rm Example \ref{exm:quot-2-fun}}
is a $G$-precovering.
\item
If $\CC$ is a locally support-finite locally bounded category,
then
\[
U(P\down, \ph\down) \colon (\mmod\CC)/\calS \to  (\mmod\CC/G)/(\calS/G)
\]
is a $G$-covering.
\end{enumerate}
\end{proposition}

\begin{proof}
Let $\sfC$ be the full 2-subcategory of $\kCAT_m$ consisting of the objects
having the form $(\Mod \calD, T')$ for some skeletally small
$k$-category $\calD$
such that $\mmod \calD$ is abelian, and
$T' = S_{\calC}(\calT')$ for some localizing subcategory
$\calT'$ of $\Mod\calC$, and consider the 2-category
$\hat{\sfC}$ defined in Definition \ref{dfn:2ext-gen}.
Then $(\Mod\calC, S') \in \hat{\sfC}$, and here
$\calC$ is a skeletally small
$G$-category
$S' = S_{\calC}(\calS')$ for some $G$-stable localizing subcategory
$\calS'$ of $\Mod\calC$.
By Lemma \ref{lem:quotient-case} (4),
$(P\down, \ph\down) \colon (\Mod \calC, S') \to (\Mod(\calC/G), S'/G)$
has an adjunction $\om \colon P\down \adj P\up$, it induces a $G$-precovering functor
$(P\down, \ph\down) \colon (\mmod \calC, S) \to (\mmod(\calC/G), S/G)$
by the adjunction $\om$.
Then we can apply Theorem \ref{thm:general-case}(2)(a)
to have a $G$-precovering
\[
U(P\down, \ph\down) \colon \cpt(\Mod\calC/\calS) \to \cpt(\Mod(\calC/G)/(\calS'/G))
\]
by the adjunction $U(\om)\colon U(P\down) \adj U(P\up)$.
It remains to show that
$\cpt(\Mod\calC) \subseteq \cpt(\Mod\calC/\calS')$,
or equivalently,
\begin{equation}
\label{eq:Q-pres-cpt}
Q(\mmod\calC) \subseteq \cpt((\Mod\calC)/\calS').
\end{equation}
Since by assumption $\Mod \calC$ is locally noetherian,
$J$ preserves small coproducts by \cite[\S 19, Corollary 13]{G1}, and hence by Lemma \ref{lem:left-adj-cpt},
the inclusion \eqref{eq:Q-pres-cpt} holds.
%
%
%
%
%
%
%
%
%
%
\end{proof}

\section{$G$-precoverings for functor categories}\label{Subsection 3.1}

In this section, we show that the canonical covering
$(P, \ph) \colon \calC \to \calC/G$ induces a $G$-precovering
$H(P\down, \ph\down) \colon H(\mmod \CC) \to H(\mmod (\CC/G))$
between morphism categories of $\mmod \CC$ and $\mmod(\CC/G)$,
respectively (Proposition \ref{Prop-pushdown-adj}).
We apply Theorem \ref{thm:factor-case} to $H(P\down, \ph\down)$
to obtain a $G$-precovering between factor categories
\[
\udl{\rmH((P\down, \phi\down))} \colon
\frac{\rmH(\Mod \CC)}{\ang{\CU_{\CC}}}
\to \frac{\rmH(\Mod(\CC/G))}{\ang{\CU_{\CC/G}}}
\]  
(see Proposition \ref{prp:H-factor}).
Moreover, we show that this yields a $G$-precovering between finitely presented functors
\[
\CF(P\down) \colon \CF(\mmod \CC) \to \CF(\mmod(\CC/G))
\]
(see Theorem \ref{thm:G-Fp-precovering}).

\subsection{$G$-precoverings for morphism category}

\begin{definition}
\label{dfn:mor-cat}
Let $\CA$ be a category. 
Then the {\em morphism category} $\rmH(\CA)$ of $\CA$ is defined as follows:
 
The objects of $\rm{H}(\CA)$ are the morphisms $f \colon x \to y$ in $\CA$.
We sometimes present it as $\vertmap{x}{f}{y}$ to visualize the situation,
and we set $\dom(f):= x, \cod(f):= y$. 

For objects $f \colon x \to y$ and $g \colon u \to v$ of $\rmH(\CA)$,
the set of morphisms $\rmH(\CA)(f,g)$ from $f$ to $g$ is defined by
\[
\rmH(\CA)(f,g):= \{(p,q) \in \CA(x,u) \times \CA(y,v)\mid qf = gp\}.
\]
The defining condition is expressed by the commutativity of the diagram
\[
\begin{tikzcd}
x & u\\
y & v
\Ar{1-1}{1-2}{"p"}
\Ar{2-1}{2-2}{"q"'}
\Ar{1-1}{2-1}{"f"'}
\Ar{1-2}{2-2}{"g"}
\end{tikzcd}.
\]
This morphism $(p,q) \colon f \to g$ is sometimes denoted by
$\vertmap{x}{f}{y} \longto{p}{q} \vertmap{u}{g}{v}$
to visualize the situation.

For morphisms $(p,q) \colon f \to g$ and $(r,s) \colon g \to h$ in
$\rmH(\CA)$, the composite $(r,s)(p,q)$ is defined by component-wise:
$(r,s)(p,q):= (rp, sq)$.
\end{definition}

\begin{example}
Consider the case that $\CA=\mmod \La$ for a finite-dimnsional algebra $\La$.
Set $T_2(\La):= \tiny \bmat{\La & \La\\ 0 & \La}$, the upper triangular
matrix algebra over $\La$.
Then $\rm{H}(\mmod \La)$ is equivalent to the category $\mmod T_2(\La)$
of finitely generated right modules over $T_2(\La)$.
Indeed, the equivalence is given by regarding each object
$f \colon X \to Y$ of $\rmH(\mmod \La)$ as a right $\La$-module
$X \ds Y$ with the right $\La$-action defined by
$(x, y) \bmat{a & b \\ 0 & c}:= (xa, f(x)b+yc)$
for all $(x,y) \in X \ds Y$ and $\bmat{a & b \\ 0 & c} \in T_2(\La)$.

\end{example}

\begin{lemma}
Let $\sfC$ be either $\kCat$ or $\kCAT$.
Then the correspondence $\rmH \colon \obj(\sfC) \to \obj(\sfC)$,
$\CC \mapsto \rmH(\CC)$ defined in Definition \ref{dfn:mor-cat}
is extended to a $2$-functor $\rmH \colon \sfC \to \sfC$.
\end{lemma}

\begin{proof}
We extend the map $\rmH$ to maps between 1-morphism and 2-morphisms as follows.

{\bf On $1$-morphisms}
Let $F \colon \calC \to \calD$ be a $1$-morphism in $\sfC$.
Then we define
$\rmH(F) \colon \rmH(\calC) \to \rmH(\calD)$ to be a $k$-functor sending
a morphism $\vertmap{x}{f}{y} \longto{p}{q} \vertmap{u}{g}{v}$ in $\rmH(\calC)$
to a morphism
$\vertmap{F(x)}{F(f)}{F(y)} \longto{F(p)}{F(q)} \vertmap{F(u)}{F(g)}{F(v)}$
in $\rmH(\calD)$.
Then it is obvious that $\rmH(F)$ in fact becomes a $k$-functor, i.e., a $1$-isomorphism in $\sfC$.

{\bf On $2$-morphisms}
Let $E, F \colon \calC \to \calD$ be $1$-morphisms in $\sfC$, and
$\al\colon E \To F$ a 2-morphism in $\sfC$.
Then we define $\rmH(\al)$ by $\rmH(\al):= (\rmH(\al)_{f})_{f \in \rmH(\calC)}$,
where we set $\rmH(\al)_{f}:= (\al_x, \al_y)$, which is visualized as
\[
\vertmap{E(x)}{E(f)}{E(y)} \longto{\al_x}{\al_y} \vertmap{F(x)}{F(f)}{F(y)}
\]
for all objects $\vertmap{x}{f}{y}$ of $\rmH(\calC)$.
Then $\rmH(\al)\colon \rmH(E) \To \rmH(F)$ is a natural transformation, i.e., a $2$-morphism in $\sfC$.
Indeed, let $\vertmap{x}{f}{y} \longto{p}{q} \vertmap{u}{g}{v}$ be a morphism in $\rmH(\calC)$.
We have to show the commutativity of the diagram
\[
\begin{tikzcd}[column sep=50pt]
E(f) & E(g)\\
F(f) & F(g)
\Ar{1-1}{1-2}{"{(E(p),E(q))}"}
\Ar{2-1}{2-2}{"{(F(p),F(q))}"'}
\Ar{1-1}{2-1}{"{(\al_x,\al_y)}"'}
\Ar{1-2}{2-2}{"{(\al_u,\al_v)}"}
\end{tikzcd}.
\]
But this immediately follows from the naturality of $\al$.
It is easy to verify that $\rmH$ preserves 2-identities, vertical compositions, and horizontal compositions.
\end{proof}

\begin{lemma}
\label{lem:H(action)}
Let $\calC = (\calC, A)$ be a $G$-category.
Then $\rmH(\calC) = (\rmH(\calC), \rmH(A))$ turns out to be a $G$-category by setting
$\rmH(A)_a$ to be the automorphism $\rmH(A_a) \colon \rmH(\calC) \to \rmH(\calC)$ for all $a \in G$.
Apply this construction to the $G$-category
$\Mod \calC = (\Mod \calC, \Mod A)$.
Then $\rmH(\Mod \calC) = (\rmH(\Mod \calC), \rmH(\Mod A))$ becomes again a $G$-category.
\end{lemma}

\begin{proof}
Straightforward.
\end{proof}

\begin{lemma}
\label{lem:H(G-eqvar)}
Let $\CC =(\calC, A)$ and $\CC'= (\calC', A')$ be $G$-categories, and $(F, \ph) \colon \calC \to \calC'$ a $G$-equivariant functor.
Then we can define a $G$-equivariant functor
$\rmH(F, \ph) \colon \rmH(\calC) \to \rmH(\calC')$ by setting
$\rmH(F, \ph):= (\rmH(F), \rmH(\ph))$. 
Here, $\rmH(\ph)_a \colon \rmH(A'_a)\rmH(F) \To  \rmH(F)\rmH(A_a)$
is defined by $\rmH(\ph)_a:= (\rmH(\ph)_{a,f})_{f \in \rmH(\CC)}$, where
\[
\rmH(\ph)_{a, f}:= (\ph_{a, x}, \ph_{a, y}) \colon
\vertmap{(A'_aF)(x)}{F(f)}{(A'_aF)(y)}
\longto{\ph_{a,x}}{\ph_{a,y}}
\vertmap{(FA_a)(x)}{(FA_a)(f)}{(FA_a)(y)}
\]
for all objects $f \colon x \to y$ of $\rmH(\calC)$.

In particular, by considering the case that $\calC'$ above has the trivial
$G$-action, $\rmH$ maps a $G$-invariant functor
$(F, \ph) \colon \calC \to \calC'$
to the $G$-invariant functor 
$\rmH(F, \ph) \colon \rmH(\calC) \to \rmH(\calC')$.
\end{lemma}

\begin{proof}
Straightforward.
\end{proof}

By Lemmas \ref{lem:H(action)} and \ref{lem:H(G-eqvar)},
we can prove the following.

\begin{proposition}
The $2$-functors $\rmH \colon \kCat \to \kCat$ and $\rmH \colon \kCAT \to \kCAT$
induce $2$-functors $\rmH \colon \GCat \to \GCat$ and $\rmH \colon \GCAT \to \GCAT$.
\qed
\end{proposition}

Let $\CC = (\calC, A)$ be a skeletely small  $G$-category. Apply Lemma \ref{lem:H(G-eqvar)} to the canonical $G$-invariant functor
$(P\down, \ph\down) \colon \Mod \calC \to \Mod (\calC/G)$ to obtain
a $G$-invariant functor
\[
\rmH((P\down, \ph\down))= (\rmH(P\down), \rmH(\ph\down)) \colon \rmH(\Mod \calC) \to \rmH(\Mod (\calC/G)).
\]
We also have a functor
\[
\rmH(P\up) \colon \rmH(\Mod (\calC/G)) \to \rmH(\Mod \calC)
\]
induced from the pull-up functor $P\up\colon \Mod (\calC/G) \to \Mod \calC$.

The following is immediate from
Lemma \ref{lem:2-fun-pres-adj}.

\begin{lemma}\label{lem. left-morphism}
Let $\calC$ and $\calD$ be $k$-categories, and
$L \colon \calC \to \calD$ a left adjoint to $R \colon \calD \to \calC$.
Then $\rmH(L) \colon \rmH(\calC) \to \rmH(\calD)$ turns out to be a left adjoint to
$\rmH(R) \colon \rmH(\calD) \to \rmH(\calC)$.
\end{lemma}



\begin{remark}\label{Rem-left-adj-morphism}
In the setting of the lemma above,
let $\om:= (\om_{x,y} \colon \calD(Lx, y) \to \calC(x, Ry))_{x\in \calC, y \in \calD}$
be the adjunction.
Then the adjunction
\[
\rmH(\om):= (\rmH(\om)_{f,g}\colon \rmH(\calD)(\rmH(L)f, g) \to \rmH(\calC)(f, \rmH(R)g))_{f\in \rmH(\calC), g \in \rmH(\calD)}
\]
is given by
\[
\rmH(\om)(p,q):= (\om(p), \om(q)),
\]
which is visualized as
\[
\left(\vertmap{L(x)}{L(f)}{L(y)} \longto{p}{q} \vertmap{u}{g}{v}\right)
\mapsto
\left(\vertmap{x}{f}{y} \longto{\om(p)}{\om(q)} \vertmap{R(u)}{R(g)}{R(v)}\right).
\]
\end{remark}

The following is an immediate consequence of the lemma above.

\begin{lemma}\label{adjointMorphim}
The functor $\rmH(P\down)$ is a left adjoint to $\rmH(P\up)$.
\end{lemma}

The following is immediate from the definition of $\rmH$, and is used without reference.

\begin{lemma}\label{lem.comp-morphism}
For each $k$-category $\calA$, we have
$\rmH(\cpt(\calA)) = \cpt(\rmH(\calA))$.
\end{lemma}

\begin{theorem}
\label{thm:mor-cat-case}
 Let  $(L, \ph) \colon \calA \to \calB$  be a  $G$-invariant functor that  has  a right adjoint $R \colon \calB \to \calA$,  with an adjunction isomorphism $\om \colon L \adj R$. 
If $(L, \ph) \colon \calA \to \calB$ induces a $G$-precovering $\cpt(\calA) \to \cpt(\calB)$
by an adjoint $\om \colon L \adj R$, then $\rmH(L, \ph) \colon \rmH(\calA) \to \rmH(\calB)$ induces a $G$-precovering
\[
\rmH(L,\ph) \colon \rmH(\cpt(\calA)) \to \rmH(\cpt(\calB))
\]
by the adjoint $H(\om) \colon H(L) \adj H(R)$.
\end{theorem}

\begin{proof}
Let $f, g$
be objects of $\rm{H}(\cpt(\calA))$.
Then we have the canonical isomorphism
\[
\Ds_{a \in G}\rmH(\cpt(\calA))(f, {}^ag) \ya{\sim} \rmH(\calA)(f, \Ds_{a\in G}{}^ag)
\]
because both $\dom(f)$ and $\cod(f)$ belong to $\cpt(\calA)$.

(PA1) Since $\rmH \colon \kCAT \to \kCAT$ is a 2-functor,
the natural isomorphism $$t \colon \Ds_{a\in G}A_a \To RL$$
yields a natural isomorpism
$\rmH(t) \colon \Ds_{a\in G}\rmH(A_a) \To \rmH(R)\rmH(L)$.

(PA2) The commutativity of the diagram
\[\begin{tikzcd}[column sep=55pt]
	{\bigoplus_{a\in G}\rmH(\cpt(\calA))(f, {}^a g)} & {\rmH(\calA)(f,\bigoplus_{a\in G}{}^a g)} \\
	{\rmH(\cpt(\calB))(\rmH(P\down)(f), \rmH(P\down)(g))} & {\rmH(\calB)(f, \rmH(P\up)\rmH(P\down)(g))}
	\arrow["\sim",from=1-1, to=1-2]
	\arrow["\sim", "{\rmH(\om)_{f,\rmH(P\down)(g)}}"', from=2-1, to=2-2]
	\arrow["{\rotatebox{90}{$\sim$}}", from=1-2, to=2-2]
	\arrow["{\rmH(L,\phi)^{(2)}_{f,g}}"', from=1-1, to=2-1]
\end{tikzcd}
\]
is shown by using the same commutative diagrams for domains and codomains of $f$ and ${}^ag$
($a \in G)$.

(PA3) $\rmH(R)$ preserves small coproducts because so does $R$.
%
%
\end{proof}

By applying Theorem \ref{thm:mor-cat-case} to
the pushdown functor $(P\down, \ph\down) \colon \Mod \calC \to \Mod (\calC/G)$ for
some small $G$-category $\calC$, we obtain the following.

\begin{proposition}
\label{Prop-pushdown-adj}
The $G$-invariant functor $\rmH(P\down, \ph\down):\rmH(\Mod \CC)\to \rmH(\Mod (\CC/G))$ induces a $G$-precovering 
$$\rmH(P\down, \ph\down):\rmH(\mmod \CC)\to \rmH(\mmod (\CC/G))$$
by the adjunction $\rmH(\theta)$ (see Subsection \ref{Pushdownfunctor} for the definition of $\theta$).
\end{proposition}

We have the following remarks, which are frequently used later.
The first one is mentioned in the proof of Theorem \ref{thm:mor-cat-case},
and the second one is immediate from the corresponding property of $P\down$.

\begin{remark}\label{RemarkFacts}
Let $f$ 
be an object of $\rm{H}(\Mod \CC)$.
Then the following hold.
\begin{enumerate}
\item
$\rmH(P\up)\rmH(P\down)(f) \cong \Ds_{a \in G}{}^af$; and
\item
$\rmH(P\down)(f) \cong \rmH(P\down)({}^af)$ for all $a \in G$.
\end{enumerate}
\end{remark}

\begin{lemma}\label{Lemma 3.1}
Let $\CC$ be a locally bounded $k$-category and $G$ a torsion-free group
acting freely on $\CC$.
Let $f \colon X \to Y$ be a nonzero object in $\rm{H}(\mmod\CC)$.
If ${}^a\!f\simeq f$ for some $a \in G$, then $a=1$. 
\end{lemma}

\begin{proof}
First note that both $X$ and $Y$ are finite-dimensional
by Remark \ref{rmk:fg-lb}
because $\CC$ is locally bounded and $X$ and $Y$ are finitely presented.
Since $f \ne 0$, either $X \ne 0$ or $Y \ne 0$.
Consider the case that $X \ne 0$.
Then we have $X(x) \ne 0$ for some $x \in \CC$.
Assume that ${}^a\!f\simeq f$.
Then ${}^{a^{-n}}X \simeq X$,
in particular, $X(a^nx) \ne 0$ for all positive integers $n$.
Then $a^n x = a^m x$ for some positive integers $m, n$ with $m < n$
because $X$ is finite-dimensional.
Since the $G$-action on $\CC$ is free,
$a^{n-m}x = x$ shows that $a^{n-m} = 1$ with $n - m >0$.
Hence $a = 1$ because $G$ is torsion-free.
The proof for the case that $Y \ne 0$ is similar.
\end{proof}

The lemma above says that $G$ acts freely on $\rm{H}(\mmod\CC)$
up to isomorphisms.

\begin{lemma}\label{lemma 3.5}
If $f:X\rt Y$ is an object of $\rmH(\Mod (\CC/G))$,
then $\rmH(P\up)(f)\cong {}^a\rmH(P\up)(f)$ for all $a \in G$.
\end{lemma}

\begin{proof}
Take $a \in G$. The natural isomorphism ${\ph\down}_a \colon P\down \to P\down \circ \ovl{A}_{a}$ follows the following isomorphism in $\rmH(\Mod \CC/G)$

\[
\vertmap{P\down X}{P\down f}{P\down Y} \longto{{\ph\down}_{a, X}}{{\ph\down}_{a, Y}} \vertmap{P\down {}^aX}{P\down{}^af}{P\down {}^aY}
\]
This implies the desired isomorphism.
\end{proof}

If $\CC$ is a locally bounded category, then it is well-known that
the categories $\mmod (\CC/G)$ and $\mmod \CC$ are
Krull-Schmidt categories, and by ?? so are the morphism categories
$\rmH(\mmod \CC)$ and $\rmH(\mmod (\CC/G))$.

\begin{proposition}\label{Prop-pushdown-ind}
Let $\CC$ be a locally bounded $k$-category,
$G$ a torsion-free group acting freely on $\CC$,
and $f, f_1, f_2$ objects of $\rmH(\mmod \CC)$.
Then the following hold.
\begin{enumerate}			
\item
If $f$ is indecomposable, then so is $\rmH(P\down)(f)$.
\item
If objects $f_1$ and $f_2$ are indecomposable,
then $\rmH(P\down)(f_1)\cong \rmH(P\down)(f_2)$	implies
$f_1\cong {}^a\!f_2$,  for some $a \in G$.
\end{enumerate}
		
\end{proposition}
\begin{proof}
(1)
Assume that $\rmH(P\down)(f)=C \oplus C'$ with $C \ne 0$.
By Remark \ref{RemarkFacts}, we have
\[
\Ds_{a \in G} {}^a\!f \cong \rmH(P\up)(\rmH(P\down)(f)) 
\cong \rmH(P\up)(C)\oplus \rmH(P\up)(C').
\]
Since ${}^a\!f$ are pairwise non-isomorphic indecomposable objects
by Lemma \ref{Lemma 3.1}, we have
$\rmH(P\up)(C)=\Ds_{b \in V}{}^b\!f$ for some $V \subseteq G$
by \cite[Theorem E.1.26]{Pe}.
We know that $\rmH(P\up)(C) \cong {}^a\rmH(P\up)(C)$ for any $a \in G$
by Lemma \ref{lemma 3.5}.
Hence $\Ds_{c \in V}{}^c\!f \cong \Ds_{b \in V}{}^{ab}f$
for all $a \in G$.
By the Krull-Remak-Schmidt-Azumaya Theorem \cite[Theorem E.1.24]{Pe}, 
for each $a \in G$ and $b \in V$, there is $c \in V$ such that ${}^{ab}f\simeq {}^{c}f$.
Thanks to Lemma \ref{Lemma 3.1}, $c=ab$.
Therefore $a = cb\inv \in V$ for all $a \in G$, and hence $V=G$,
which shows that $\rmH(P\up)(C') = 0$.
Therefore $C'=0$.

 (2)
 Assume $\rmH(P\down)(f_1)\cong \rmH(P\down)(f_2)$.
 Applying the functor $\rmH(P\up)$, we obtain
 \[
\Ds_{a \in G} {}^a\!f_1 \cong \Ds_{a \in G} {}^a\!f_2
 \]
by Remark \ref{RemarkFacts}.
Hence $f_1$ is a direct summand of $\Ds_{a \in G} {}^a\!f_2$.
Since ${}^a\!f_2$ are indecomposable module with local endomorphism algebras, the Krull-Schmidt property of $\rm{H}(\mmod \CC)$ shows the existence of $a \in G$ such that $f_1\simeq {}^af_2$, as desired.
\end{proof}

\subsection{A relative version}

Throughout the rest of this paper,
we assume that $\CC$ is a \emph{skeletally small $G$-category}.

\begin{definition}
Let $\CM$ be a full subcategory of $\mmod \CC$, and
$\CN$ a full subcategory of $\mmod (\CC/G)$.
\begin{enumerate}
\item
$\CM$ is called {\it $G$-stable} if ${}^a\!\CM = \CM$ for all $a \in G$.
If this is the case,
then we can define a $G$-action on $\CM$.
This $G$-action is said to be \emph{free up to isomorphisms} 
(or $G$ \emph{acts freely} on $\CM$ \emph{up to isomorphisms}) in case
${}^aM \cong M$ for some $0 \ne M \in \CM$ and $a \in G$ implies that $a = 1$.
\item
We denote by $P\down(\CM)$ the full subcategory of $\mmod (\CC/G)$ consisting of all objects $V$ such that $V$ is isomorphic to $P\down(U)$ for some $U \in \CM$.
Note that $P\down(\CM)$ is closed under isomorphisms.
\item
We denote by $P^{-1}\down(\CN)$ the full subcategory of $\mmod \CC$ consisting of all objects $U$ such that $P\down(U)$ is isomorphic to an object in $\CN$.
Note that $P\down\inv(\CN)$ is closed under isomorphisms.

\end{enumerate}
\end{definition}

Based on Example \ref{ex-G-stable}, 
our results above concerning the morphism category
of $\Mod \CC$
can be generalized to the morphism category $\rmH(\CK)$ of $G$-stable subcategories $\CK$ of $\Mod \CC$.
Since their proofs are similar to those already written above, (just restrict the arguments to the subcategories), we state them here without proofs.
The following lemma is useful:

\begin{lemma}\label{lem-inverse-image}
The following hold.
\begin{enumerate}
\item
Let $\CN$ be a full subcategory of $\mmod (\CC/G)$.
Then $P\down^{-1}(\CN)$ is a $G$-stable subcategory of $\mmod \CC$.
\item
Assume that $\CC$ is a locally bounded $k$-category and
$G$ acts freely on $\ind \CC$ (Definition \ref{dfn:ind-C}).
If $\CM$ is a full subcategory of $\mmod \CC$
closed under direct summands, then so is $P\down(\CM)$.
\end{enumerate}
\end{lemma}

\begin{proof}
(1)
Take $X \in P\down^{-1}(\CN)$ and $a \in G$.
Then, by definition, $P\down(X) \cong N$ for some $N \in \CN$.
On the other hand, we have an isomorphism ${\phi\down}_{a, X}: P\down(X) \to P\down({}^aX)$ (see Sect.~\ref{Pushdownfunctor}).
Therefore, $P\down({}^aX) \cong N$, and ${}^aX \in P\down^{-1}(\CN)$, as desired.

(2)
Let $N \in P\down(\CM)$ with a decomposition
$N = U \oplus V$ in $\mmod (\CC/G)$.
By assumption there exists an $M \in \CM$
such that  $P\down(M) \cong N$.
Since $\mmod \CC$ is a Krull-Schmidt category,
we can write $M =\Ds_{i \in I} M_i$,
where $M_i$ are indecomposable objects and $I$ is a finite set.
Hence $ U \oplus V = N \cong P\down(M)=\Ds_{i \in I} P\down(M_i)$.
By Lemma \ref{Inde-Pushdown}, all $P\down(M_i)$ are indecomposable.
Here, since $\mmod (\CC/G)$ is also a Krull-Schmidt category,
there exists a subset $J \subseteq I$ such that
$U = \Ds_{i \in J}P\down(M_i) \cong P\down(\Ds_{i \in J}M_i)$.
Since $\Ds_{i \in J}M_i$ is a direct summand of $M$ and
$\CM$ is closed under direct summand, $\Ds_{i \in J}M_i$ is in
$\CM$.
Hence $U \in P\down(\CM)$, as desired.
\end{proof}

\begin{proposition}\label{Proposition3last}
Let $\CK$ be $G$-stable subcategories of $\mmod \CC$, and $\CK'$ a subcategory of $\mmod \CC/G$.
Assume $P\down$ sends $\CK$ to $\CK'$ and $P\up$ sends $\CK'$ to ${\rm Add}\mbox{-}\CK,$ e.g., $\CK'=P\down(\CK)$.  Then the following hold.
\begin{enumerate}
\item 
The $G$-invariant functor $\rmH(P\down|, \ph\down|)\colon \rmH({\rm  Add}\mbox{-}\CK)\rt \rmH({\rm Add}\mbox{-}\CK')$ induces a $G$-precovering $$\rmH(P\down|, \ph\down |)\colon \rmH(\CK)\rt \rmH(\CK')$$ with the adjuction  $\rmH(\theta|)$.
\item 
Assume that $\CC$ is a locally bounded category
and $G$ is a torsion-free group
acting freely on $\CC$. Then
\begin{enumerate}			
\item
If $f$ is an indecomposable object in $\rm{H}(\CK)$,
then so is $\rmH(P\down)(f)$ in $\rmH(\CK')$.
\item
If objects $f_1$ and $f_2$ are indecomposable objects in $\rmH(\CK)$, then $\rmH(P\down)(f_1) \cong \rmH(P\down)(f_2)$
implies $f_1 \cong {}^a\!f_2$ for some $a \in G$. \qed
\end{enumerate}
\end{enumerate}
\end{proposition}

Actually, Propositions \ref{Prop-pushdown-adj} and \ref{Prop-pushdown-ind} are special cases of the above proposition when we take $\CK=\mmod \CC$ and $\CK'=\mmod \CC/G$. Here are   some additional  examples of $\CK \subseteq \mmod \CC$ that satisfy the assumptions of Proposition \ref{Proposition3last}.

\begin{example}
\begin{enumerate}
\item
Let  $\CK = \add\{{}^aM \mid a \in G\}$ for some
$M \in \mmod \CC$, and $\CK':= P\down(\CK)$.

\item
Assume that $\CC$ is a locally support-finite $G$-category with $G$-action free on $\ind \CC$. Let  $\CK=\CG\!p\mbox{-}\CC$ and  $\CK'=\CG\!p\mbox{-}(\CC/G)$ 
be the full subcategory
 $\Mod \CC$ and  $\Mod(\CC/G)$, respectively,  consisting of  all finitely generated  Gorenstein projective functors. Note that $\CG\!p\mbox{-}\CC$ is $G$-stable since    the equivalence $\ovl{A}_{a}$ for every $a \in G$ preserves  Gorenstein projective functors. Moreover,  due to  \cite[Theorem 4.5]{AHV}, we get $P\down(\CG\!p\mbox{-}\CC) \subseteq \CG\!p\mbox{-}(\CC/G)$, and by \cite[Lemma 4.2 and Theorem 4.9]{AHV}, we can deduce that $P \up(\CG\!p\mbox{-}(\CC/G)) \subseteq {\rm Add}\mbox{-}\CG\!p\mbox{-}\CC$. Therefore, the pushdown and pullup functors have the required conditions respect to the finitely generated Gorenstein projective functors.

 \item 
Let $\CA$ be an abelian category. A full subcategory $\CU$ of $\CA$
  is \emph{$d$-cluster-tilting} if it is functorially finite in $\CA$ and
\[
\begin{aligned}
\mathcal{U}&=\{X\in\CA \mid\Ext_{\CA}^i(U,X)=0\ \mbox{ for all }\ 1\le i\le d-1\},\\
&=\{X\in\CA \mid\Ext_{\CA}^i(X,U)=0\ \mbox{ for all }\ 1 \le i\le d-1\}.
\end{aligned}
\]
 Moreover, a finitely generated module $M\in \CA$ is called a \emph{$d$-cluster-tilting
  module} if its additive closure $\text{add}\mbox{-} M$ forms a $d$-cluster-tilting subcategory of
  $\CA$.
Assume that $\CC$ is as in $(2)$. 
Let $\CK':=\CM$ be a  $d$-cluster-tilting subcategory of $\mmod \CC/G$. By Lemma \ref{lem-inverse-image} (1), $P\down^{-1}(\CM)$ is a $G$-stable subcategory in $\mmod \CC$. Take $\CK:=P\down^{-1}(\CM)$. Then,  $\CK$ and $\CK'$ have the required conditions in the proposition. In addition, by \cite[Theorem 2.14]{DI}, $\CK$ and $\CK'$ are $d$-cluster-tilting subcategories in $\mmod \CC$ and $\mmod \CC/G$, respectively.
 
\end{enumerate}
\end{example}

\subsection{$G$-precoverings for factor categories of morphism categories }\label{Subsection 3.2}\label{sub-G-precovering}

The most important purpose of this section is to explain how the 
functor $\rmH(P\down, \ph\down) \colon \rmH(\mmod \CC) \to \rmH(\mmod (\CC/G))$ defined in Proposition \ref{Prop-pushdown-adj} is applied to obtain $G$-precovering functors from certain  factor categories of $\rmH(\mmod \CC)$ to those of $\rmH(\mmod(\CC/G))$.

First, we need the following definition to construct an additive factor category of the given morphism category.

\begin{definition}
For a category $\CM$, we denote by $\CU_\CM$ the class of objects of $\rmH(\CM)$ having the form $(X \to 0)$ or $(X \ya{1} X)$ for some object $X$ of $\CM$.
\end{definition}

We intend to apply Theorem \ref{thm:general-case} $(2)(b)$ to obtain our main result in this subsection. For this purpose, we consider $\sfC_h$ as a full $2$-subcategory of $\kCAT_s$, where the objects are of the form $(\rmH(\CM), \CU_{\CM})$, where  $\CM$ is an additive category having small coproducts.   We define 2-functor $U_h:\sfC_h\rt \kCAT$ by sending an object $(\rmH(\CM), \CU_{\CM})$ to the factor category $\frac{\rmH(\CM)}{\ang{\CU_{\CM}}}$. Furthermore, we denote by $\pi_{\CM}$ the canonical functor $\rmH(\CM) \to \frac{\rmH(\CM)}{\ang{\CU_{\CM}}} $,  and  the family
$\pi_h:= (\pi_{\CM})_{(\rmH(\CM), \CU_{\CM}) \in \sfC_h}$
defines a strict $2$-natural transformation
\[
\begin{tikzcd}
\makebox[1em][r]{$\sfC_h$} & \makebox[1em][l]{$\kCAT$}
\Ar{1-1}{1-2}{"\si", ""'{name=a}, bend left}
\Ar{1-1}{1-2}{"U_h"', ""{name=b}, bend right}
\Ar{a}{b}{"\pi_h", Rightarrow}
\end{tikzcd}\qquad,
\]
where $\si$ as usual  is the inclusion $2$-functor. The above diagram is extended to the following diagram
\[
\begin{tikzcd}
\widehat{\sfC_h} & \makebox[1em][l]{$\GCAT$}
\Ar{1-1}{1-2}{"\hat{\si}", ""'{name=a}, bend left}
\Ar{1-1}{1-2}{"\widehat{U_h}"', ""{name=b}, bend right}
\Ar{a}{b}{"\widehat{\pi_h}", Rightarrow}
\end{tikzcd}\qquad,
\]
where $\widehat{\sfC_h}, \widehat{U_h}, \widehat{\si}$ are defined as in Lemma \ref{lem:ext-2-fun-gen}.

To provide a better understanding, let us explain the construction of $\widehat{U_h}$ in the case where $\CM=\Mod \CC$, and $\CC=(\CC, A)$ is a skeletally small $G$-category. For simplicity, we use $\CU_{\CC}$ and $\CU_{\CC/G}$ instead of $\CU_{\Mod \CC}$ and $\CU_{\Mod(\CC/G)}$, respectively. 
By applying Lemma \ref{lem:H(G-eqvar)}
to the $G$-invariant functor $(P\down, \ph\down) \colon \Mod\CC \to \Mod(\CC/G)$, we obtain a $G$-invariant funcor
\[
\rmH((P\down, \ph\down)) = (\rmH(P\down), \rmH(\ph\down)) \colon \rmH(\Mod\CC) \to \rmH(\Mod(\CC/G)),
\]
where 
\[
\rmH({\ph\down}_a) \colon \rmH(P\down) \To \rmH(P\down)\circ \rmH(A_a)
\]
is a natural isomorphism for all $a \in G$.
Let $a \in G$. Then the automorphism $A_a$ of $\CC$ defines an automorphism
${}^a(\blank)$ of $\Mod \CC$, which induces an automorphism $\rmH({}^a(\blank))$ of $\rmH(\Mod\CC)$.
Furthermore, the latter induces an automorphism $\udl{\rmH({}^a(\blank))}$ of the factor category
$\frac{\rmH(\Mod \CC)}{\ang{\CU_{\CC}}}$
as it preserves the ideal $\ang{\CU_{\CC}}$.
Then by sending $a$ to $\udl{\rmH({}^a(\blank))}$ we can define a $G$-action on $\frac{\rmH(\Mod \CC)}{\ang{\CU_{\CC}}}$, namely
we have ${}^a\udl{f} = \udl{{}^af}$
for all $a \in G$ and $f \in \rmH(\Mod\calC)$. Since every additive functor preserves isomorphisms and zero objects, the pushdown functor $\rmH(P\down)$ sends the ideal $\ang{\CU_{\CC}}$ into $\ang{\CU_{\CC/G}}$. Hence, we obtain the induced functor
\[
\udl{\rmH((P\down, \phi\down))} \colon
\frac{\rmH(\Mod \CC)}{\ang{\CU_{\CC}}}
\to \frac{\rmH(\Mod(\CC/G))}{\ang{\CU_{\CC/G}}}.
\]
The induced functor $\udl{\rmH((P\down, \phi\down))}$ is precisely $\widehat{U_h}(\rmH((P\down, \phi\down))$.
We can continue the same observation for the pullup functor $P\up:\Mod \CC/G\rt \Mod \CC $. Similarly, we obtain the induced functor 
\[
\udl{\rmH(P\up)} \colon
\frac{\rmH(\Mod \CC/G)}{\ang{\CU_{\CC/G}}}
\to \frac{\rmH(\Mod(\CC))}{\ang{\CU_{\CC}}},
\]
which is indeed $\widehat{U_h}(\rmH(P\up))$.

By Lemma \ref{lem. left-morphism} in conjunction with Remark \ref{Rem-left-adj-morphism}, we can infer  that $G$-invariant functor $\rmH(P\down)\colon \rmH(\Mod\CC) \to \rmH(\Mod \CC/G)$ has the right adjoint $\rmH(P\up)\colon \rmH(\Mod \CC/G) \to \rmH(\Mod \CC)$ with the adjucntion $\rmH(\theta) \colon \rmH(P\down) \adj \rmH(P\up)$ in $\kCAT$. As explained earlier, since
\[
\rmH(P\down)(\CU_{\CC}) \subseteq  \CU_{\CC/G} \ \  \text{and} \ \ \rmH(P\up)(\CU_{\CC/G})\subseteq \CU_{\CC},
\]
the $G$-invariant functor $\rmH(P\down, \phi\down) \colon (\rmH(\Mod \CC),  \rmH(A), \CU_{\CC}) \to (\rmH(\Mod (\CC/G)), \CU_{\CC/G})$ has the right adjoint $\rmH(P\up): (\rmH(\Mod (\CC/G))), \CU_{\CC/G})\to  (\rmH(\Mod \CC),  \CU_{\CC})$ with the adjucntion $\rmH(\theta) \colon \rmH(P\down) \adj \rmH(P\up)$ in $\sfC_h.$

Now, we are ready to apply Theorem \ref{thm:general-case} to the 2-functor $U_h$ defined above in order to prove our main result in this subsection. For the sake of simplicity, we denote $\widehat{U_h}(F)$ by $\underline{F}$
for a  1-morphism or a 2-morphism $F$ in $\widehat{\sfC_h}$.

\begin{proposition}
\label{prp:H-factor}
    The $G$-invariant functor  
    \[
\udl{\rmH((P\down, \phi\down))} \colon
\frac{\rmH(\Mod \CC)}{\ang{\CU_{\CC}}}
\to \frac{\rmH(\Mod(\CC/G))}{\ang{\CU_{\CC/G}}},
\]     
    restricts to a $G$-precovering 
    $$\udl{\rmH(P\down, \ph\down)}:\frac{\rmH(\mmod \CC)}{\ang{\CU_{\mmod \CC}}}\to \frac{\rmH(\mmod (\CC/G))}{\ang{\CU_{\mmod \CC/G}}}.$$    
\end{proposition}

\begin{proof}
  Consider 2-subcategory $\sfC_h$, 2-functor $U_h$ and strict 2-natural transformation $\pi_h.$ It is plain that the conditions (0), (i) of Theorem \ref{thm:general-case} are  satisfied. Since $\pi_{\CM}$ is full (as stated in  Example \ref{exm:ess-surj-repbl}), the condition (ii) of the theorem is also satisfied. Therefore  the primary conditions required to apply Theorem \ref{thm:general-case} are fulfilled. Take the $G$-invariant functor $$\rmH(P\down, \phi\down) \colon (\rmH(\Mod \CC),  \rmH(A), \CU_{\CC}) \to (\rmH(\Mod (\CC/G)), \CU_{\CC/G}),$$ 
in $\sfC_h.$ Based on the previous observation before the proposition, it satisfies the required condition to use the theorem. Moreover, applying Proposition \ref{Prop-pushdown-adj}, we can conclude that it induces a $G$-precovering 
\begin{equation}
\label{eq: G-precoverin.mor}
\rmH(P\down, \phi\down) \colon
(\rmH(\mmod \CC),  \rmH(A), \cpt(\CU_{\CC}))
\to (\rmH(\mmod (\CC/G)), \cpt(\CU_{\CC/G})),
\end{equation}
by the adjunction $\rmH(\theta)$. 
Note that by Lemma \ref{lem.comp-morphism}:
\[
\begin{aligned}
\cpt(\rmH(\Mod \CC), \CU_{\CC})&=(\rmH(\mmod \CC), \CU_{\mmod \CC})\\
\cpt(\rmH(\Mod (\CC/G)), \CU_{\CC/G})&=(\rmH(\mmod (\CC/G)), \CU_{\mmod (\CC/G)}).
\end{aligned}
\]
Finally, we apply the 2-functor $U_h$ to $ \rmH(P\down, \phi\down)$ in \eqref{eq: G-precoverin.mor} to give  the desired $G$-precovering.   
\end{proof}

Based on the observation mentioned before, we can replace  the $G$-invariant functor $(P\down, \ph\down) \colon \Mod\CC \to \Mod(\CC/G)$ by the $G$-invariant functor $(P\down|, \ph\down|)\colon {\rm  Add}\mbox{-}\CK\rt {\rm Add}\mbox{-}\CK'$ (as given in Example \ref{ex-G-stable}),
which proves  the  first item of the following proposition.

\begin{proposition}\label{factorGpre}
Let $\CK$ and $\CK'$ be $G$-stable subcategories of $\mmod \CC$ and $\mmod (\CC/G)$. Assume $P\down$ sends $\CK$ to $\CK'$ and $P\up$ sends $\CK'$ to ${\rm Add}\mbox{-}\CK,$ e.g., $\CK'=P\down(\CK)$.  Then the following hold.
\begin{enumerate}
\item
The $G$-invariant functor  
    \[
\udl{\rmH((P\down|, \phi\down|))} \colon
\frac{\rmH(\Add\CK)}{\ang{\CU_{\CK}}}
\to \frac{\rmH(\Add\CK')}{\ang{\CU_{\CK'}}}.
\] 
restricts to a $G$-precovering
$$\udl{\rmH(P\down|, \phi\down |)} \colon 
\frac{\rm{H}(\CK)}{\ang{\CU_{\CK}}} \to
\frac{\rm{H}(\CK')}{\ang{\CU_{\CK'}}}$$
\item
Assume that $\CC$ is a locally bounded category,
$G$ is a torsion-free group acting freely on $\CC$,
and both $\CK$ and $\CK'$ are closed under direct sums
and direct summands.
Then
\begin{enumerate}			
\item
If $\udl{f}$ is an indecomposable object of
$\frac{\rm{H}(\CK)}{\ang{\CU_{\CK}}}$, then so is $\udl{\rmH(P\down)}(\udl{f})$.
\item
Let $\udl{f_1}$ and $\udl{f_2}$ be indecomposable objects of $\frac{\rm{H}(\CK)}{\ang{\CU_{\CK}}}$.
If
$\udl{\rmH(P\down)}(\udl{f_1})\cong \udl{\rmH(P\down)}(\udl{f_2})$,
then $\udl{f_1} \cong {}^a\!\udl{f_2}$ for some $a \in G$.
\end{enumerate}
\end{enumerate}
\end{proposition}

\begin{proof}
The proof of $(1)$ has been previously explained before the proposition.

(2)(a)
Let $\udl{f}$ be an indecomposable object in  $\frac{\rmH(\CK)}{\ang{\CU_{\CK}}}$.
Then because $\rmH(\CK)$ is a Krull--Schmidt category,
we can write $f = \Ds_{i=1}^n f_i$ for some
indecomposable objects $f_i$ in $\rmH(\CK),\ (i=1,\dots, n)$, and we have
$\udl{f} = \Ds_{i=1}^n \udl{f_i}$.
Since $\udl{f}$ is indecomposble, 
there is a unique $i$ such that $\udl{f_i} \ne 0$.
Take $f':= f_i$ and $g:= \Ds_{j \ne i} f_j$.
Then
$f=f'\oplus g$ in $\rm{H}(\CK)$,  where $f'$ is indecomposable in $\rm{H}(\CK)$ 
with $\udl{f'} \ne 0$, and $\udl{g} = 0$.
Hence, we have
\begin{equation}
\label{eq:HPf}
\rmH(P\down)(f)= \rmH(P\down)(f')\oplus \rmH(P\down)(g).
\end{equation}
According to Proposition \ref{Proposition3last}, $\rmH(P\down)(f')$ is indecomposable in $\rm{H}(\CK')$, and by
\eqref{eq:HPf} we have an isomorphism $\udl{\rmH(P\down)}(\udl{f}) \cong \udl{\rmH(P\down)}(\udl{f'})$ in the factor category $\frac{\rmH(\CK')}{\ang{\CU_{\CK'}}}$.
Hence it is enough to show that $\udl{\rmH(P\down)}(\udl{f'})$ is indecomposable.
We have a canonical surjective algebra homomorphsim
\[
\rmH(\CK')(\rmH(P\down)(f'),\rmH(P\down)(f'))
\to
\frac{\rmH(\CK')}{\ang{\CU_{\CK'}}}(\udl{\rmH(P\down)}(\udl{f'}),\udl{\rmH(P\down)}(\udl{f'})),
\]
where $\rmH(\CK')(\rmH(P\down)(f'),\rmH(P\down)(f'))$ is
a local algbebra because $\rmH(P\down)(f')$ is
indecompoable.
If we show that 
$\udl{\rmH(P\down)}(\udl{f'}) \ne 0$,
then it has a local endomorphism algebra, and hence
it turns out to be indecomposable, and the proof of (a) is completed.
Assume that $\udl{\rmH(P\down)}(\udl{f'}) = 0$.
Then $\id_{\rmH(P\down)(f')} \in \ang{\CU_{\CK'}}$.
Thus we have a factorization
\[
\begin{tikzcd}
\Nname{1}\rmH(P\down)(f') && \Nname{2}\rmH(P\down)(f')\\
&\Nname{3}(\Ds_{i=1}^m(X_i\to 0)) \ds (\Ds_{j=1}^n(X_j \ya{1}X_j))
\Ar{1}{2}{"\idty_{\rmH(P\down)(f')}"}
\Ar{1}{3}{"u" '}
\Ar{3}{2}{"{v}"}
\end{tikzcd}
\]
We apply $\rmH(P\up)$ to this diagram to have
\[
\begin{tikzcd}
\Nname{a}f' && \Nname{b}f'\\
\Nname{1}\Ds_{a\in G}{}^af' && \Nname{2}\Ds_{a\in G}{}^af'\\
&\Nname{3}(\Ds_{i=1}^m P\up(X_i)\to 0) \ds (\Ds_{j=1}^nP\up(X_j) \ya{1}P\up(X_j)))
\Ar{1}{2}{"\idty_{\Ds_{a\in G}{}^af'}"}
\Ar{1}{3}{"\rmH(P\up)(u)" '}
\Ar{3}{2}{"{\rmH(P\up)(v)}"}
\Ar{a}{b}{"\idty_{f'}"}
\Ar{a}{1}{"\si_{1_G}" '}
\Ar{2}{b}{"\pi_{1_G}" '}
\end{tikzcd}
\]
Therefore, $\idty_{f'}$ factors through $\ang{\CU_{\CK}}$,
hence $\udl{f'} = 0$, a contradiction.



(2)(b) Assume $\udl{f_1}$ and  $\udl{f_2}$ are indecomposable objects in  $\frac{\rm{H}(\CM)}{\ang{\CU_{\CM}}}$  such that  we have an isomorphism $\udl{\rmH(P\down)}(\udl{f_1})\cong \udl{\rmH(P\down)}(\udl{f_2})$ in $\frac{\rm{H}(\CK')}{\ang{\CU_{\CK'}}}$.  Hence,  we have $ \rmH(P\down)(f_1)\oplus {\rm h}\cong \rmH(P\down)(f_2)\oplus {\rm  g}$ for some objects ${\rm h}$ and ${\rm g}$ in $\rm{add}\mbox{-}\CU_{\CK'}$. We  can have  decompositions $f_1=f'_1\oplus m$ and $f_2=f'_2\oplus n$ in $\rm{H}(\CM)$, where $f'_1$ and $f'_2$ are indecomposable, and $m, n \in \rm{add}\mbox{-}\CU_{\CK}$. Then by applying the functor $\rmH (P\down)$ on the decompositions,  we obtain  the following isomorphisms:
\[
\begin{aligned}
 \rmH(P\down)(f'_1)\oplus \rmH(P\down)(m)\oplus {\rm h} &=\rmH(P\down)(f_1)\oplus {\rm h} \\ &\cong   \rmH(P\down)(f_2)\oplus {\rm g}  \\ &= \rmH(P\down)(f'_2)\oplus \rmH(P\down)(n)\oplus {\rm g}
\end{aligned}
\]
Since  $\rmH(P\down)(f'_1)$ and $\rmH(P\down)(f'_2)$ are only direct summands do  not belong to $\rm{add}\mbox{-}\CK'$, by the Krull-Schmidt property, we have $\rmH(P\down)(f'_1)\cong \rmH(P\down)(f'_2)$. It follows from Proposition \ref{Proposition3last} that $f'_1\cong {}^g (f'_2)$ for some $g \in G$. Therefore,
$$f_1\oplus {}^gn= f_1'\oplus m\oplus{}^g n\cong 
{}^g (f'_2) \oplus{}^g n\oplus m={}^gf_2\oplus m.$$
Hence, we obtain  the isomorphism
$\udl{f_1}\cong \udl{{}^gf_2} = {}^g\udl{f_2}$ in $\frac{\rm{H}(\CK)}{\ang{\CU_{\CK}}}.$
\end{proof}

\subsection{$G$-precoverings for the categories of finitely presented functors}


In this subsection, we will keep the notations introduced in Proposition \ref{factorGpre} unless otherwise stated.

It is worth noticing that if  both $\CK$ and $\CK'$ have weak kernels, and that $\CC$ is a locally bounded $k$-category, then 
we have $\mmod \CC=\CF(\CC)$, $\mmod \CK=\CF(\CK)$,
and
$\mmod \CC/G=\CF(\CC/G)$, $\mmod \CK'=\CF(\CK')$.  For instance, if we assume that $\CK':=P\down(\CK)$ is contravariantly finite in $\mmod \CC/G$ (and hence has weak kernels), then according to  Lemma \cite[Lemma 3.7]{DI} the subcategory $\CK$ is also contravarintly finite in $\mmod \CC$. As a result, $\CK$ inherits the property of having weak kernels.

The category of finitely presented functors from a category $\CM$
to $\Mod k$
can be expressed as a factor category of the morphism category of $\CM$
as follows.

\begin{construction}
\label{const:Theta}
Let $\CM$ be a linear category.
We define a functor
$\Th_{\CM} \colon \rmH(\CM) \to \CF(\CM)$ as follows.

{\bf On objects.}
By Applying the Yoneda functor, we can define
the correspondence between objects as follows:
\[
\Th_{\CM} \colon \rmH(\CM) \to \CF(\CM),
\quad \vertmap{X_1}{f}{X_2}\mapsto
\Th_{\CM}(f):=\Coker(\CM(\blank, X_1)\ya{\CM(\blank, f)}
\CM(\blank, X_2)).
\]
	
{\bf On morphisms.}
Let $\vertmap{X_1}{f}{X_2} \longto{h_1}{h_2} \vertmap{X_1'}{f'}{X_2'}$ be a morphism in $\rmH(\CM)$.
Define $\Th_{\CM}(h_1, h_2)$ as the unique morphism $\sigma$ that makes the following diagram commutative:
\[
\xymatrix{
\CM(\blank, X_1) \ar[r]^{\CM(\blank, f)} \ar[d]_{\CM(\blank, h_1)} & \CM(\blank, X_2) \ar[r] \ar[d]^{\CM(\blank, h_2)} & \Theta_{\CM}(f)  \ar[r] \ar[d]^{\sigma} & 0 \\
\CM(\blank, X'_1)  \ar[r]^{\CM(\blank,f')} & \CM(\blank, X'_2) \ar[r] & \Theta_{\CM}(f')  \ar[r]  & 0.}
\]
Since $\Th_{\CM}({\CU_{\CM}})=0$, we have the induced functor $\overline{\Theta}_{\CM}:\frac{\rm{H}(\CM)}{\ang{\CU_{\CM}}}\rt \CF(\CM)$ satisfying the following strictly commutative diagram
\[
\xymatrix {\rm{H}(\CM)\ar[d]_{\Theta_{\CM}} \ar[r]^{\pi_{\CC}} & \frac{\rm{H}(\CM)}{\ang{\CU_{\CM}}}\ar[ld]^{\overline{\Theta}_{\CM}} \\ \CF(\CM)   & }.
\]
\end{construction}

\begin{lemma}\label{Gequivalence}
Assume that $\calC$ is a locally bounded category, and that
both $\CK$ and $\CK'$ have weak kernels.
Then the following statements hold. 
\begin{itemize}
\item[$(1)$] The functors $\overline{\Theta}_{\CK'}$ and $\overline{\Theta}_{\CK}$ are equivalences.
\item[$(2)$] The functor $\overline{\Theta}_{\CK}$ is a $G$-equivariant equivalence.
\end{itemize}
\end{lemma}
\begin{proof}
$(1)$ We refer to \cite[Theorem 4.2]{HE} for a proof.

(2)
We first show that $\Theta_{\CK}$ is $G$-equivariant.
Let $(X\st{f}\rt Y)$ 
be an object in $\rmH(\CK)$ and $g \in G$.
By definitions, we have the following exact sequence:
$$
\CK(\blank, {}^gX) \ya{\CK(\blank, {}^g f)} \CK(\blank, {}^gY)\rt \Theta_{\CK}({}^g f) \rt 0.
$$
For any $Z \in \CK$, in view of the above sequence, we have the following commutative diagram:
\[
\xymatrix@C=50pt{
\CK(Z,{}^gX_1) \ar[r]^{\CK(Z, {}^gf)} \ar[d]^{\simeq} & \CK(Z, {}^gX_2) \ar[r] \ar[d]^{\simeq} & \Theta_{\CK}({}^g f) (Z) \ar[r] \ar[d]^{\eta} & 0 \\
\CK({}^{g^{\blank1}}Z, X_1)  \ar[r]^{\CK({}^{g^{-1}}Z,f)} & \CK({}^{g^{-1}}Z, X_2) \ar[r] & \Theta_{\CK}(f)({}^{g\inv}Z)  \ar[r]  & 0.}
\]
Since $\CK$ is a $G$-stable subcategory of the $G$-category $\mmod\CC$,
$\CK$ has the induced $G$-action,
which induces the $G$-action on $\CF(\CK)$
by $({}^g M)(Z):= M({}^{g\inv}Z)$ for all
$M \in \CF(\CK)$, $g \in G$ and $Z \in \CK$.
Hence we have ${}^g\Theta_{\CK}(f)(Z)=\Theta_{\CK}(f)({}^{g^{-1}}Z)$.
The above diagram yields a natural isomorphism $\eta \colon \Th_{\CK}({}^gf)\cong {}^g\Th_{\CK}(f)$, as we wanted to show.
Now since ${}^g\CU_{\CM}=\CU_{({}^g{\CM})}$,
we obtain the desired isomorphism $\overline{\Theta}_{\CK}({}^gf)\cong {}^g\overline{\Theta}_{\CK}(f)$.
\end{proof}

Define $\CF(P\down|):\CF(\CK)\rt \CF(\CK')$ to be the composite of the following functors:
$$
\CF(P\down|):\CF(\CK) \xrightarrow{(\ovl{\Theta}_{\CK})^{-}} \frac{\rmH(\CK)}{\ang{\CU_{\CK}}}
\xrightarrow{\udl{\rmH(P\down|)}} \frac{\rm{H}(\CK')}{\ang{\CU_{\CK'}}}
\xrightarrow{\ovl{\Th}_{\CK'}} \CF(\CK'),
$$
where we have chosen a quasi-inverse $(\ovl{\Theta}_{\CK})^{-}$ of $\ovl{\Theta}_{\CK}$.
Here is the commutative diagram that illustrates the relationship between these functors (up to a natural isomorphism in the middle square):
\[\begin{tikzcd}
	&& {\frac{\rmH(\CK)}{\ang{\CU_{\CK}}}} \\
	{} & {\rmH(\CK)} && {\CF(\CK)} & {} \\
	& \rm{H}(\CK') && \CF(\CK') \\
	&& {\frac{\rmH(\CK')}{\ang{\CU_{\CK'}}}}
	\arrow["{\pi_{\CK}}", from=2-2, to=1-3]
	\arrow["{\overline{\Theta}_{\CK}}", from=1-3, to=2-4]
	\arrow["{\Theta_{\CK}}", from=2-2, to=2-4]
	\arrow["{\rmH(P\down|)}"', from=2-2,  to=3-2]
	\arrow["{\pi_{\CK'}}"', from=3-2, to=4-3]
	\arrow["{\Theta_{\CK'}}", from=3-2, to=3-4]
	\arrow["{\CF(P\down|)}", from=2-4, to=3-4]
	\arrow["{\overline{\Theta}}_{\CK'}"', from=4-3, to=3-4]
\end{tikzcd}.
\]

The following is a direct consequence of Lemma \ref{Gequivalence} and Proposition \ref{factorGpre}:

\begin{proposition}
\label{G-precove-functor-catgeory}
Under the above notations, we have the following.
\begin{enumerate}
\item
The functor $\CF(P\down|):\CF(\CK)\rt \CF(\CK')$ is a $G$-precovering.
\item
Assume that $\CC$ is a locally bounded category,
$G$ is a torsion-free group acting freely on $\CC$,
and both $\CK$ and $\CK'$ are closed under direct sums
and direct summands.
Then
\begin{enumerate}			
\item
If $T$ is indecomposable in $\CF(\CK)$, then  so is $\CF(P\down|)(T)$.
\item
If objects $T_1$ and $T_2$ are indecomposable in $\CF(\CK)$, then $\CF(P\down|)(T_1)\cong \CF(P\down|)(T_2)$ implies $T_1\cong {}^gT_2$, for some $g \in G.$
\end{enumerate}
\end{enumerate}
\end{proposition}

As is mentioned in the introduction,
a functor $\Phi:\CF(\mmod \CC)\rt \CF(\mmod (\CC/G))$ is defined in  \cite[Section 5]{P} and is proved to be $G$-precovering
in \cite[Theorem 5.5]{P}.
In fact, we have a natural isomorphism $\CF(P\down)\cong \Phi$
by \cite[Proposition 5.3]{P}.

\begin{lemma}
\label{lem:cpt-fp}
Let $\CA$  be a category that has small coproducts. Then the following statement hold.  
\begin{enumerate}
\item
$\cpt(\frac{\rmH(\CA)}{\ang{\CU_{\CA}}})=\frac{\rmH(\CA)}{\ang{\CU_{\CA}}}$, and
\item
$\cpt(\CF(\CA))=\CF(\CA).$
\end{enumerate}
Moreover, in this case, $\CF(\CA)$ has small coproducts.
\end{lemma}

\begin{proof}
From Construction \ref{const:Theta}, we obtain the equivalence  $\frac{\rmH(\CA)}{\ang{\CU_{\CA}}}\simeq \CF (\CA)$. Therefore, the statements (1) and (2) are equivalent, and it is enough to show only the statement (2).
The latter  is clear as every finitely presented functor is compact.

To show that $\CF(\CA)$ has small coproducts, let $(F_i)_{i \in I}$ be a small family of objects in $\CF(\CA).$ For each $i \in I$, take a projective presentation $\CA(-, X_i)\st{\CA(-, f_i)}\to \CA(-, Y_i)\to F_i\rt 0$. Consider the morphism $(\oplus_{i \in I}X_i\st{\oplus f_i}\to \oplus_{i \in I} Y_i)$ in $\CA$. Let $F=\Theta_{\CA}(\oplus_{i \in I}f_i)$. It is not difficult to prove that $(F, \Th_{\CA}(\delta_i):F_i\to F)$ is a coproduct of the family $(F_i)_{i \in I}$, where $\delta_i:(X_i\st{f_i}\rt Y_i)\rt (\oplus_{i \in I}X_i\st{\oplus f_i}\to \oplus_{i \in I} Y_i)$ is the canonical injection.
\end{proof}

\begin{remark}
Later we will use the lemma above for the cases that
$\CA = \Mod\CC, \linebreak[3] \Add\CK, \Add\CK'$.
\end{remark}

By using the same argument as in Lemma \ref{Gequivalence}, we observe that the functors $\Theta_{\Mod \CC}, \Theta_{\Mod(\CC/G)}, \Theta_{{\rm Add}\mbox{-}\CK}$ and $\Theta_{{\rm Add}\mbox{-}\CK'}$ induce the equivalences
in the following diagrams:
\[
\begin{tikzcd}[column sep=50pt]
\Nname{HK} \frac{\rm{H}(\Mod \CC)}{\ang{\CU_{\Mod \CC}}} & \Nname{fK}\CF(\Mod \CC)\\
\Nname{HK'}\frac{\rm{H}(\Mod(\CC/G))}{\ang{\CU_{\Mod(\CC/G}}} & \Nname{fK'}\CF(\Mod(\CC/G)
\Ar{HK}{fK}{"\ovl{\Th}_{\Mod\CC}"}
\Ar{HK'}{fK'}{"\ovl{\Th}_{\Mod(\CC/G)}" '}
\Ar{HK}{HK'}{"\udl{\rmH(P\down)}" '}
\Ar{fK}{fK'}{"\mathrm{Fp}(P\down)", dashed}
\end{tikzcd},
\quad
\begin{tikzcd}[column sep=50pt]
\Nname{HK} \frac{\rm{H}(\Add\CK)}{\ang{\CU_{\Add\CK}}} & \Nname{fK}\CF(\Add\CK)\\
\Nname{HK'}\frac{\rm{H}(\Add\CK')}{\ang{\CU_{\Add\CK'}}} & \Nname{fK'}\CF(\Add\CK')
\Ar{HK}{fK}{"\ovl{\Th}_{\Add\CK}"}
\Ar{HK'}{fK'}{"\ovl{\Th}_{\Add\CK'}" '}
\Ar{HK}{HK'}{"\udl{\rmH(P\down|)}" '}
\Ar{fK}{fK'}{"\mathrm{Fp}(P\down|)", dashed}
\end{tikzcd}.
\]
By choosing quasi-inverses $(\ovl{\Th}_{\Mod\CC})^-$ 
and $(\ovl{\Th}_{\Add\CK})^-$
of $\ovl{\Th}_{\Mod\CC}$ and $\ovl{\Th}_{\Add\CK}$,
we can define $\mathrm{Fp}(P\down)$ and $\mathrm{Fp}(P\down|)$
(dashed arrows) in the diagrams above.
Similarly, by choosing quasi-inverses $(\ovl{\Th}_{\Mod(\CC/G})^-$ 
and $(\ovl{\Th}_{\Add\CK'})^-$
of $\ovl{\Th}_{\Mod(\CC/G)}$ and $\ovl{\Th}_{\Add\CK'}$,
we define $\mathrm{Fp}(P\up)$ and $\mathrm{Fp}(P\up|)$
by the following strictly commutative diagrams:
\[
\begin{tikzcd}[column sep=60pt]
\Nname{HK} \frac{\rm{H}(\Mod \CC)}{\ang{\CU_{\Mod \CC}}} & \Nname{fK}\CF(\Mod \CC)\\
\Nname{HK'}\frac{\rm{H}(\Mod(\CC/G))}{\ang{\CU_{\Mod(\CC/G}}} & \Nname{fK'}\CF(\Mod(\CC/G)
\Ar{HK}{fK}{"\ovl{\Th}_{\Mod\CC}"}
\Ar{fK'}{HK'}{"(\ovl{\Th}_{\Mod(\CC/G)})^-"}
\Ar{HK'}{HK}{"\udl{\rmH(P\up)}"}
\Ar{fK'}{fK}{"\mathrm{Fp}(P\up)" ', dashed}
\end{tikzcd},
\quad
\begin{tikzcd}[column sep=50pt]
\Nname{HK} \frac{\rm{H}(\Add\CK)}{\ang{\CU_{\Add\CK}}} & \Nname{fK}\CF(\Add\CK)\\
\Nname{HK'}\frac{\rm{H}(\Add\CK')}{\ang{\CU_{\Add\CK'}}} & \Nname{fK'}\CF(\Add\CK')
\Ar{HK}{fK}{"\ovl{\Th}_{\Add\CK}"}
\Ar{fK'}{HK'}{"(\ovl{\Th}_{\Add\CK'})^-"}
\Ar{HK'}{HK}{"\udl{\rmH(P\up|)}"}
\Ar{fK'}{fK}{"\mathrm{Fp}(P\up|)" ', dashed}
\end{tikzcd}.
\]

\begin{theorem}
\label{thm:G-Fp-precovering}
The following statements hold. 
\begin{enumerate}
\item
The $G$-invariant functor  ${\rm Fp}(P\down)$ is  a $G$-precovering,
and has the right adjoint ${\rm Fp}(P\up)$
with the adjucntion $\Theta \colon {\rm Fp}(P\down) \adj {\rm Fp}(P\up)$.

Moreover, there is the following commutative diagram
\begin{equation}
\begin{tikzcd}
\CF(\mmod \CC) & \CF(\Mod \CC)\\
\CF(\mmod (\CC/G))& \CF(\Mod (\CC/G))
\Ar{1-1}{1-2}{"", }
\Ar{1-2}{2-2}{"{\rm Fp}(P\down)", ""'}
\Ar{1-1}{2-1}{"{{\rm fp}(P\down)}"'}
\Ar{2-1}{2-2}{"",""'}
\end{tikzcd}
\end{equation}
where the horizontal functors are embeddings.
\item 
The $G$-invariant functor  ${\rm Fp}(P\down |)$ is  a $G$-precovering,
and has the right adjoint ${\rm Fp}(P\up|)$
with the adjucntion $\Theta \colon {\rm Fp}(P\down|) \adj {\rm Fp}(P\up|)$. 

Moreover, there is the following commutative diagram
\begin{equation}
\begin{tikzcd}
\CF(\CK) & \CF({\rm Add}\mbox{-}\CK)\\
\CF(\CK')& \CF({\rm Add}\mbox{-}\CK')
\Ar{1-1}{1-2}{"", }
\Ar{1-2}{2-2}{"{\rm Fp}(P\down|)", ""'}
\Ar{1-1}{2-1}{"{{\rm fp}(P\down|)}"'}
\Ar{2-1}{2-2}{"",""'}
\end{tikzcd}
\end{equation}
where the horizontal functors are embeddings.
\end{enumerate}
\end{theorem}
 
\begin{proof}
We only prove (1) since proof of (2) is similar.
The statement (1) is immediate from Theorem \ref{thm:general-case}(b) and Lemma \ref{lem:cpt-fp} because it is easy to verify that
the adjoint pair $({\rm Fp}(P\down), {\rm Fp}(P\up)))$ satisfies the conditions {\rm (PA1)}, {\rm (PA2)} and {\rm (PA3)} in Definition \ref{dfn:precov-adj}, i.e., 
\begin{enumerate}
\item[{\rm (PA1)}]
There exists a natural isomorphism $t \colon \Ds_{a\in G}\overline{A}_a \To {\rm Fp}(P\up){\rm Fp }(P\down)$.
\item[{\rm (PA2)}]
For any $x,y \in \CF(\Mod \CC)$, the following diagram commutes:
\begin{equation}
\begin{tikzcd}
\Ds_{a \in G} \calA(x, ay) & \calA(x, \Ds_{a \in G} ay)\\
\calB({\rm Fp}(P\down)(x), {\rm Fp}(P\down)(y)) & \calA(x, {\rm Fp}(P\up){\rm Fp}(P\down)y)
\Ar{1-1}{1-2}{"\nu", }
\Ar{1-2}{2-2}{"{\calA(x, t_y)}", "\rotatebox{90}{$\sim$}"'}
\Ar{1-1}{2-1}{"{({\rm Fp}(P\down))^{(2)}_{x,y}}"'}
\Ar{2-1}{2-2}{"\sim","\Theta_{x, {\rm Fp}(P\down)y}"'}
\end{tikzcd}
\end{equation}
where $\nu$ is the canonical morphism, $\CA:=\CF(\Mod \CC)$ and $\CB:=\CF(\Mod \CC/G).$

\item[{\rm (PA3)}]
${\rm Fp}(P\up)$ preserves small coproducts.
\end{enumerate}
%
\end{proof}


\begin{thebibliography}{9999}
\bibitem[As97]{As}
{H. Asashiba,} {A covering technique for derived equivalence,}  J. Algebra {\bf 191} (1997), 382-415. DOI 10.1006/jabr.1997.6906. MR MR1444505.

\bibitem[As11]{A}
{H. Asashiba,} {A generalization of Gabriel's Galois covering functors and derived equivalences,}  J. Algebra {\bf 334} (2011), 109-149. DOI 10.1016/j.jalgebra.2011.03.002. MR 2787656.

\bibitem[As22]{Asa-book}
H. Asashiba,
{Categories and representation theory--with a focus on 2-categorical covering theory}, Math. Surveys Monogr., {\bf 271}. American Mathematical Society, Providence, RI, 2022. xviii+240 pp. 
DOI 10.1090/surv/271.
MR 4486377.

\bibitem[As13]{A3}
{H. Asashiba,} 
{Gluing derived equivalences together,}
Adv. Math. {\bf 235} (2013), 134-160. DOI 10.1016/j.aim.2012.10.021. MR 3010054.

\bibitem[AHV18]{AHV}
H. Asashiba, R. Hafezi, and R. Vahed, Gorenstein versions of covering techniques for linear categories and their applications, J. Algebra {\bf 507} (2018), 320–361. DOI 10.1016/j.jalgebra.2018.04.017. MR 3807051.


\bibitem [Au66]{Au1} {M. Auslander,} {Coherent functors,} in Proc. Conf. Categorical Algebra (La Jolla, Calif., 1965), 189-231, Springer, New York, 1966. MR 0212070.


\bibitem[BL14]{BL}
{R. Bautista and S. Liu,} {Covering theory for linear categories with application to derived categories,} J. Algebra {\bf 406} (2014), 173-225. DOI 10.1016/j.jalgebra.2014.02.016. MR 3188335.


\bibitem[BG82]{BG}
{K. Bongartz and P. Gabriel,} {Covering spaces in representation-theory,} Invent. Math. {\bf 65} (1982) 331-378.
DOI 10.1007/BF01396624. MR 0643558.


\bibitem [DI20]{DI}
{E. Darp$\rm{\ddot{o}}$ and O. Iyama,} $d$-representation-finite selfinjective algebras. Adv. Math. {\bf 362} (2020), 106932, DOI 10.1016/j.aim.2019.106932. MR 4052555.


\bibitem[DLS86]{DLS}
{P. Dowbor, H. Lenzing and A. Skowro\'nski,} {Galois coverings by algebras of locally support-finite categories,} (Ottawa, Ont., 1984), 91-93, Lecture Notes in Math., 1177, Springer, Berlin, 1986. DOI 10.1007/BFb0075260. MR 0842461.


\bibitem[Ga62]{G1}
P. Gabriel, Des Cat$\acute{e}$gories Ab$\acute{e}$liennes, Bull. Soc. Math. France {\bf 90} (1962), 323-448. DOI 10.24033/bsmf.1583. MR 0232821.

\bibitem[Ga80]{G} {P. Gabriel,}
{The universal cover of a representation-finite algebra,} Representations of Algebras (Puebla, 1980), Lecture Notes in Math., 903, Springer, 1981, 68-105. MR 0654725.

	

\bibitem[HE23]{HE} {R. Hafezi and H. Eshraghi,} {From Morphism Categories to Functor Categories,} avaialabe on arXiv:2301.00534.


\bibitem[HM20]{HM} {R. Hafezi and E. Mahdavi,} {Covering theory, (mono)morphism categories and stable Auslander algebras,} available at arXiv:2011.08646.


\bibitem[HZha]{HZhao} {Y. Hu and T. Zhao,} {Morphisms determined by objects under Galois G-covering theory,}  J. Algebra {\bf 620} (2023), 225-256. DOI 10.1016/j.jalgebra.2022.11.034. MR 4531546.

 \bibitem[HZho]{HZou} {Y. Hu and P. Zhou,}{ Galois G-covering of quotients of linear categories,}  J. Pure Appl. Algebra {\bf 227} (2023), 107244. DOI 10.1016/j.jpaa.2022.107244. MR 4510795.


\bibitem[Le94]{Le94}
{Z. Leszczy\'nski,}
{On the representation type of tensor product algebras},
Fund. Math.
{\bf 144} (1994), 143-161. DOI 10.4064/fm-144-2-143-161. MR 1273693.



\bibitem[LS00]{LSk}
{Z. Leszczy\'nski and A. Skowro{\'n}ski,} {Tame triangular matrix algebras,} Colloq. Math. {\bf 86} (2000), 259-303. DOI 10.4064/cm-86-2-259-303. MR 1808681.




\bibitem[Pa19]{P}
{G. Pastuszak,} { On Krull-Gabriel dimension and Galois coverings,}  Adv. Math. {\bf 349} (2019), 959-991. DOI 10.1016/j.aim.2019.04.035. MR 3945585.

\bibitem[Pr09]{Pe}
{M. Prest,} {Purity, Spectra and Localization,} Encyclopedia Math. Appl., 121, Cambridge University Press, Cambridge, 2009. DOI 10.1017/CBO9781139644242. MR 2530988.

\bibitem[Po73]{Po73}
{N. Popescu,} Abelian categories with applications to rings and modules, London Math. Soc. Monogr., Vol. {\bf 3}. Academic Press, London { 1973}. MR 0340375.


\bibitem[Ri80]{Ri}
{C. Riedtmann,} {Algebren, Darstellungsk{\"o}cher, {\"U}berlagerungen and zur{\"u}ck,} Comment. Math. Helv. 55 (1980), 199-224. DOI 10.1007/BF02566682. MR 0576602.



\end{thebibliography}
\end{document}